\documentclass{article}
\makeatletter
\usepackage{eso-pic} % used by \AddToShipoutPicture
\RequirePackage{fancyhdr}
\RequirePackage{natbib}

\setcitestyle{authoryear,round,citesep={;},aysep={,},yysep={;}}

\newif\ificlrfinal
\iclrfinalfalse
\def\iclrfinalcopy{\iclrfinaltrue}
\font\iclrtenhv  = phvb at 8pt

\def\addcontentsline#1#2#3{}

\def\maketitle{\par
\begingroup
   \def\thefootnote{\fnsymbol{footnote}}
   \def\@makefnmark{\hbox to 0pt{$^{\@thefnmark}$\hss}} % for perfect author
   \long\def\@makefntext##1{\parindent 1em\noindent
                            \hbox to1.8em{\hss $\m@th ^{\@thefnmark}$}##1}
   \@maketitle \@thanks
\endgroup
\setcounter{footnote}{0}
\let\maketitle\relax \let\@maketitle\relax
\gdef\@thanks{}\gdef\@author{}\gdef\@title{}\let\thanks\relax}

\usepackage{fancyhdr}
\def\@maketitle{\vbox{\hsize\textwidth
{\LARGE\sc \@title\par}
\ificlrfinal
    \lhead{}
    \def\And{\end{tabular}\hfil\linebreak[0]\hfil
            \begin{tabular}[t]{l}\bf\rule{\z@}{24pt}\ignorespaces}%
  \def\AND{\end{tabular}\hfil\linebreak[4]\hfil
            \begin{tabular}[t]{l}\bf\rule{\z@}{24pt}\ignorespaces}%
    \begin{tabular}[t]{l}\bf\rule{\z@}{24pt}\@author\end{tabular}%
\else
       \lhead{Under review as a conference paper at ICLR 2027}
   \def\And{\end{tabular}\hfil\linebreak[0]\hfil
            \begin{tabular}[t]{l}\bf\rule{\z@}{24pt}\ignorespaces}%
  \def\AND{\end{tabular}\hfil\linebreak[4]\hfil
            \begin{tabular}[t]{l}\bf\rule{\z@}{24pt}\ignorespaces}%
    \begin{tabular}[t]{l}\bf\rule{\z@}{24pt}Anonymous authors\\Paper under double-blind review\end{tabular}%
\fi
\vskip 0.3in minus 0.1in}}

\renewenvironment{abstract}{\vskip.075in\centerline{\large\sc
Abstract}\vspace{0.5ex}\begin{quote}}{\par\end{quote}\vskip 1ex}

\def\section{\@startsection {section}{1}{\z@}{-2.0ex plus
    -0.5ex minus -.2ex}{1.5ex plus 0.3ex
minus0.2ex}{\large\sc\raggedright}}

\def\subsection{\@startsection{subsection}{2}{\z@}{-1.8ex plus
-0.5ex minus -.2ex}{0.8ex plus .2ex}{\normalsize\sc\raggedright}}
\def\subsubsection{\@startsection{subsubsection}{3}{\z@}{-1.5ex
plus      -0.5ex minus -.2ex}{0.5ex plus
.2ex}{\normalsize\sc\raggedright}}
\def\paragraph{\@startsection{paragraph}{4}{\z@}{1.5ex plus
0.5ex minus .2ex}{-1em}{\normalsize\bf}}
\def\subparagraph{\@startsection{subparagraph}{5}{\z@}{1.5ex plus
  0.5ex minus .2ex}{-1em}{\normalsize\sc}}

\skip\footins 9pt plus 4pt minus 2pt
\def\footnoterule{\kern-3pt \hrule width 12pc \kern 2.6pt }
\leftmargini\leftmargin \leftmarginii 2em
\def\@listi{\leftmargin\leftmargini}
\def\@listii{\leftmargin\leftmarginii
   \labelwidth\leftmarginii\advance\labelwidth-\labelsep
   \topsep 2pt plus 1pt minus 0.5pt
   \parsep 1pt plus 0.5pt minus 0.5pt
   \itemsep \parsep}
\def\@listiii{\leftmargin\leftmarginiii
    \labelwidth\leftmarginiii\advance\labelwidth-\labelsep
    \topsep 1pt plus 0.5pt minus 0.5pt
    \parsep \z@ \partopsep 0.5pt plus 0pt minus 0.5pt
    \itemsep \topsep}
\def\@listiv{\leftmargin\leftmarginiv
     \labelwidth\leftmarginiv\advance\labelwidth-\labelsep}
\def\@listv{\leftmargin\leftmarginv
     \labelwidth\leftmarginv\advance\labelwidth-\labelsep}
\def\@listvi{\leftmargin\leftmarginvi
     \labelwidth\leftmarginvi\advance\labelwidth-\labelsep}

\belowdisplayskip \abovedisplayskip
\def\normalsize{\@setsize\normalsize{11pt}\xpt\@xpt}
\def\small{\@setsize\small{10pt}\ixpt\@ixpt}
\def\footnotesize{\@setsize\footnotesize{10pt}\ixpt\@ixpt}
\def\scriptsize{\@setsize\scriptsize{8pt}\viipt\@viipt}
\def\tiny{\@setsize\tiny{7pt}\vipt\@vipt}
\def\large{\@setsize\large{14pt}\xiipt\@xiipt}
\def\Large{\@setsize\Large{16pt}\xivpt\@xivpt}
\def\LARGE{\@setsize\LARGE{20pt}\xviipt\@xviipt}
\def\huge{\@setsize\huge{23pt}\xxpt\@xxpt}
\def\Huge{\@setsize\Huge{28pt}\xxvpt\@xxvpt}

\makeatletter
\newbox\iclrrulerbox
\newcount\iclrrulercount
\newdimen\iclrruleroffset
\newdimen\cv@lineheight
\newdimen\cv@boxheight
\newbox\cv@tmpbox
\newcount\cv@refno
\newcount\cv@tot
\newcount\cv@tmpc@ \newcount\cv@tmpc
\def\fillzeros[#1]#2{\cv@tmpc@=#2\relax\ifnum\cv@tmpc@<0\cv@tmpc@=-\cv@tmpc@\fi
\cv@tmpc=1 %
\loop\ifnum\cv@tmpc@<10 \else \divide\cv@tmpc@ by 10 \advance\cv@tmpc by 1 \fi
   \ifnum\cv@tmpc@=10\relax\cv@tmpc@=11\relax\fi \ifnum\cv@tmpc@>10 \repeat
\ifnum#2<0\advance\cv@tmpc1\relax-\fi
\loop\ifnum\cv@tmpc<#1\relax0\advance\cv@tmpc1\relax\fi \ifnum\cv@tmpc<#1 \repeat
\cv@tmpc@=#2\relax\ifnum\cv@tmpc@<0\cv@tmpc@=-\cv@tmpc@\fi \relax\the\cv@tmpc@}%
\def\makevruler[#1][#2][#3][#4][#5]{\begingroup\offinterlineskip
\textheight=#5\vbadness=10000\vfuzz=120ex\overfullrule=0pt%
\global\setbox\iclrrulerbox=\vbox to \textheight{%
{\parskip=0pt\hfuzz=150em\cv@boxheight=\textheight
\cv@lineheight=#1\global\iclrrulercount=#2%
\cv@tot\cv@boxheight\divide\cv@tot\cv@lineheight\advance\cv@tot2%
\cv@refno1\vskip-\cv@lineheight\vskip1ex%
\loop\setbox\cv@tmpbox=\hbox to0cm{{\iclrtenhv\hfil\fillzeros[#4]\iclrrulercount}}%
\ht\cv@tmpbox\cv@lineheight\dp\cv@tmpbox0pt\box\cv@tmpbox\break
\advance\cv@refno1\global\advance\iclrrulercount#3\relax
\ifnum\cv@refno<\cv@tot\repeat}}\endgroup}%
\makeatother
\def\iclrruler#1{\makevruler[12pt][#1][1][3][0.993\textheight]\usebox{\iclrrulerbox}}
\AddToShipoutPicture{%
\ificlrfinal\else
\iclrruleroffset=\textheight
\advance\iclrruleroffset by -3.7pt
  \color[rgb]{.7,.7,.7}
  \AtTextUpperLeft{%
    \put(\LenToUnit{-35pt},\LenToUnit{-\iclrruleroffset}){%left ruler
      \iclrruler{\iclrrulercount}}
  }
\fi
}
\makeatother
\iclrfinalcopy

\usepackage{times}

\usepackage[T1]{fontenc}
\usepackage{amsmath,amssymb,amsthm,mathtools}
\usepackage{microtype,booktabs,tabularx,array,enumitem}
\usepackage{graphicx,placeins}
\usepackage[pass]{geometry}
\usepackage{pdflscape,xcolor,longtable}
\newif\ifshowrevisions
\showrevisionsfalse
\newcommand{\revisioncolor}{\ifshowrevisions\color{red}\fi}
\DeclareRobustCommand{\revised}[1]{{\revisioncolor#1}}
\newenvironment{revisedblock}{\par\revisioncolor}{\par}
\definecolor{pmblue}{HTML}{000000}
\definecolor{pmgray}{HTML}{000000}
\usepackage{hyperref}
\hypersetup{hidelinks,
  pdftitle={An Asymptotic One-Sample Surrogate for Non-Gaussian ERM with General Convex Regularization},
  pdfauthor={Tianjian Zhang, Cosme Louart, Malik Tiomoko, Zhenyu Liao}}
\allowdisplaybreaks[2]
\newcommand{\E}{\mathbb E}
\newcommand{\R}{\mathbb R}
\newcommand{\Pp}{\mathbb P}
\newcommand{\cE}{\mathcal E}
\newcommand{\cH}{\mathcal H}

\newcommand{\tr}{\operatorname{tr}}
\newcommand{\Cov}{\operatorname{Cov}}
\newcommand{\Var}{\operatorname{Var}}
\newcommand{\prox}{\operatorname{prox}}
\newcommand{\argmin}{\operatorname*{arg\,min}}
\newcommand{\op}{\mathrm{op}}
\newcommand{\hs}{\mathrm{hs}}
\newcommand{\dd}{\mathrm d}

\newcommand{\ind}{\mathbf 1}
\newcommand{\soft}{\operatorname{soft}}
\newcommand{\Law}{\mathcal L}
\newcommand{\Wass}{d_{\mathrm W}}
\newtheorem{assumption}{Assumption}
\newtheorem{theorem}{Theorem}
\newtheorem{proposition}[theorem]{Proposition}
\newtheorem{lemma}[theorem]{Lemma}
\newtheorem{corollary}[theorem]{Corollary}
\theoremstyle{definition}
\newtheorem{example}{Example}
\theoremstyle{remark}
\newtheorem{remark}{Remark}

\title{An Asymptotic One-Sample Surrogate \\ for Non-Gaussian ERM \\with General Convex Regularization}
\author{Tianjian Zhang\textsuperscript{1}\quad
Cosme Louart\textsuperscript{1}\quad
Malik Tiomoko\textsuperscript{2}\quad
Zhenyu Liao\textsuperscript{3}}
\begin{document}
\maketitle
\begingroup
\renewcommand{\thefootnote}{}
\footnotetext{\textsuperscript{1}School of Data Science, The Chinese University
of Hong Kong, Shenzhen, China.
\textsuperscript{2}Huawei Noah's Ark Lab, Huawei Technologies, Paris, France.
\textsuperscript{3}School of Electronic Information and Communications,
Huazhong University of Science \& Technology, China.\\
Correspondence to: Cosme Louart
\textless\href{mailto:cosmelouart@cuhk.edu.cn}{cosmelouart@cuhk.edu.cn}\textgreater.}
\endgroup
\begin{abstract}
We study convex empirical risk minimization (ERM) with non-Gaussian features and nonsmooth penalties through a surrogate driven by a single Gaussian vector. The surrogate preserves the penalty and calibrates its coefficients using the actual joint law of the feature--label model. 
We prove that this population calibration exists and is unique. 
\revised{We prove a Gaussian approximation for fresh prediction scores when the data satisfy a Poincaré inequality.} 

Together with regularity conditions on the losses, this approximation yields agreement between the ERM minimizer and its surrogate in normalized means and prediction variances, as well as formulas for prediction and classification risks. The result covers general uniformly strongly convex penalties and pure Lasso under projection anti-concentration. 
The core proof relies on a cavity Hamilton--Jacobi argument that compares the main statistics of the two minimizers without differentiating the penalty.
\end{abstract}

\section{Introduction}
\label{sec:intro}
We study high-dimensional convex empirical risk minimization:
\begin{equation}
\label{eq:original}
 \widehat W:= \widehat W((x_i, y_i)_{i\in[n]}) = \argmin_w\left\{\rho(w)+\sum_{i=1}^n
 L_{y_i}(x_i^\top w/\sqrt n)\right\}.
\end{equation}
The training pairs $(x_i,y_i)$ are independent, with $x_i\in\R^p$ and
$p/n\to c\in(0,\infty)$. Features $(x_i)_{i\in[n]}$ may be non-Gaussian with dependent
coordinates. The loss is smooth and convex; the finite convex penalty
$\rho:\R^p\to\R$ may be nonsmooth, including Lasso, which can set
coefficients exactly to zero.
Gaussian theory and its universality extensions give precise descriptions
under additional structure. They do not justify a general Gaussian
replacement with matching covariances: the joint feature--label law can
affect learning errors.

Inspired by earlier surrogate approaches
\citep{Loureiro2021,CelentanoMontanariWei2023}, we introduce a surrogate
designed to approximate the prediction behavior of $\widehat W$.
Write $m_x:=\E\left[x\right]$, $\bar x:=x-m_x$, and
$C_x:=\E\left[\bar x\bar x^\top\right]$ for the feature mean, centered feature, and
covariance. The second-moment matrix is
$\Sigma_x:=\E\left[xx^\top\right]=C_x+m_xm_x^\top$. Define
\begin{equation}
\label{eq:intro-prox}
 W_\ast(G):=\argmin_w\left\{\rho(w)+\frac{\nu_\ast}{2}w^\top \Sigma_xw
       +\beta_\ast G^\top w-\sqrt n\,\eta_\ast^\top w\right\},
 \qquad G\sim N(0,\Sigma_x).
\end{equation}
The coefficients $\nu_\ast,\beta_\ast>0$ and $\eta_\ast\in\R^p$ are three
\emph{calibration parameters} (we will further introduce $\alpha_\ast, \kappa_\ast\geq 0, \mu_\ast\in\mathbb R^p$). For the given penalty, they are calibrated
from the actual joint law of $(x,y)$ and the loss by solving the fixed-point
equations \eqref{eq:losscoeff}--\eqref{eq:fixedpoint}.
We call $W_\ast(G)$ a \emph{one-sample surrogate} because, once calibrated,
it replaces a minimization problem depending on $n$ training pairs by one
whose only random input is a single $p$-dimensional Gaussian vector $G$.
The law of $(x,y)$ remains encoded in $\Sigma_x$ and the calibration
parameters. This substantially simplifies the random optimization problem
while preserving the original penalty $\rho$.

\begin{revisedblock}
\section{Main contributions and comparison with prior work}
\label{sec:contributions}
\paragraph{A Hamilton--Jacobi equation for $\ell_1$ regularization under general assumptions.}
In one sentence, our main result extends the non-Gaussian ERM analysis of
\citet{Yaakoubi2026} to pure $\ell_1$ regularization. Their framework
requires a smooth, uniformly strongly convex penalty. The nonsmooth Lasso
penalty calls for a different proof strategy. We construct an interpolation
path connecting the original empirical problem to the one-sample surrogate,
while keeping the original penalty throughout. The normalized expected
minimum satisfies an asymptotic Hamilton--Jacobi equation
(Theorem~\ref{thm:main-interior-HJ}). A calibrated path then transports
prediction moments between the two problems without differentiating the
penalty. To our knowledge, this Hamilton--Jacobi interpolation between
non-Gaussian ERM and its one-sample surrogate, with the original nonsmooth
penalty preserved, is new.

\paragraph{Non-Gaussian features with dependent coordinates.}
Gaussian AMP and CGMT analyses already handle nonsmooth penalties with
independent design entries \citep{BayatiMontanari2012,Thrampoulidis2018};
\citet{HanShen2023} extend universality to independent non-Gaussian entries.
\citet{Vilucchio2025} also treat nonconvex objectives, using independent
Gaussian features, single-index responses, and separable penalties.
Correlated Gaussian covariates and mixtures are covered by
\citet{Loureiro2021,LoureiroMixtures2021}, and Gaussian-mixture features
with single-index targets by \citet{Pesce2023}.

Our advantage is in the feature assumptions: neither Gaussianity nor
coordinate independence nor a prescribed mixture representation is required.
We use convex objectives under the joint-data Poincar\'e and spectral
conditions stated below and the penalty and population conditions of
Theorem~\ref{thm:endpoint}. The uniformly strongly convex penalty may also
be nonsmooth and nonseparable.

Existing non-Gaussian results also impose specific model structures:
\citet{Bosch2023} treat random-feature least squares with nonsmooth
penalties, and \citet{LahirySur2024} treat Lasso and ridge in linear models
with block-dependent coordinates. Our theorem instead covers the stated
class of smooth convex losses, including logistic loss, under concentration
conditions on the feature--label law. It requires no random-feature
construction, block partition, or linear response model.
\citet{Mallory2025} go beyond our
setting by allowing dependence between observations, but focus on ridge
logistic regression. We retain independent training pairs while allowing
pure $\ell_1$ regularization and the broader penalty class described above.

\paragraph{A proved Gaussian approximation for fresh prediction scores.}
Theorem~\ref{thm:prediction-directions} gives qualitative uniform Gaussian
approximation along the complete calibrated path, including the unprobed
empirical endpoint, under the stated penalty/design and population
conditions. Its quantitative input is
Proposition~\ref{prop:prediction-positive-probes} in
Appendix~\ref{app:pred-proof}, which gives a $Cn^{-1/8}$ bound at every
interpolation time with fixed positive probe bounds under
Assumptions~\ref{ass:analytic}--\ref{ass:feature-concentration}. The independent test feature $x$ is held fixed;
the optimizer varies over the training data, selectors, and Gaussian probe.
The approximation error is then averaged
over $x$. This resembles Klartag's theorem on Gaussian marginals of
isotropic convex bodies in most directions \citep{Klartag2007}, with
uniform spherical averaging replaced here by averaging over a feature law
satisfying a Poincar\'e inequality, using the Stein kernel from
Lemma~\ref{lem:pred-proof-kernel} in Appendix~\ref{app:pred-proof}.
In \citet[Assumption~4, Theorem~6.1]{Yaakoubi2026}, Gaussian
approximation of estimator projections is an additional hypothesis for the
test-score formula. The hypothesis in
\citet[Assumption~5, equations~(11)--(12)]{MontanariSaeed2022} instead
concerns feature projections along deterministic coefficient vectors.
Broader Gaussian equivalence is proposed as Conjecture~1 in
\citet{Loureiro2021}; their Section~3.4 discusses its limits.
\citet{Dandi2023} extend universality to non-Gaussian mixtures under a
Gaussian projection hypothesis within each cluster. Their comparison
replaces the data by an equivalent Gaussian mixture.
\citet[Assumptions~4--7, Theorem~2.1]{Louart2026} proves a standardized CLT for nondegenerate deterministic projections, with a
quadratic penalty, a uniformly convex smooth loss, and further derivative
bounds. Our result permits nonsmooth penalties and averages the
approximation error over test features. Unlike the nondegenerate
coordinate projections covered by \citet{Louart2026}, fitted coordinates
in our setting need not be asymptotically Gaussian.
Proposition~\ref{prop:coordinate-nongaussian} in
Appendix~\ref{app:directional-counterexample} gives a pure-Lasso
counterexample with positive probes, even for Gaussian data.

\paragraph{A test-score law calibrated to the actual data distribution.}
Under the loss and penalty conditions stated below,
Theorem~\ref{thm:population-uniqueness} proves existence and uniqueness of
the calibration parameters $\alpha_\ast, \kappa_\ast, \mu_\ast,\nu_\ast,\beta_\ast,\eta_\ast$. Set
$\mu_\ast:=\E_G\left[W_\ast\right]/\sqrt n$ and
$\alpha_\ast^2:=\tr(\Sigma_xC_{W_\ast})/n$ with $C_{W_\ast}=\Cov(W_\ast)$.
Theorem~\ref{thm:endpoint} identifies these quantities with the normalized
mean and averaged prediction variance of $\widehat W$.
Corollaries~\ref{cor:conditional-prediction} and~\ref{cor:classification}
give the corresponding test-score law and prediction and classification
risks. For an independent test feature $x$ held fixed, the distribution of
$x^\top\widehat W/\sqrt n$ over training samples is approximated by that of
\[
 U_\ast:=x^\top\mu_\ast+\alpha_\ast Z,
 \qquad Z\sim N(0,1),\quad Z\perp(x,y).
\]
The approximation error tends to zero after averaging over $x$.
The signal $x^\top\mu_\ast$ retains the original feature--label law and may
be non-Gaussian. \citet{Yaakoubi2026} obtain a related test-score formula
for smooth, uniformly strongly convex penalties, and \citet{MaiLiao2025}
obtain one for ridge-regularized linear factor mixtures. Neither analysis
covers pure $\ell_1$ regularization. To our knowledge, the works discussed
above do not identify the mean, prediction variance, and test-score law
for pure Lasso under the non-Gaussian, dependent-feature assumptions used
here. Our formula preserves the predictive signal even when a Gaussian
replacement based on matching covariances gives different learning errors.
A finite Gaussian mixture can approximate this score law by discretizing
the distribution of $x^\top\mu_\ast$. Such an approximation can require
many components \citep{MaWuYang2025}, limiting descriptions with only a few
clusters. The Gaussian-mixture theory of
\citet[Assumption~(A4)]{LoureiroMixtures2021} treats a fixed finite number
of clusters as the sample size and dimension grow.
Figure~\ref{fig:diamonds-prediction} compares our score and classification
predictions with the closest competing approaches: the smooth-penalty
equations of \citet{Yaakoubi2026} and the Gaussian-mixture equations of
\citet{LoureiroMixtures2021}. The diamonds-based experiment in
Example~\ref{sec:penalties} illustrates why preserving both the exact
Lasso penalty and the original feature--label law matters for predicting
the observed behavior.

\paragraph{Organization.}
Section~\ref{sec:main-results} states the population calibration,
Gaussian score approximation, and prediction results, with numerical
prediction and classification comparisons in
Example~\ref{sec:penalties}. Example~\ref{sec:sparse-applications} gives
explicit Lasso and elastic-net formulas.
Section~\ref{sec:HJ} explains the Hamilton--Jacobi mechanism.
Appendix~\ref{app:pred-proof} proves
Proposition~\ref{prop:prediction-positive-probes} and the strongly convex
case of Theorem~\ref{thm:prediction-directions}.
Appendices~\ref{app:four} and~\ref{app:value-transport} use only that
proposition; Appendix~\ref{app:proof-main-theorem} completes
Theorem~\ref{thm:prediction-directions} for pure Lasso and proves
Theorem~\ref{thm:endpoint}.
\end{revisedblock}

\section{Setting and theoretical results}
\label{sec:main-results}
\subsection{Definition and well-posedness of the surrogate}
\label{sec:population}
Fix $n,p$ and a feature--label law $(x,y)\in\R^p\times\mathcal Y$,
where $\mathcal Y$ is the label space, with finite second moment
$\Sigma_x:=\E\left[xx^\top\right]$. We use $|\cdot|$ for the Euclidean
norm of a vector and the absolute value of a scalar. 
We use the Moreau-envelope notation of \citet{Thrampoulidis2018}.
For a convex function $\omega:\R^d\to\R$ and $\tau>0$, write
\begin{equation}
\label{eq:moreau}
 e_\omega(u;\tau):=\min_{z\in\R^d}
 \left\{\omega(z)+\frac{|u-z|^2}{2\tau}\right\}.
\end{equation}
The minimizing point is $\prox_{\tau\omega}(u)$, and
$e_\omega(u;0):=\omega(u)$ is the right-limit extension.
For the scalar losses, set $\prox_{0L_y}(u):=u$.
In \eqref{eq:losscoeff}--\eqref{eq:fixedpoint}, loss expectations use the
actual joint law of $(x,y)$ and an independent $Z\sim N(0,1)$; penalty
expectations use an independent $G\sim N(0,\Sigma_x)$.
For a vector $\mu\in\R^p$ and $\alpha,\kappa\ge0$, define the loss
coefficients by \eqref{eq:losscoeff}.
The population system is considered on the branch $\beta,\nu>0$:
\begin{equation}
\label{eq:losscoeff}
 \begin{gathered}
 U:=x^\top\mu+\alpha Z,\qquad \zeta:=\prox_{\kappa L_y}(U),\qquad
 \xi:=L_y'(\zeta),\\
 \beta^2:=\E\left[\xi^2\right],\qquad
 \nu:=\E\left[\frac{L_y''(\zeta)}{1+\kappa L_y''(\zeta)}\right],\qquad
 \eta:=\nu \Sigma_x\mu-\E\left[x\xi\right].
 \end{gathered}
\end{equation}
\begin{subequations}
\label{eq:fixedpoint}
\begin{gather}
\label{eq:effective}
 W(G):=\argmin_w\{\rho(w)+\nu w^\top \Sigma_xw/2+\beta G^\top w
                          -\sqrt n\,\eta^\top w\},\\
 \mu=\frac{\E_G\left[W(G)\right]}{\sqrt n},\qquad
 \alpha^2=\frac1n\tr\!\left(\Sigma_x C_W\right),\qquad
 \kappa=-\frac{\E_G\left[G^\top W(G)\right]}{\beta n}.
\end{gather}
\end{subequations}

Equations~\eqref{eq:losscoeff}--\eqref{eq:fixedpoint} contain auxiliary
Gaussian expectations and the population law, but no realization of the
training sample.

\begin{theorem}[Existence and uniqueness of the population solution]
\label{thm:population-uniqueness}
Fix $n,p$. Assume $\E\left[|x|^2\right]<\infty$ and
$\Sigma_x=\E\left[xx^\top\right]\succ0$. Suppose $\rho:\R^p\to\R$ is
convex with $0\in\partial\rho(0)$, and every $L_y$ is $C^2$ with
\[
 |L_y'(u)|\le L_1,\qquad 0<L_y''(u)\le L_2,\qquad u\in\R.
\]
Then \eqref{eq:losscoeff}--\eqref{eq:fixedpoint} has at most one finite
solution with $\alpha,\kappa\ge0$ and $\beta,\nu>0$; any such solution
has $\alpha,\kappa>0$.

If, in addition, the losses are nonnegative with
$\E\left[L_y(0)\right]<\infty$, and
\begin{equation}
\label{eq:population-existence-conditions}
 \operatorname*{arg\,min}_{w\in\R^p}\rho(w)=\{0\},
 \qquad \E\left[L_y'(0)^2\right]>0,
\end{equation}
then this solution exists.
\end{theorem}
Theorem~\ref{thm:population-uniqueness} is proved in
Appendix~\ref{app:population-stability}. Under the existence assumptions of
Theorem~\ref{thm:population-uniqueness},
denote the unique solution by
$(\mu_\ast,\alpha_\ast,\beta_\ast,\nu_\ast,\kappa_\ast,\eta_\ast)$ and its
optimizer by $W_\ast(G)$, as in \eqref{eq:intro-prox}.
The penalty and loss conditions in Assumption~\ref{ass:analytic} imply the
existence assumptions of Theorem~\ref{thm:population-uniqueness}. Positivity
at each finite size does not
imply the uniform positive-variance condition in Theorem~\ref{thm:endpoint}.
Lemma~\ref{lem:population-stationarity} in Appendix~\ref{app:population-stationarity}
gives an equivalent stationary-point formulation.

\begin{example}[Lasso and elastic regularization]
\label{sec:sparse-applications}
Our aim here is to give easy-to-use formulas for the surrogate and its
prediction moments. Take
\[
 \rho(w):=\lambda_1\|w\|_1+\frac{\lambda_2}{2}|w|^2,\qquad
 \|w\|_1:=\sum_{j=1}^p|w_j|,\qquad
 \lambda_1>0,\quad\lambda_2\ge0,
\]
where $|\cdot|$ is the Euclidean norm. If
$\Sigma_x=\operatorname{diag}(c_1,\ldots,c_p)$ in the penalty's
coordinates, put $\Gamma_{\ast,j}:=\nu_\ast c_j+\lambda_2$ and
$q_{\ast,j}:=\sqrt n\,\eta_{\ast,j}-\beta_\ast\sqrt{c_j}Z_j$, with
independent $Z_j\sim N(0,1)$. Then
\begin{equation}
\label{eq:soft}
 W_{\ast,j}=\frac{\soft(q_{\ast,j},\lambda_1)}{\Gamma_{\ast,j}},\qquad
 \soft(z,\lambda_1):=\operatorname{sign}(z)(|z|-\lambda_1)_+,
\end{equation}
where $a_+:=\max\{a,0\}$ and $\operatorname{sign}(0):=0$.
The moment equations reduce to one-dimensional Gaussian expectations:
\begin{equation}
\label{eq:diagonal-sparse-moments}
 \begin{gathered}
 \sqrt n\mu_{\ast,j}=\frac{\E\left[\soft(q_{\ast,j},\lambda_1)\right]}{\Gamma_{\ast,j}},\quad
 \alpha_\ast^2=\frac1n\sum_j\frac{c_j\Var(\soft(q_{\ast,j},\lambda_1))}{\Gamma_{\ast,j}^2},\\
 \kappa_\ast=\frac1n\sum_j\frac{c_j\Pp(|q_{\ast,j}|>\lambda_1)}{\Gamma_{\ast,j}}.
 \end{gathered}
\end{equation}
Proposition~\ref{prop:diagonal-sparse-formulas} proves these formulas.
The loss supplies curvature $\nu_\ast c_j>0$ even for pure Lasso
($\lambda_2=0$), which Theorem~\ref{thm:endpoint} covers under
\eqref{eq:lasso-strength-bounds} and
\eqref{eq:projection-anticoncentration}.
\end{example}

%  If $\inf_n\lambda_2>0$, the
% elastic penalty satisfies \eqref{eq:strongly-convex-penalty}, without
% anti-concentration or a Lasso threshold, also for size-dependent
% $\lambda_1$; a vanishing ridge does not meet this condition.

% Diagonal second moments suffice for separation of the surrogate and do
% not require independent features. For general $\Sigma_x$, the surrogate
% remains a coupled convex problem, with response given by
% Corollary~\ref{cor:correlated-sparse-response}. For pure Lasso, it yields
% $\nu_\ast\kappa_\ast=\E_G[|\operatorname{supp}W_\ast|]/n$.
% This identity concerns the surrogate; transfer to training support counts
% remains open. Appendix~\ref{app:sparse-specialization} gives the derivations.

\subsection{Asymptotic setting and assumptions}
\label{sec:setting}
For comparison with the trained estimator, consider a sequence with
\[
 n\longrightarrow\infty,\qquad p=p(n),\qquad
 \frac pn\longrightarrow c\in(0,\infty).
\]
For each $n$, the training pairs $(x_i,y_i)_{i=1}^n$ are independent copies
of $(x,y)$.
The population law, loss, penalty, and all derived quantities may depend
on $n$. We suppress this dependence, including for $\Sigma_x$, $\widehat W$,
and the population order parameters. Independence is required within each
training sample; no coupling of different sample sizes is assumed.
We assume
\begin{equation}
\label{eq:spectral-bounds}
 \underline c_x I_p\preceq \Sigma_x\preceq KI_p,
 \qquad \underline c_x>0.
\end{equation}
Also $p/n\le K$ for all sufficiently large $n$.
All constants in \eqref{eq:spectral-bounds}, the eventual bound $p/n\le K$,
and the loss, data-Poincar\'e, and penalty/design conditions in
Assumptions~\ref{ass:analytic}--\ref{ass:penalty-structure}, are uniform along
the sequence unless dependence on a fixed auxiliary parameter is stated
explicitly.
Write $\|A\|_{\hs}:=(\tr(A^\top A))^{1/2}$ for the
Hilbert--Schmidt norm and $\|\cdot\|_{L^q}$ for a moment norm.
Feature coordinates may be dependent and $\E\left[x\right]$ need not vanish.
We study the estimator $\widehat W$ in \eqref{eq:original};
the score normalization remains $x^\top w/\sqrt n$ after thinning or deletion.
The one-sample surrogate uses the same dimension and population law as the
training problem. Our conclusions compare these size-matched objects;
they do not require population coefficients to have dimension-independent
limits.

\begin{assumption}[Penalty and loss]
\label{ass:analytic}
The penalty $\rho:\R^p\to\R$ is convex,
with $\operatorname*{arg\,min}_{w\in\R^p}\rho(w)=\{0\}$.
The losses are nonnegative convex $C^4$ functions, with
\[
 \sup_{y,u}|L_y^{(j)}(u)|\le L_j\quad(j=1,2,3,4),
 \qquad 0\le L_y(0)\le L_0\quad\text{almost surely}.
\]
For every $M<\infty$, 
$ c_L(M):=\inf_{y,\,|u|\le M}L_y''(u)>0,$ and
 $\E\left[L_y'(0)^2\right]>0$ at every problem size. The constants $L_0,\ldots,L_4$ and $c_L(M)$, for each fixed $M$, are uniform along the
sequence.
\end{assumption}
Binary logistic loss satisfies the loss conditions in
Assumption~\ref{ass:analytic}, with
$L_y(0)=\log2$ and $\E\left[L_y'(0)^2\right]=1/4$;
squared loss has unbounded slope and needs additional moment estimates that we do not provide here for simplicity.
Theorem~\ref{thm:population-uniqueness} retains its broader finite-size
assumptions. Theorems~\ref{thm:population-uniqueness}
and~\ref{thm:main-interior-HJ} allow general convex penalties;
Theorem~\ref{thm:prediction-directions} and
Theorem~\ref{thm:endpoint} cover pure Lasso and uniformly strongly convex
penalties under their additional population conditions.

\begin{assumption}[Poincar\'e inequality for the data]
\label{ass:feature-concentration}
The joint law of $(x,y)$ satisfies, for a constant $C_{\rm P}<\infty$
independent of dimension, for all bounded test functions $F$ continuously
differentiable in their continuous arguments
\begin{equation}
\label{eq:data-poincare}
 \Var_{x,y}\!\left(F(x,y)\right)
 \le C_{\rm P}\E_{x,y}\left[
       |\nabla_xF(x,y)|^2+\Gamma_yF(x,y)\right].
\end{equation}
where, for Euclidean labels, $\Gamma_yF(x,y):=|\nabla_yF(x,y)|^2$.
For a fixed finite label space, use
$\Gamma_yF(x,y):=\frac12\sum_{y'\in\mathcal Y}
|F(x,y')-F(x,y)|^2$.
\end{assumption}
The finite-label definition uses differences, so binary labels are included
without differentiating them. Condition~\eqref{eq:data-poincare} concerns
the joint law; it requires neither independence of $x$ and $y$ nor centered
or independent feature coordinates. Applying it to functions of $x$ alone
gives, by Lemma~\ref{lem:qf-moments},
\begin{equation}
\label{eq:feature-quadratic-concentration}
 \E_x\left[\left|x^\top Ax-\tr(\Sigma_xA)\right|^2\right]
 \le K_Q^2\|A\|_{\hs}^2,\qquad K_Q^2:=4C_{\rm P}K,
\end{equation}
for every deterministic symmetric $A$.

Appendix~\ref{app:poincare-criteria} gives tensorization and
Lipschitz-transformation criteria for
Assumption~\ref{ass:feature-concentration}.
For the synthetic model in
Appendix~\ref{app:synthetic-prediction-experiments}, the fixed-size
Gaussian-mixture blocks have finite Poincar\'e constants. Tensorization
and the uniformly Lipschitz response give a joint bound independent of
$n,p$ by Lemma~\ref{lem:poincare-stability}.

\begin{assumption}[Penalty and design alternatives]
\label{ass:penalty-structure}
One of the following holds:
\begin{enumerate}
\item \emph{Pure Lasso with projection anti-concentration.} For constants
$0<\lambda_-\le\lambda_+<\infty$ independent of dimension,
\begin{equation}
\label{eq:lasso-strength-bounds}
 \rho(w)=\lambda\|w\|_1,\qquad
 0<\lambda_-\le\lambda\le\lambda_+<\infty.
\end{equation}
The strength $\lambda$ may depend on size, and there is a constant
$A_{\rm ac}<\infty$, independent of size and direction, such that
\begin{equation}
\label{eq:projection-anticoncentration}
 \Pp\!\left(|x^\top u|\le r\sqrt{u^\top \Sigma_xu}\right)\le A_{\rm ac}r,
 \qquad u\ne0,\quad 0<r<1.
\end{equation}
\item \emph{Uniform strong convexity.} The penalty satisfies
\begin{equation}
\label{eq:strongly-convex-penalty}
 w\longmapsto\rho(w)-\frac{\gamma_{\rm sc}}2|w|^2
 \quad\text{is convex},\qquad \gamma_{\rm sc}>0
 \quad\text{independent of }n.
\end{equation}
\end{enumerate}
\end{assumption}
Anti-concentration is required only in the pure-Lasso alternative, with no
smallness restriction on $A_{\rm ac}$.
The strongly convex alternative allows nonsmooth and nonseparable penalties,
with their minimum at zero as in Assumption~\ref{ass:analytic}.
It includes the elastic penalty
\begin{equation}
\label{eq:elastic-penalty}
 \rho(w)=\lambda\|w\|_1+\rho_0(w),\qquad
 \rho_0\in C^2(\R^p),\qquad
 \nabla^2\rho_0(w)\succeq\gamma_{\rm sc}I_p\quad\text{for all }w,
\end{equation}
with $\lambda,\gamma_{\rm sc}>0$ independent of $n$.
The remainder may depend on size and be nonseparable; no upper Hessian
bound is required.
% Assumption~\ref{ass:analytic} retains the minimum of the total penalty at
% zero, equivalent here to $\|\nabla\rho_0(0)\|_\infty\le\lambda$.
% The empirical objective is coercive in both cases. Its minimizer is unique
% in the strongly convex case; for pure Lasso, $\widehat W$ denotes the
% minimum-Euclidean-norm minimizer.

\subsection{Perturbed optimizers and prediction quantities}
\label{sec:perturbed-optimizers}
For $t\in[0,1]$, let $B_i(t)\sim\operatorname{Bernoulli}(t)$ retain
observations independently, and write $B_i:=B_i(t)$.
The selectors, training pairs, and probe
$G\sim N(0,\Sigma_x)$ are mutually independent. For $v,b>0$ and
$h\in\R^p$, write $z:=(v,b,h)$ and define
\begin{equation}
\label{eq:energy}
\begin{aligned}
 \cE_{t,z}(w):={}&\rho(w)
 +\sum_{i=1}^n B_i(t)L_{y_i}(x_i^\top w/\sqrt n)+\frac v2w^\top \Sigma_xw+\sqrt{2b}\,G^\top w-\sqrt n\,h^\top w.
\end{aligned}
\end{equation}
\begin{equation}
\label{eq:interior-estimator}
 W_{t,z}:=\argmin_w\cE_{t,z}(w).
\end{equation}
The parameters $v,b,h$ add quadratic curvature, Gaussian noise, and a
linear field to the retained empirical objective.
At $t=0$, no observation is retained, and the choice
$z=z_\ast:=(\nu_\ast,\beta_\ast^2/2,\eta_\ast)$ gives exactly the surrogate
$W_\ast(G)$ in \eqref{eq:intro-prox}.
At $t=1$, every observation is retained. Extending \eqref{eq:energy} to
$z=0$ recovers \eqref{eq:original} and its selected minimizer
$\widehat W$.

The Hamilton--Jacobi equation \eqref{eq:value-equation} in
Section~\ref{sec:HJ}, together with \eqref{eq:Hgradient}, determines the path
\begin{equation}
\label{eq:calibrated-path}
 z(t):=(v(t),b(t),h(t))=(1-t)z_\ast,
\end{equation}
with $z(0)=z_\ast$ and $z(1)=0$. This path is used to prove Theorem~\ref{thm:endpoint} below (more details are provided in Section~\ref{sec:HJ}).

For any square-integrable coefficient vector $W$, write
\begin{equation}
\label{eq:target-moments}
 \mu:=\E\left[W\right]/\sqrt n,\qquad
 \alpha^2:=\tr(\Sigma_x C_W)/n.
\end{equation}
Here $\mu$ is the normalized mean, and $\alpha^2$ is the prediction variance
over the randomness defining $W$, averaged over a fresh feature $x$.
For an interior optimizer, define
\begin{equation}
\label{eq:orders}
 \kappa:=\frac1{n\sqrt n}\E\left[\tr(\Sigma_xD_hW)\right].
\end{equation}
Here $D_hW$ is the Jacobian of the optimizer with respect to $h$, and
$\kappa$ is the average score response to the physical linear field
$\sqrt n\,h$. This Jacobian exists almost everywhere because positive
quadratic curvature makes the solution map Lipschitz; the Gaussian field
avoids its exceptional set (Lemma~\ref{lem:interior-wellposedness}).

\subsection{Conditional Gaussian prediction}
To compare score laws along the interpolation, we first specify the relevant
conditional law and distance. For a real random variable $u$ and a random
vector $x$, write $\Law(u\mid x)$ for a version of the regular conditional
law of $u$ given $x$, characterized by
\begin{equation}
\label{eq:conditional-law-definition}
 \int_{\R}f(a)\,\Law(u\mid x)(\dd a)
 =\E\left[f(u)\mid x\right]\quad\text{almost surely}
\end{equation}
for every bounded Borel function $f:\R\to\R$. For Borel probability laws
$P,Q$ on $\R$ with finite first moments, define their $1$-Wasserstein distance
by
\begin{equation}
\label{eq:wasserstein-definition}
 \Wass(P,Q):=\inf_{\pi\in\Pi(P,Q)}
 \int_{\R^2}|a-b|\,\pi(\dd a,\dd b)
 =\sup_{\operatorname{Lip}(f)\le1}
 \left|\int_{\R}f\,\dd P-\int_{\R}f\,\dd Q\right|,
\end{equation}
where $\Pi(P,Q)$ is the set of couplings of $P$ and $Q$, and the supremum is
over all $1$-Lipschitz functions $f:\R\to\R$.

The following theorem gives a conditional Gaussian approximation along the
complete calibrated path, including the empirical endpoint, with each
optimizer's actual mean and covariance.
\begin{revisedblock}
\begin{theorem}[Conditional Gaussian prediction under Poincar\'e concentration]
\label{thm:prediction-directions}
Along a sequence in which the training pairs are independent copies of
$(x,y)$, $p/n\to c\in(0,\infty)$,
and $\underline c_xI_p\preceq\Sigma_x\preceq KI_p$, suppose
Assumptions~\ref{ass:analytic}--\ref{ass:penalty-structure} hold and
\begin{equation}
\label{eq:population-moment-regime}
 \sup_n\left\{\frac1n\E_G\left[|W_\ast|^2\right]+\kappa_\ast\right\}<\infty,
 \qquad \liminf_{n\to\infty}\alpha_\ast^2>0.
\end{equation}
Let $x$ be an independent fresh feature. For $0\le t<1$, put
$W_t:=W_{t,(1-t)z_\ast}$, and put $W_1:=\widehat W$, using the
minimum-norm selection for pure Lasso. Define the actual moments
\[
 \mu_t:=\E\left[W_t\right]/\sqrt n,\qquad
 \alpha_t^2:=\tr(\Sigma_xC_{W_t})/n.
\]
Then
\begin{equation}
\label{eq:prediction-directions-path}
 \sup_{t\in[0,1]}
 \E_x\left[\Wass\!\left(
   \Law\!\left(\frac{x^\top W_t}{\sqrt n}\,\middle|\,x\right),
   N(x^\top\mu_t,\alpha_t^2)\right)\right]\longrightarrow0.
\end{equation}
\end{theorem}

\begin{figure}[!tbp]
 \centering
 \includegraphics[width=\linewidth]{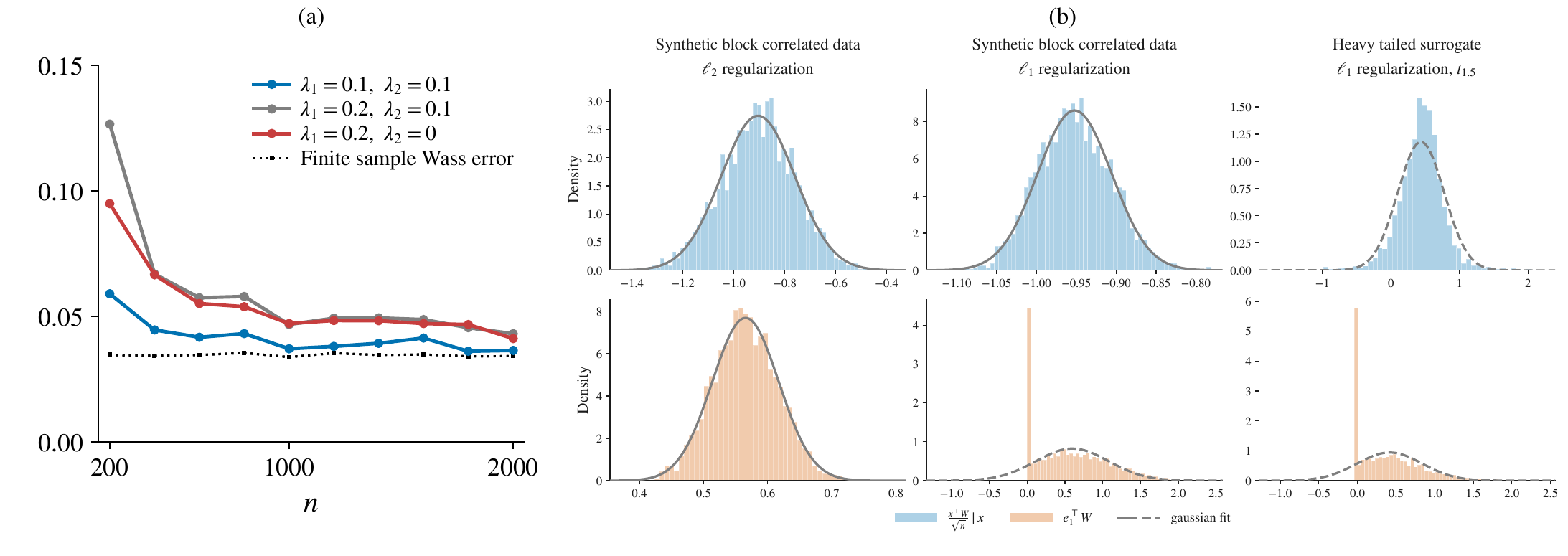}
 \revisioncolor
 \caption{(a) Standardized conditional Wasserstein distance, maximized over
 the interpolation grid, for the indicated elastic-net penalties; dotted
 curves are Gaussian finite-sampling benchmarks. (b) Conditional scores
 (top, blue) and first coefficients (bottom, orange) at $(n,p)=(200,400)$:
 pure ridge, pure Lasso, and heavy-tailed Lasso designs. Gray curves are
 moment-matched Gaussian densities. See
 Appendix~\ref{app:synthetic-prediction-experiments} for the settings.}
 \label{fig:synthetic-prediction}
\end{figure}

Proposition~\ref{prop:prediction-positive-probes} in
Appendix~\ref{app:pred-proof} gives a $Cn^{-1/8}$ conditional Gaussian
bound for full and unforced cavity optimizers under
Assumptions~\ref{ass:analytic}--\ref{ass:feature-concentration}, uniformly
in $t\in[0,1]$ with fixed positive lower bounds on both probes, including
nonsmooth penalties. For uniformly strongly convex penalties,
Appendix~\ref{app:pred-proof} adds and removes an auxiliary quadratic probe
and an independent Gaussian probe, uniformly over $t\in[0,1]$.
For pure Lasso, Appendix~\ref{app:proof-main-theorem} combines the
positive-probe comparison on $[0,1-\epsilon]$ with a uniform cutoff
coupling over the terminal interval $[1-\epsilon,1]$.
Both arguments yield uniform full-path convergence without an asserted rate.
The Hamilton--Jacobi and transport arguments use only the positive-probe
proposition.

Condition~\eqref{eq:population-moment-regime} holds for most practical scenarios. For instance, the following proposition gives a
sufficient condition for logistic Lasso.

\begin{proposition}[Population moment bounds for logistic Lasso]
\label{prop:logistic-lasso-moments}
Suppose $p/n\to c\in(0,\infty)$ and
$\underline c_xI_p\preceq\Sigma_x\preceq KI_p$.
Assume that Assumptions~\ref{ass:analytic}--\ref{ass:penalty-structure} hold with
$L_y(u)=\log(1+e^{-yu})$, $y\in\{-1,1\}$, and
$\rho(w)=\lambda_1\|w\|_1$ satisfying \eqref{eq:lasso-strength-bounds}.
If, for some $\epsilon\in(0,1/2]$,
\begin{equation}
\label{eq:logistic-label-overlap}
 \epsilon\le\Pp(y=1\mid x)\le1-\epsilon
 \qquad\text{almost surely},
\end{equation}
then the population solution satisfies
\eqref{eq:population-moment-regime}.
\end{proposition}
The proof is given in Appendix~\ref{app:logistic-lasso-moments}.
Condition~\eqref{eq:logistic-label-overlap} is a mild assumption requiring
that labels retain a nonvanishing degree of uncertainty given the features.

Figure~\ref{fig:synthetic-prediction} illustrates
Theorem~\ref{thm:prediction-directions} on synthetic data
(Appendix~\ref{app:synthetic-prediction-experiments}). In panel~(a), the
standardized conditional Wasserstein error approaches the Gaussian
finite-sampling benchmark as $n$ grows with $p=2n$, reaching the resolution
of the Monte Carlo estimates.

Panel~(b) shows score histograms over independently resampled training
datasets with $x$ held fixed. The conditional scores are close to Gaussian
for the Poincar\'e data, whereas the heavy-tailed design shows visible
departures. Gaussian fresh scores can nevertheless coexist with
non-Gaussian coefficient projections under pure $\ell_1$ regularization,
as the coordinate histograms illustrate; a rigorous example is given in
Proposition~\ref{prop:coordinate-nongaussian}. With pure $\ell_2$
regularization, the displayed coordinate is also close to Gaussian.
This agrees qualitatively with \citet[Theorem~2.1]{Louart2026}, which
establishes Gaussian limits for every nondegenerate deterministic
projection under quadratic regularization, a uniformly convex smooth
loss, and additional derivative bounds.
\end{revisedblock}

\subsection{Prediction moments and consequences}
\begin{theorem}[Prediction moments for Lasso and strongly convex penalties]
\label{thm:endpoint}
Along a sequence in which the training pairs are independent copies of
$(x,y)$, $p/n\to c\in(0,\infty)$,
and $\underline c_xI_p\preceq\Sigma_x\preceq KI_p$, suppose
Assumptions~\ref{ass:analytic}--\ref{ass:penalty-structure} and the
population regime \eqref{eq:population-moment-regime} hold. Then
\begin{equation}
\label{eq:main-moment-comparison}
 \left|\frac{\E\left[\widehat W\right]}{\sqrt n}-\mu_\ast\right|
 +\left|\frac1n\tr\!\left(\Sigma_x C_{\widehat W}\right)-\alpha_\ast^2\right|
 \longrightarrow0.
\end{equation}
\end{theorem}
Theorem~\ref{thm:endpoint} identifies the empirical optimizer's normalized
mean and prediction variance with those of the population surrogate.

\paragraph{Fresh prediction and classification.}
\begin{corollary}[Asymptotic prediction law]
\label{cor:conditional-prediction}
\label{cor:test-loss}
Under the assumptions of Theorem~\ref{thm:endpoint}, for an independent fresh
pair $(x,y)$ and for every measurable family of uniformly Lipschitz tests
$f_y$ for which both expectations in \eqref{eq:main-risk} are finite,
\begin{equation}
\label{eq:main-risk}
 \left|\E\left[f_y(x^\top\widehat W/\sqrt n)\right]
 -\E\left[f_y(x^\top\mu_\ast+\alpha_\ast Z)\right]\right|\longrightarrow0,
\end{equation}
where $Z\sim N(0,1)$ is independent
of $(x,y)$. The Lipschitz bound is common along the sequence.
\end{corollary}
The score distribution is taken over training samples, with the test
feature $x$ held fixed; the approximation error is then averaged over $x$.
The parameters $\mu_\ast$ and $\alpha_\ast^2$ are evaluated at each problem
size; neither $x^\top\mu_\ast$ nor $\alpha_\ast^2$ is required to converge.

\begin{corollary}[Asymptotic binary classification error]
\label{cor:classification}
Under the assumptions of Theorem~\ref{thm:endpoint}, let $(x,y)$ be an
independent fresh pair
with $y\in\{-1,1\}$ and set
$\widehat y(x):=\operatorname{sign}_+(x^\top\widehat W/\sqrt n)$,
where $\operatorname{sign}_+(z):=1$ for $z\ge0$ and $-1$ otherwise.
With $\Phi_{\rm N}$ the standard normal distribution function,
\begin{equation}
\label{eq:classification-risk}
 \left|\Pp\{\widehat y(x)\ne y\}
 -\E_{x,y}\left[\Phi_{\rm N}\!\left(
 -\frac{yx^\top\mu_\ast}{\alpha_\ast}\right)\right]\right|
 \longrightarrow0.
\end{equation}
The probability averages the training sample and the fresh pair.
\end{corollary}
The condition $\liminf_n\alpha_\ast^2>0$ in
\eqref{eq:population-moment-regime} bounds the comparison Gaussian's density,
so no additional margin or fitted-score density assumption is needed.
Corollaries~\ref{cor:conditional-prediction} and~\ref{cor:classification} are proved in Appendix~\ref{app:prediction-laws}.

\begin{revisedblock}
\begin{example}[Comparison with prior theoretical predictions on diamonds classification]
\label{sec:penalties}
\label{sec:diamonds-experiments}
We illustrate Corollaries~\ref{cor:conditional-prediction}
and~\ref{cor:classification} with pure-Lasso logistic classification on a
semisynthetic population derived from \texttt{diamonds}: predict whether
price exceeds \$3,000 from the other attributes. A fixed feature map,
Gaussian smoothing, and independent noise coordinates give
$(n,p)=(512,256)$. The fits use no ridge term or external probes;
the population surrogate uses the same penalty and is calibrated
independently of ERM moments. Appendix~\ref{app:diamonds-experiments}
gives the preprocessing and numerical methods, and
Example~\ref{sec:sparse-applications} gives the Lasso formulas.

\textit{Prediction and classification.}
Figure~\ref{fig:diamonds-prediction}(a) compares the classification formula
in Corollary~\ref{cor:classification} with population errors averaged
over independent training samples as $\lambda_1$ varies. Panels~(b,c)
illustrate Corollary~\ref{cor:conditional-prediction} at $\lambda_1=0.5$:
they compare the pooled fresh-score law with that of
$x^\top\mu_\ast+\alpha_\ast Z$, and isolate the centered fluctuations
$x^\top(\widehat W-\E[\widehat W])/\sqrt n$.
Here both the training sample and the independent test feature vary.
The predictions closely match the measured classification errors and
the score histograms, including the fluctuation scale revealed by centering.

\textit{Comparison with existing theories.}
The closest theoretical candidates for comparison are
\citet{Yaakoubi2026}, for smooth convex regularized ERM with non-Gaussian
data, and \citet{LoureiroMixtures2021}, for convex ERM on Gaussian-mixture
models. They address the two main difficulties of this example, but
neither directly explains the observed pure-Lasso behavior on the
original population.

The framework of \citet{Yaakoubi2026} requires a smooth, uniformly
strongly convex regularizer and does not cover exact Lasso. We evaluate
its equations with smooth approximations to $\|w\|_1$, using two smoothing
levels. These zero-ridge comparisons remain outside its hypotheses.
The resulting curves miss the pure-Lasso classification errors in
panel~(a). Even for the smoothed estimator, the equations underestimate
the centered score fluctuations in panel~(c), despite closely matching
the uncentered score law in panel~(b). The term
$x^\top\E[\widehat W]/\sqrt n$ dominates the uncentered score and can
therefore hide an inaccurate fluctuation variance.

The theory of \citet{LoureiroMixtures2021} permits Lasso, but requires
a Gaussian-mixture population. We fit one or ten Gaussian components
per class and apply its equations to those fitted populations.
For one component per class, the equations agree with ERM simulations
generated from the fitted mixture, yet the classification error differs
from that on the original data. Increasing the number to ten still
leaves a visible discrepancy. Thus the displayed mixture approximations
do not retain the aspects of the original feature--label law needed
to predict its learning curve. Our calibration retains that law and
the exact nonsmooth penalty.

\begin{figure}[!htbp]
 \centering
 \includegraphics[width=\linewidth]{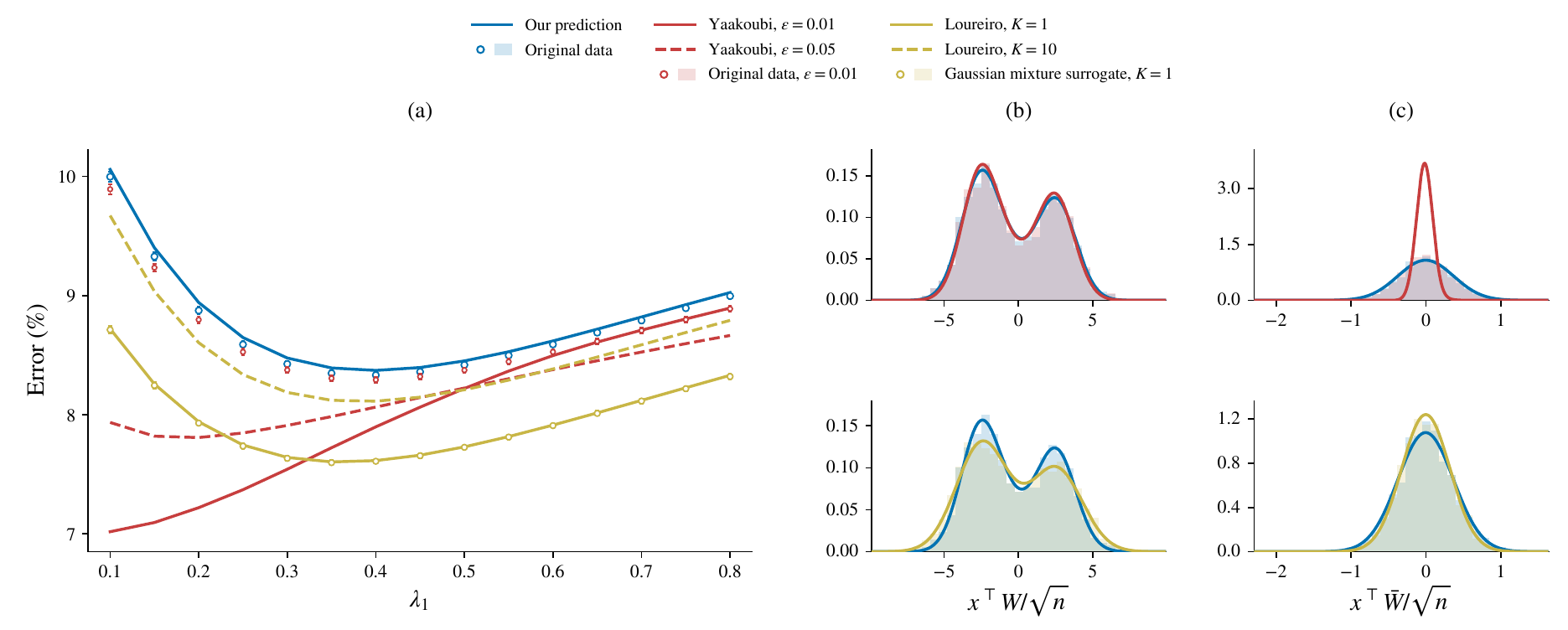}
 \revisioncolor
 \caption[Classification and pooled fresh-score laws on diamonds]{
 Classification and fresh-score comparisons on diamonds. Blue: pure Lasso;
 red: smooth-penalty equations~\citep{Yaakoubi2026}; yellow:
 Gaussian-mixture equations~\citep{LoureiroMixtures2021}.
 (a) Population classification errors with $\pm1.96$ Monte Carlo standard
 errors. (b,c) Pooled and centered scores at $\lambda_1=0.5$.
 Appendix~\ref{app:diamonds-experiments} specifies the line styles,
 comparison settings, and sampling counts.}
 \label{fig:diamonds-prediction}
\end{figure}

\end{example}
\end{revisedblock}

\section{Hamilton--Jacobi mechanism}
\label{sec:HJ}
\subsection{Pseudo-equation}
\label{sec:HJ-pseudo-equation}
Write the expected normalized minimum of the objective formula \eqref{eq:energy} as
\begin{equation}
\label{eq:value}
\phi(t,z):=\frac1n\E\left[\min_w\cE_{t,z}(w)\right].
\end{equation}
The parameter values $(t,z)=(1,0)$ give the empirical problem,
while $(t,z)=(0,z_\ast)$ give the calibrated surrogate.
The expected minimum $\phi$ records the prediction quantities through
the envelope identities
\begin{equation}
\label{eq:main-envelope}
 \mu=-\nabla_h\phi,\qquad
 \alpha^2=2\partial_v\phi-\mu^\top \Sigma_x\mu,\qquad
 \kappa=-\partial_b\phi.
\end{equation}
Thus the deterministic linear, quadratic, and Gaussian perturbations
reveal the mean, prediction variance, and response. Positive probes give
the required curvature without differentiating the penalty; the identities
are proved in Lemma~\ref{lem:interior-envelope}. We keep $v,b>0$ during the
comparison: $v=0,b>0$ can destroy Lasso coercivity.

Let $W_-$ be the optimizer with one observation omitted, and let $\mathcal I$
be the increase in the minimum objective when that independent observation
is reinserted. Lemma~\ref{lem:interior-envelope} gives
$\partial_t\phi=\E\left[\mathcal I\right]$.
Proposition~\ref{prop:r1}, using
Lemma~\ref{lem:derived-scalar-response}, replaces $\mathcal I$ by the envelope
$e_{L_y}(u;\kappa_-)$ from \eqref{eq:moreau}, where
$u:=x^\top W_-/\sqrt n$. Proposition~\ref{prop:r2} then applies
\eqref{eq:prediction-directions} in
Proposition~\ref{prop:prediction-positive-probes} to replace the fresh fluctuation
by a Gaussian with deterministic averaged variance.
Proposition~\ref{prop:r3} replaces the cavity moments by full-sample moments.
The actual feature--label law remains in the resulting Hamiltonian:
\begin{align}
\label{eq:H}
 \cH(\mu,\alpha,\kappa)&:=\E_{x,y,Z}\left[e_{L_y}(x^\top\mu+\alpha Z;\kappa)\right],\\
\label{eq:Hhat}
 \widetilde\cH(\pi)&:=\cH\!\left(-\pi_h,
 \sqrt{2\pi_v-\pi_h^\top \Sigma_x\pi_h},-\pi_b\right).
\end{align}
Here $\pi=(\pi_v,\pi_b,\pi_h)$ denotes the Hamiltonian argument,
which is evaluated at $\pi=\nabla_z\phi(t,z)$. The domain of
$\widetilde\cH$ is given by $2\pi_v\ge\pi_h^\top\Sigma_x\pi_h$
and $-\pi_b\ge0$.

For the Hamilton--Jacobi estimates, let $\mathcal D$ be a fixed interior
domain: $t$ ranges over a closed interval inside $(0,1)$, $v,b$ range over
compact subintervals of $(0,\infty)$, and $|h|$ is bounded, with all bounds
independent of $n$. Whenever finite differences are taken, slightly wider
bounds specify an enlarged domain.

\begin{theorem}[Asymptotic interior value equation]
\label{thm:main-interior-HJ}
Along a sequence in which the training pairs are independent copies of
$(x,y)$, suppose $p/n\to c\in(0,\infty)$,
$\underline c_xI_p\preceq\Sigma_x\preceq KI_p$, and
Assumptions~\ref{ass:analytic} and~\ref{ass:feature-concentration} hold.
On every fixed interior domain $\mathcal D$ as above,
\begin{equation}
\label{eq:value-equation}
 \sup_{\mathcal D}\left|
 \partial_t\phi-\widetilde\cH(\nabla_z\phi)\right|\longrightarrow0.
\end{equation}
Convergence in \eqref{eq:value-equation} does not use terminal quadratic growth or a population
calibration.
\end{theorem}

\subsection{Choice of the correct path that preserves statistics}
\label{sec:transport}
Write $\nabla:=\nabla_z$. By \eqref{eq:main-envelope}, fixed
$\nabla\phi$ means fixed $(\mu,\alpha^2,\kappa)$.
The gradient $\nabla_\pi\widetilde\cH$ is taken with respect to
$\pi=(\pi_v,\pi_b,\pi_h)$; evaluated at
$\pi=\nabla\phi(t,z)$, Lemma~\ref{lem:derivative-calculus}\textup{(\ref{item:hamiltonian-derivative-calculus})} gives
\begin{equation}
\label{eq:Hgradient}
 \nabla_\pi\widetilde\cH\!\left(\nabla\phi(t,z)\right)
 =(\nu,\beta^2/2,\eta),\qquad
 \nu\alpha=\E\left[Z\xi\right].
\end{equation}
These coefficients use the current moments in \eqref{eq:losscoeff},
before identification with the starred population solution.
First suppose $\phi$ is smooth and satisfies the exact equation
$\partial_t\phi=\widetilde\cH(\nabla\phi)$. Along a path, the chain rule gives
\[
\begin{aligned}
 \frac{\dd}{\dd t}\nabla\phi(t,z(t))
 &=\nabla\widetilde\cH(\nabla\phi)
   +\nabla^2\phi\,\dot z(t)\\
 &=\nabla^2\phi\,
   \bigl(\nabla_\pi\widetilde\cH(\nabla\phi)
            +\dot z(t)\bigr).
\end{aligned}
\]
By \eqref{eq:Hgradient}, choosing
$\dot z(t)=-(\nu,\beta^2/2,\eta)$ cancels the two terms,
so $\dd(\nabla\phi)/\dd t=0$.
For the calibrated initial point $z(0)=z_\ast$, the population calibration
equations \eqref{eq:losscoeff}--\eqref{eq:fixedpoint}, together with
\eqref{eq:Hgradient}, give
$\nabla_\pi\widetilde\cH(\nabla\phi(0,z_\ast))=z_\ast$.
Because \eqref{eq:losscoeff} determines $(\nu,\beta,\eta)$ from
$(\mu,\alpha,\kappa)$, the constancy of the latter along this formal
characteristic keeps all six parameters equal to their initial calibrated
values $(\mu_\ast,\alpha_\ast,\kappa_\ast,\nu_\ast,\beta_\ast,\eta_\ast)$.
Consequently, the characteristic is $z(t)=(1-t)z_\ast$ as in
\eqref{eq:calibrated-path}. Increasing $t$ adds observations while reducing
the probes and preserving the moments.

\section{Final discussion}
Our prediction results do not establish agreement between surrogate and
ERM support statistics: selection is discontinuous at zero. We nevertheless
tested whether the surrogate's coordinate-selection and support-recovery
predictions agree with ERM simulations, and found close agreement
(Appendix~\ref{app:support-simulation}). Rigorous coefficient and model-size
results exist for squared-loss Lasso with Gaussian designs
\citep{BayatiMontanari2012,CelentanoMontanariWei2023} and, under additional
conditions, independent non-Gaussian entries \citep{HanShen2023}.
Transfer under our broader hypotheses remains open; finer coordinate
comparisons and control of selection margins near the penalty threshold
may provide a route.

Further directions include extending
Theorem~\ref{thm:prediction-directions} to weaker concentration hypotheses
and describing all empirical minimizers when uniqueness fails, through
extensions of the Kac--Rice methods of
\citet[Assumptions~3--4]{AsgariMontanariSaeed2025} to non-Gaussian
covariates and nonsmooth regularizers.

\clearpage
\section*{Code availability}
The code used to reproduce the numerical experiments is available at
\url{https://anonymous.4open.science/r/ERM_non_smooth_penalty-7437/}.

\section*{AI use statement}
Generative AI tools assisted with mathematical development, literature
search, drafting including proof writing, and code and verification checks.
The authors are responsible for reviewing the mathematical claims, proofs
and cited sources, and for the final manuscript.
\section*{Reproducibility statement}
The appendices prove the stated asymptotic comparisons using convex
response calculus, value comparisons, calibrated transport and cutoff
coupling. Proposition~\ref{prop:prediction-positive-probes} proves the
positive-probe Gaussian comparison in Appendix~\ref{app:pred-proof}.
The qualitative full calibrated-path conclusion of
Theorem~\ref{thm:prediction-directions} uses the additional penalty/design
and population hypotheses and is proved in
Appendices~\ref{app:pred-proof} and~\ref{app:proof-main-theorem}.
Theorem~\ref{thm:endpoint} identifies the empirical endpoint moments with
the population moments.
Lemma~\ref{lem:lasso-bounded-root-localization} proves localization;
Propositions~\ref{prop:strongly-convex-calibrated-cutoff}
and~\ref{prop:lasso-calibrated-cutoff-MSE} prove cutoff mean-square
comparison. The source package includes compilation
instructions, a proof-status ledger, a numerical protocol, and executable
checks of the diagonal proximal formulas.
Appendix~\ref{app:prediction-experiments} specifies the feature laws,
preprocessing, calibration, and replication protocols for the prediction
experiments. The accompanying notebook reproduces the two main-text figures
from saved numerical data. Appendix~\ref{app:support-simulation}
gives the additional support-recovery results, penalty grid, and
support metrics.
\begingroup
\raggedright
% BEGIN bibliography

% END bibliography
\endgroup
\clearpage
% Continuous landscape proof map; page breaks occur automatically between rows.
\clearpage
\newgeometry{left=0.35in,right=0.35in,top=0.35in,bottom=0.35in}
\begin{landscape}
\begingroup
\iclrfinaltrue % Suppress the review line ruler on these landscape pages only.
\pagestyle{empty}
\fontsize{9.5}{10.5}\selectfont
\setlength{\parindent}{0pt}
\setlength{\tabcolsep}{5pt}\newlength{\pmwidth}\setlength{\pmwidth}{\dimexpr\linewidth-30pt\relax}
\setlength{\abovedisplayskip}{1.5pt}
\setlength{\belowdisplayskip}{1.5pt}
\setlength{\abovedisplayshortskip}{1pt}
\setlength{\belowdisplayshortskip}{1pt}
\renewcommand{\arraystretch}{1.0}
\newcommand{\pmnode}[2]{\textbf{\color{pmblue}#1}\par #2}
\newcommand{\pmref}[1]{{\footnotesize\color{pmgray}#1}}
\phantomsection\label{sec:key-proof-map}
\pdfbookmark[1]{Key Proof Table}{proof-map}
{\Large\bfseries Key Proof Table}\par\vspace{5pt}
{\small This table summarizes the key steps in the proofs of Theorem~\ref{thm:population-uniqueness}, Theorem~\ref{thm:prediction-directions}, Theorem~\ref{thm:endpoint} and its Corollaries~\ref{cor:conditional-prediction}--\ref{cor:classification}, and Theorem~\ref{thm:main-interior-HJ}.}\par
\setlength{\LTpre}{5pt}\setlength{\LTpost}{0pt}
\begin{longtable}{@{}>{\raggedright\arraybackslash}p{0.22\pmwidth}p{0.55\pmwidth}>{\raggedright\arraybackslash}p{0.11\pmwidth}>{\raggedright\arraybackslash}p{0.12\pmwidth}@{}}
\toprule
\textbf{Key Results}&\textbf{Proof Idea}&\textbf{Technical Tools}&\textbf{Reference}\\\midrule
\endfirsthead
\toprule
\textbf{Key Results}&\textbf{Proof Idea}&\textbf{Technical Tools}&\textbf{Reference}\\\midrule
\endhead
\bottomrule
\endfoot
\bottomrule
\endlastfoot
\multicolumn{4}{@{}p{\linewidth}@{}}{\textbf{Appendix~\ref{app:population-stability}}\enspace\textit{Establish existence and uniqueness of the solution to the population system \eqref{eq:losscoeff}--\eqref{eq:fixedpoint}.}}\\*[2pt]
\noalign{\nobreak\vskip\aboverulesep\hrule height\lightrulewidth\nobreak\vskip\belowrulesep\nobreak}
\pmnode{A1. Existence of a solution to the population system}{A finite solution
$(\mu_\ast,\alpha_\ast,\beta_\ast,\nu_\ast,\kappa_\ast,\eta_\ast)$
to \eqref{eq:losscoeff}--\eqref{eq:fixedpoint} exists, with $\alpha_\ast,\beta_\ast,\nu_\ast,\kappa_\ast>0$.}
& Minimize the auxiliary convex objective in \eqref{eq:population-convex-program} over $(W,\zeta,\alpha)$. Attainment and Slater's condition yield KKT multipliers; positivity and stationarity identify the surrogate optimizer and the complete population system:
\par\noindent\makebox[\linewidth][c]{$\displaystyle\text{attainment}+\text{KKT}+\text{positivity}\ \Longrightarrow\ \eqref{eq:losscoeff}\text{--}\eqref{eq:fixedpoint}.$}
& Convex optimization,\par Slater's condition,\par KKT conditions
& Theorem~\ref{thm:population-uniqueness}\par Appendix~\ref{app:population-stability}.\\
\midrule
\pmnode{A2. Uniqueness of the solution to the population system}{If
$(\mu,\alpha,\beta,\nu,\kappa,\eta)$ and
$(\mu',\alpha',\beta',\nu',\kappa',\eta')$
are finite solutions to \eqref{eq:losscoeff}--\eqref{eq:fixedpoint}
with $\alpha,\kappa\ge0$ and $\beta,\nu>0$ (likewise for the primed solution), then
\par\noindent\makebox[\linewidth][c]{$\displaystyle\thinmuskip=2mu (\mu,\alpha,\beta,\nu,\kappa,\eta)=(\mu',\alpha',\beta',\nu',\kappa',\eta').$}}
& Both solutions minimize \eqref{eq:population-convex-program}. Their feasible midpoint gives
\[\E\!\left[\frac{L_y(\zeta)+L_y(\zeta')}{2}-L_y\!\left(\frac{\zeta+\zeta'}2\right)\right]=0
 \ \Rightarrow\ \xi=\xi'.\]
The proximal identity and the parameter dependencies in the population system then identify the remaining parameters.
& Midpoint optimality,\par Strict Jensen equality,\par Proximal identity
& Theorem~\ref{thm:population-uniqueness}\par Eq.~\eqref{eq:population-uniqueness-jensen-gap}.\\
\midrule
\multicolumn{4}{@{}p{\linewidth}@{}}{\textbf{Appendix~\ref{app:calculus}}\enspace\textit{Derive formulas for the derivatives of $\phi$ with respect to $z$ and $t$.}}\\*[2pt]
\noalign{\nobreak\vskip\aboverulesep\hrule height\lightrulewidth\nobreak\vskip\belowrulesep\nobreak}
\pmnode{B1. Value gradient identities}{
$\mu=-\nabla_h\phi$,\quad $\kappa=-\partial_b\phi$,\par
$\alpha^2=2\partial_v\phi-\mu^\top\Sigma_x\mu$.}
& The envelope theorem gives
\[\nabla_h\phi=-\frac{\E\left[W\right]}{\sqrt n},\qquad \partial_v\phi=\frac{\E\left[W^\top\Sigma_xW\right]}{2n},\qquad
\partial_b\phi=\frac{\E\left[G^\top W\right]}{n\sqrt{2b}}=-\kappa.\]
The last equality follows by Gaussian integration by parts; centering the quadratic moment yields $\alpha^2$.
& Envelope theorem,\par Gaussian integration by parts
& Lemma~\ref{lem:normalized-identities}\par \eqref{eq:envelope}.\\
\midrule
\pmnode{B2. Time derivative identity}{$\partial_t\phi=\E\left[\mathcal I\right]$.}
& Differentiating the product Bernoulli law, then using exchangeability, gives
\[\partial_t\phi=\frac1n\sum_{i=1}^n\E\left[\mathcal I_i\right]=\E\left[\mathcal I\right],\]
where $\mathcal I_i$ is the minimum-value increment upon reinserting observation $i$.
& Bernoulli differentiation,\par Exchangeability
& Lemma~\ref{lem:interior-envelope}\par insertion cost \eqref{eq:increment}.\\
\midrule
\multicolumn{4}{@{}p{\linewidth}@{}}{\textbf{Appendix~\ref{app:pred-proof}}\enspace\textit{Prove Proposition~\ref{prop:prediction-positive-probes} on $\mathcal D_+=[0,1]\times\mathcal Z_n$, and Theorem~\ref{thm:prediction-directions} for uniformly strongly convex penalties. The pure-Lasso full-path conclusion is completed in Appendix~\ref{app:proof-main-theorem}.}}\\*[2pt]
\noalign{\nobreak\vskip\aboverulesep\hrule height\lightrulewidth\nobreak\vskip\belowrulesep\nobreak}
\pmnode{C1. Optimizer concentration}{For $\Delta:=W-\E\left[W\right]$, the covariance operator is bounded and deterministic quadratic forms of $\Delta$ concentrate after division by $n$.}
& Remove one observation and integrate the response over the full bounded force interval. A fixed-radius feature-tail bound and Gaussian translation control the directional and quadratic insertion increments. Efron--Stein over observations and Gaussian Poincar\'e over the probe give
\[\|C_W\|_{\op}\le C_{\mathcal D_+},\qquad
\sup_{\substack{A=A^\top\\\|A\|_{\op}\le1}}
\left\|\frac{\Delta^\top A\Delta-\tr(AC_W)}n\right\|_{L^2}
\le C_{\mathcal D_+}n^{-1/2}.\]
Since $\tr C_W=O(n)$, this also yields $\|C_W\|_{\op}/n\to0$ and $\|C_W\|_{\hs}/n\to0$.
& Leave-one-out,\par Gaussian translation,\par Efron--Stein,\par Gaussian Poincar\'e
& Lemma~\ref{lem:pred-proof-concentration}\par \eqref{eq:pred-proof-concentration}.\\
\midrule
\pmnode{C2. A Stein kernel from a spectral gap}{With $m_x:=\E\left[x\right]$, $\bar x:=x-m_x$, and $C_x:=\Cov(x)$, construct a Stein kernel $\mathcal T_x$ satisfying
\[\E\left[\mathcal T_x\right]=C_x,\]
\[\E\left[\|\mathcal T_x\|_{\hs}^2\right]\le C_{\rm P}\tr C_x.\]}
& The marginal Poincar\'e inequality makes $\nabla f\mapsto\E\left[\bar x_i f(\bar x)\right]$ continuous on the $L^2$ closure of gradients. Riesz representation gives each row of $\mathcal T_x$ and the identity
\[\E\left[\bar x f(\bar x)\right]=\E\left[\mathcal T_x(\bar x)\nabla f(\bar x)\right].\]
Coordinate tests identify its expectation. No feature density is required.
& Poincar\'e inequality,\par Riesz representation,\par Stein identity
& Lemma~\ref{lem:pred-proof-kernel}\par \eqref{eq:pred-proof-kernel}.\\
\midrule
\pmnode{C3. Stein's Gaussian approximation strategy}{For $a>0$ and a smooth $1$-Lipschitz function $F$, the Stein equation has a solution with
\[\|f'\|_\infty\le\sqrt{2/\pi}/a,\]
\[\|f''\|_\infty\le2/a^2.\]}
& Solve $a^2 f'(u)-uf(u)=F(u)-\E\left[F(aZ)\right]$. Set $T:=\bar x^\top\Delta/\sqrt n$ and $a^2:=\tr(C_xC_W)/n$. The Stein kernel from C2, its spectral truncation, and C1 yield
\[|\E\left[Tf(T)\right]-a^2\E\left[f'(T)\right]|
\le C_{\mathcal D_+}\bigl(n^{-1/4}\|f'\|_\infty+n^{-1/2}\|f''\|_\infty\bigr).\]
Apply the solution bounds and Wasserstein duality when $a\ge n^{-1/8}$; for smaller $a$, use the second moment directly. Thus $\Wass(\Law(T),N(0,a^2))\le C_{\mathcal D_+}n^{-1/8}$.
& Chatterjee's Stein equation,\par Spectral truncation,\par Wasserstein duality
& Lemma~\ref{lem:pred-proof-stein-equation}\par \eqref{eq:pred-proof-stein}.\\
\midrule
\pmnode{C4. Gaussian approximation of fresh scores}{With $x$ held fixed, approximate the score distribution over training data, subsampling, and Gaussian probes by $N(x^\top\mu,\alpha^2)$. The Wasserstein error averaged over $x$ is $O(n^{-1/8})$ uniformly on $\mathcal D_+$, for full and unforced cavity optimizers.}
& Smooth the conditional CDF with Gaussian noise. Poincar\'e and C1 control its fluctuations; integrate on a bounded interval and control the tails by second moments. Letting the smoothing scale decrease controls the mean Wasserstein distance between the score law with the feature held fixed and the law averaged over features. Apply C3, then restore the feature mean and the variance $\alpha^2=\tr(\Sigma_xC_W)/n$ by coupling. The argument also covers $a=0$.
& Gaussian smoothing,\par Poincar\'e inequality,\par CDF formula for Wasserstein distance,\par Coupling
& Proposition~\ref{prop:prediction-positive-probes}\par \eqref{eq:prediction-directions}.\\
\midrule
\pmnode{C5. Strongly convex full-path comparison}{The full optimizer's Gaussian error tends to zero uniformly on $t\in[0,1]$ along the calibrated path.}
& Add independent positive quadratic and Gaussian probes of size $\delta$. Strong monotonicity and uniform second moments give
\[\sup_{t\in[0,1]}\E|W_t^\delta-W_t|^2/n\le C\delta.\]
Apply Proposition~\ref{prop:prediction-positive-probes} at fixed $\delta$, then the Gaussian-comparison coupling from B; take $n\to\infty$ before $\delta\downarrow0$.
& Strong monotonicity,\par Independent Gaussian perturbation,\par $L^2$ coupling
& Theorem~\ref{thm:prediction-directions}\par \eqref{eq:strongly-convex-probe-coupling}.\\
\midrule
\multicolumn{4}{@{}p{\linewidth}@{}}{\textbf{Appendix~\ref{app:four}}\enspace\textit{Prove Theorem~\ref{thm:main-interior-HJ} by decomposing the approximate Hamilton--Jacobi equation into three errors:}
\[
\begin{aligned}
&r(t,z):=\partial_t\phi-\widetilde\cH(\nabla_z\phi)
 \overset{\mathrm{B2}}{=}\E\left[\mathcal I\right]-\widetilde\cH(\nabla_z\phi)=\E\left[\mathcal I-e_{L_y}(x^\top W_-/\sqrt n;\kappa_-)\right]
 +\big\{\E\left[e_{L_y}(x^\top W_-/\sqrt n;\kappa_-)\right]-\cH(\mu_-,\alpha_-,\kappa_-)\big\}+\big\{\cH(\mu_-,\alpha_-,\kappa_-)-\widetilde\cH(\nabla_z\phi)\big\}:=r^{(1)}+r^{(2)}+r^{(3)}
\end{aligned}
\]
\textit{D1, D2, and D3 show that $r^{(1)}$, $r^{(2)}$, and $r^{(3)}$ converge to zero uniformly on $\mathcal D$, respectively.}
}\\*[2pt]
\noalign{\nobreak\vskip\aboverulesep\hrule height\lightrulewidth\nobreak\vskip\belowrulesep\nobreak}
\pmnode{D1. Scalar-response approximation}{
$r^{(1)}=\E\left[\mathcal I\right]-\E\left[e_{L_y}(x^\top W_-/\sqrt n;\kappa_-)\right]$,\par
$\sup_{\mathcal D}|r^{(1)}|\to0$.}
& Let $W_-$ be the unforced cavity optimizer, with moments $(\mu_-,\alpha_-^2,\kappa_-)$. Let $W_-(r)$ minimize the cavity objective plus $-rx^\top w/\sqrt n$. Response concentration and Gaussian change of measure yield
\[\E\sup_{|r|\le L_1}\left|\frac{x^\top(W_-(r)-W_-)}{\sqrt n}-\kappa_-r\right|\longrightarrow0.\]
Integrate this response in $r$. Fenchel duality expresses both insertion cost and the Moreau envelope as scalar suprema over $[-L_1,L_1]$; their difference is bounded by $L_1$ times the displayed error.
& Leave-one-out,\par Response concentration,\par Change of measure,\par Fenchel duality
& Proposition~\ref{prop:r1}\par finite-force linearization and Fenchel comparison.\\
\midrule
\pmnode{D2. Gaussian approximation of the cavity score}{
$r^{(2)}=\E\left[e_{L_y}(x^\top W_-/\sqrt n;\kappa_-)\right]-\cH(\mu_-,\alpha_-,\kappa_-)$,\par
$\sup_{\mathcal D}|r^{(2)}|\to0$.}
& The Moreau envelope is $L_1$-Lipschitz. Condition on the deleted pair $(x,y)$, independent of $W_-$, and apply Wasserstein duality:
\[|r^{(2)}|\le L_1\E_x\Wass\!\left(\Law\!\left(\frac{x^\top W_-}{\sqrt n}\,\middle|\,x\right),
 N(x^\top\mu_-,\alpha_-^2)\right)\longrightarrow0.\]
This applies only Proposition~\ref{prop:prediction-positive-probes} in the three-error decomposition.
& Proposition~\ref{prop:prediction-positive-probes},\par Lipschitz envelope,\par Wasserstein duality
& Proposition~\ref{prop:r2}.\\
\midrule
\pmnode{D3. Cavity-to-full comparison}{
$r^{(3)}=\cH(\mu_-,\alpha_-,\kappa_-)-\widetilde\cH(\nabla_z\phi)$,\par
$\sup_{\mathcal D}|r^{(3)}|\to0$.}
& Couple $W$ and $W_-$ with the same observations and probe. Strong monotonicity gives $D:=\|W-W_-\|_{L^2}\le C_{\mathcal D}$. The moment bounds and Gaussian-probe correlation identity yield
\[\begin{aligned}
|\mu-\mu_-|&\le D/\sqrt n,\qquad
|\alpha^2-\alpha_-^2|\le C_{\mathcal D}D/\sqrt n,\\
|\kappa-\kappa_-|&=\frac{|\E\left[G^\top(W-W_-)\right]|}{n\sqrt{2b}}
\le\frac{\sqrt{\tr(\Sigma_x)/n}\,D}{\sqrt{2bn}}.
\end{aligned}\]
Hamiltonian stability then gives $|r^{(3)}|\le C_{\mathcal D}D/\sqrt n\to0$ uniformly on $\mathcal D$.
& Strong monotonicity,\par Optimizer coupling,\par Moment and response stability
& Proposition~\ref{prop:r3}.\\
\midrule
\multicolumn{4}{@{}p{\linewidth}@{}}{\textbf{Appendix~\ref{app:value-transport}}\enspace\textit{For fixed $\epsilon>0$, transport calibrated prediction moments along the characteristic path up to the cutoff $t=1-\epsilon$.}}\\*[2pt]
\noalign{\nobreak\vskip\aboverulesep\hrule height\lightrulewidth\nobreak\vskip\belowrulesep\nobreak}
\pmnode{E1. Auxiliary value error}{Set
\[\pi_\ast:=\nabla_z\phi(0,z_\ast).\]
For $A$ sufficiently large, define
\[\begin{aligned}
&\mathcal E(t):=\max_{|d|\le r_0}\{\phi(t,z(t)+d)\\[-2pt]
&\quad-\phi(t,z(t))-\pi_\ast^\top d-A|d|^2/2\}.
\end{aligned}\]
Here $r_0>0$ is a fixed radius chosen so that $z(t)+d$ remains in the positive-probe neighborhood whenever $|d|\le r_0$.
The auxiliary value bounds the squared gradient error from above and below.}
& Let $L$ be the Lipschitz constant of $\nabla_z\phi$ on this neighborhood. Quadratic Taylor bounds give
\[\begin{aligned}
\langle\nabla_z\phi(t,z(t))-\pi_\ast,d\rangle-\tfrac{A+L}{2}|d|^2
&\le\phi(t,z(t)+d)-\phi(t,z(t))-\pi_\ast^\top d-\tfrac A2|d|^2\\
&\le\langle\nabla_z\phi(t,z(t))-\pi_\ast,d\rangle-\tfrac{A-L}{2}|d|^2.
\end{aligned}\]
Maximize the upper quadratic; for the lower bound test
$d=(\nabla_z\phi(t,z(t))-\pi_\ast)/(A+L)$, which lies in the ball. Thus
\[\frac{|\nabla_z\phi(t,z(t))-\pi_\ast|^2}{2(A+L)}
\le\mathcal E(t)\le\frac{|\nabla_z\phi(t,z(t))-\pi_\ast|^2}{2(A-L)}.\]
& Lipschitz gradient,\par Quadratic Taylor bounds,\par Auxiliary maximization
& Lemma~\ref{lem:lasso-calibrated-transport}\par \eqref{eq:direct-transport-energy-bounds}.\\
\midrule
\pmnode{E2. Calibrated gradient transport}{For each fixed $\epsilon>0$,
\[\nabla_z\phi(1-\epsilon,z(1-\epsilon))-\pi_\ast\to0,\]
so
\[\begin{aligned}
\big|\mu^\epsilon-\mu_\ast\big|&\to0,\\
| (\alpha^\epsilon)^2-\alpha_\ast^2|+|\kappa^\epsilon-\kappa_\ast|&\to0.
\end{aligned}\]}
& The Hamiltonian derivative gives $\dot z=-\nabla\widetilde\cH(\pi_\ast)=-z_\ast$.
For fixed $s,\epsilon>0$, set $e_{n,s,\epsilon}:=\sup_{s\le t\le1-\epsilon,\,|d|\le r_0}|r(t,z(t)+d)|$. The interior maximizer in E1 and the envelope theorem yield
$\mathcal E'(t)\le C_0\mathcal E(t)+2e_{n,s,\epsilon}$, without differentiating $r$. Apply Gr\"onwall between the two sides of E1:
\[\begin{aligned}
|\nabla_z\phi(t,z(t))-\pi_\ast|^2
&\le 2(A+L)\mathcal E(t)\\
&\le 2(A+L)e^{C_0(t-s)}\mathcal E(s)+C_\epsilon e_{n,s,\epsilon}\\
&\le\frac{A+L}{A-L}e^{C_0(t-s)}|\nabla_z\phi(s,z(s))-\pi_\ast|^2
 +C_\epsilon e_{n,s,\epsilon}.
\end{aligned}\]
Direct coupling \eqref{eq:initial-short-time-coupling} gives the initial squared error $\le Cs$; D1--D3 give $e_{n,s,\epsilon}\to0$ for fixed $s,\epsilon$. Hence
$\limsup_n|\nabla_z\phi(1-\epsilon,z(1-\epsilon))-\pi_\ast|^2\le C_\epsilon s$.
Let $s\downarrow0$; B1 converts gradient convergence into moment convergence.
& Hamiltonian characteristic,\par Envelope theorem,\par Gr\"onwall inequality,\par Direct coupling
& Lemma~\ref{lem:lasso-calibrated-transport}\par \eqref{eq:lasso-calibrated-moments}.\\
\midrule
\multicolumn{4}{@{}p{\linewidth}@{}}{\textbf{Appendix~\ref{app:proof-main-theorem}}\enspace\textit{Uniform terminal-interval coupling completes Theorem~\ref{thm:prediction-directions} for pure Lasso. Endpoint cutoff removal then proves the moment comparison in Theorem~\ref{thm:endpoint}.}}\\*[2pt]
\noalign{\nobreak\vskip\aboverulesep\hrule height\lightrulewidth\nobreak\vskip\belowrulesep\nobreak}
\pmnode{F1. Strongly convex cutoff removal}{For $0<\epsilon<1$,
$\displaystyle\frac1n\E|\widehat W-W_{1-\epsilon,\,\epsilon z_\ast}|^2\le C\epsilon$.}
& Let $\overline W:=\argmin_w F_I(w)$ be the unprobed retained-data optimizer, with $I=I(1-\epsilon)$. Strong convexity and independence of the omitted pairs give
\[\frac1n\E|\widehat W-\overline W|^2\le C\epsilon,
\qquad \frac1n\E|W_{1-\epsilon,\,\epsilon z_\ast}-\overline W|^2\le C\epsilon.\]
The first bound restores observations; the second removes probes using strong monotonicity. Combine by the squared triangle inequality.
& Strong convexity,\par Objective comparison,\par Optimizer coupling
& Proposition~\ref{prop:strongly-convex-calibrated-cutoff}.\\
\midrule
\pmnode{F2. Pure-Lasso terminal-interval comparison}{For $V_\delta:=W_{1-\delta,\delta z_\ast}$ and $V_0:=\widehat W$,
\[\limsup_n\sup_{0\le\delta\le\epsilon}\E|V_\delta-V_\epsilon|^2/n\le C\epsilon^{1/4}.\]}
& Use common nested selectors and the same Gaussian probe. Square completion and common localization control exceptional moments uniformly in $\delta$. At the fixed reference, unconditional omitted-loss and signed Gaussian estimates give $\sup_\delta\E D_{\epsilon,\delta}/n\le C\epsilon+o_n(1)$. Reference support/slack and uniform curvature give $D_{\epsilon,\delta}\ge ca|V_\delta-V_\epsilon|^2$ on a common event; its remaining failure probability satisfies $\limsup_n\Pp(\text{failure})\le Ca/\sqrt\epsilon$ at fixed $\epsilon,a$. Hence
\[\limsup_n\sup_\delta\E|V_\delta-V_\epsilon|^2/n
 \le C(\epsilon/a+a/\sqrt\epsilon).\]
After the dimension limit choose $a=\epsilon^{3/4}$. Transfer the Gaussian error through B, combine with Proposition~\ref{prop:prediction-positive-probes} on $[0,1-\epsilon]$, then let $\epsilon\downarrow0$ to prove Theorem~\ref{thm:prediction-directions}.
& Fixed-reference convexity remainder,\par Gaussian sign,\par Curvature or slack,\par Uniform exceptional moments
& Proposition~\ref{prop:lasso-calibrated-cutoff-MSE}\par Theorem~\ref{thm:prediction-directions}.\\
\midrule
\pmnode{F3. Endpoint moment identification}{The empirical mean and prediction variance agree asymptotically with the calibrated population moments.}
& The endpoint specialization of F2 or the strongly convex comparison F1 removes the cutoff in normalized $L^2$. At each fixed cutoff, E2 identifies the actual mean and prediction variance with the population moments. The mean bound and reverse triangle inequality for centered weighted $L^2$ norms pass these moments to the empirical endpoint; take $n\to\infty$ before $\epsilon\downarrow0$.
& $L^2$ coupling,\par Fixed-cutoff transport,\par Moment stability
& Theorem~\ref{thm:endpoint}\par \eqref{eq:main-moment-comparison}.\\
\midrule
\multicolumn{4}{@{}p{\linewidth}@{}}{\textbf{Appendix~\ref{app:prediction-laws}}\enspace\textit{Derive prediction and classification laws from the endpoint Gaussian comparison and identified moments.}}\\*[2pt]
\noalign{\nobreak\vskip\aboverulesep\hrule height\lightrulewidth\nobreak\vskip\belowrulesep\nobreak}
\pmnode{G1. Prediction and classification}{Transfer uniformly Lipschitz tests and binary classification risk to $x^\top\mu_\ast+\alpha_\ast Z$.}
& Theorem~\ref{thm:prediction-directions} at $t=1$ gives the Gaussian comparison with actual moments. Theorem~\ref{thm:endpoint} replaces them by population moments. Lipschitz tests transfer directly; half-line ramps of width $\ell$ give classification discrepancy $\le\Wass/\ell+\ell/(\sqrt{2\pi}\,\underline\alpha)$. The normalized Hilbert--Schmidt covariance bound separately uses cutoff coupling and C1.
& Gaussian comparison,\par Moment identification,\par Lipschitz tests,\par Threshold smoothing
& Corollaries~\ref{cor:conditional-prediction}--\ref{cor:classification}\par Remark~\ref{rem:boundary-covariance}.\\
\end{longtable}
\clearpage
\endgroup
\end{landscape}
\restoregeometry

\appendix
% Continue the shared result sequence within each lettered appendix.
\numberwithin{theorem}{section}
\numberwithin{remark}{section}

\section{Population well-posedness and coefficient bounds}
\label{app:population-stability}
\label{app:fixedpoint-stability}

Appendices~\ref{app:population-existence}--\ref{app:population-uniqueness-proof}
prove Theorem~\ref{thm:population-uniqueness} using the convex program
\eqref{eq:population-convex-program}: attainment and the
Karush--Kuhn--Tucker (KKT) conditions give existence, and a feasible
midpoint gives uniqueness. Proposition~\ref{prop:logcosh-nonexistence}
shows why score nondegeneracy is needed for existence.
Appendix~\ref{app:population-stationarity} gives an equivalent stationary-point
formulation. Appendix~\ref{app:population-coefficients} collects population moment
identities and uniform coefficient bounds;
Appendix~\ref{app:logistic-lasso-moments} proves the logistic-Lasso
sufficient condition in Proposition~\ref{prop:logistic-lasso-moments}.

\subsection{Existence}
\label{app:population-existence}
\begin{proof}[Proof of existence in Theorem~\ref{thm:population-uniqueness}]
Fix $n,p$. The variables in the auxiliary minimization problem are exactly
\[
 (W,\zeta,\alpha)\in
 L^2(G;\R^p)\times L^2(x,y,Z)\times\R.
\]
Thus $W$ is a measurable function of $G$, whereas $\zeta$ is a measurable
function of $(x,y,Z)$. For each $W$, define
\[
 \mu:=\frac{\E_G\left[W\right]}{\sqrt n},\qquad
 C_W:=\E_G\left[(W-\E_G\left[W\right])
                    (W-\E_G\left[W\right])^\top\right].
\]
In particular, $\mu$ and $C_W$ are determined by $W$ and are not additional
variables. The variance constraint in
\eqref{eq:population-convex-program} forces $\alpha\ge0$. Consider
\begin{equation}
\label{eq:population-convex-program}
\begin{aligned}
 \text{minimize}\quad&
 \frac1n\E_G\left[\rho(W)-\rho(0)\right]
       +\E\left[L_y(\zeta)-L_y(0)\right]\\
 \text{subject to}\quad&
 \sqrt{\tr(\Sigma_x C_W)/n}\le\alpha,\\
 &\|x^\top\mu+\alpha Z-\zeta\|_{L^2}
       +\frac1n\E_G\left[G^\top W\right]\le0.
\end{aligned}
\end{equation}
The constraints are convex and continuous; $(W,\zeta,\alpha)=(0,0,0)$ is
feasible with value zero. This Hilbert-space formulation is related to
\citet{BellecKoriyamaExistence}. The coefficients $\beta$, $\nu$, $\kappa$ and
$\eta$ are not primal variables of \eqref{eq:population-convex-program}.

\paragraph{Attainment.}
The finite convex penalty has a unique minimum at zero, hence grows at
least linearly in $|w|$. On each feasible objective sublevel of height $M$,
nonnegative loss and the variance constraint therefore give
$\E_G|W|/\sqrt n+|\mu|\le C_M$ and
$\|W\|_{L^2}/\sqrt n\le C_M+C\alpha$.
Strict loss convexity and symmetry of $Z$ give
$\E[L_y(\alpha Z)-L_y(0)]\ge c\alpha-C$ for some $c>0$.
Put $\varepsilon_R:=\|(G/\sqrt n)\mathbf1_{\{|G|>R\sqrt n\}}\|_{L^2}\to0$ as $R\to\infty$.
Truncating the Gaussian correlation at radius $R\sqrt n$ and using
feasibility and the $L_1$-Lipschitz loss yields the coercive bound
\[
 \E[L_y(\zeta)-L_y(0)]\ge(c-C\varepsilon_R)\alpha-C_{M,R}.
\]
Choose a fixed $R$ with $C\varepsilon_R<c/2$. Since the centered loss is
bounded above on the sublevel, this bounds $\alpha$, then
$\|W\|_{L^2}$ and, by the residual constraint, $\|\zeta\|_{L^2}$.
All bounds here are at fixed size. The feasible set is weakly closed;
Fatou and convexity make the penalty integral weakly lower semicontinuous,
and the centered loss integral is convex and Lipschitz on $L^2$.
Weak compactness therefore gives a minimizer.

\paragraph{KKT multipliers and positivity.}
Let $\mathcal F(W,\zeta,\alpha)$ denote the objective in
\eqref{eq:population-convex-program}, and write its two constraint functions as
\[
 c_{\rm var}(W,\alpha)
 :=\sqrt{\tr(\Sigma_x C_W)/n}-\alpha,
 \qquad
 c_{\rm res}(W,\zeta,\alpha)
 :=\|x^\top\mu+\alpha Z-\zeta\|_{L^2}
       +\frac1n\E_G\left[G^\top W\right].
\]
Slater's condition (a strictly feasible point with finite objective) holds for
$W:=-G/(1+|G|)$, $\alpha>\sqrt{\|\Sigma_x\|_{\rm op}/n}$ and $\zeta:=\alpha Z$:
$\mu=0$ and $\E_G\left[G^\top W\right]<0$.
Thus the KKT theorem supplies a multiplier $\tau_\nu\ge0$ for
$c_{\rm var}\le0$ and a multiplier $\beta\ge0$ for $c_{\rm res}\le0$.
Concretely, the minimizing triple is stationary for the Lagrangian
\begin{equation}
\label{eq:population-existence-lagrangian}
\begin{aligned}
 \mathcal K(W,\zeta,\alpha;\tau_\nu,\beta)
 &:=\mathcal F(W,\zeta,\alpha)
   +\tau_\nu c_{\rm var}(W,\alpha)
   +\beta c_{\rm res}(W,\zeta,\alpha),\\
 0&\in\partial_{(W,\zeta,\alpha)}
       \mathcal K(W,\zeta,\alpha;\tau_\nu,\beta).
\end{aligned}
\end{equation}
Complementary slackness means
\begin{equation}
\label{eq:population-existence-complementarity}
 \tau_\nu c_{\rm var}(W,\alpha)=0,
 \qquad \beta c_{\rm res}(W,\zeta,\alpha)=0.
\end{equation}
At a minimizer, set
$U:=x^\top\mu+\alpha Z$, $\tau_\kappa:=\|U-\zeta\|_{L^2}$ and
$\xi:=L_y'(\zeta)$. Notice that $\tau_\kappa$ is the residual norm, not a
Lagrange multiplier. Choose
$s_{\rm res}\in\partial\|\cdot\|_{L^2}(U-\zeta)$; then
$\|s_{\rm res}\|_{L^2}\le1$, with
$s_{\rm res}=(U-\zeta)/\tau_\kappa$ when $\tau_\kappa>0$.
The $\zeta$- and $\alpha$-stationarity conditions for
\eqref{eq:population-existence-lagrangian} are, respectively,
\begin{equation}
\label{eq:population-existence-multipliers}
 \xi=\beta s_{\rm res},\qquad
 \tau_\nu=\beta\E\left[Zs_{\rm res}\right]
          =\E\left[Z\xi\right].
\end{equation}
In particular, $\|\xi\|_{L^2}\le\beta$, and if $\tau_\kappa>0$, then
$\xi=\beta(U-\zeta)/\tau_\kappa$.
Equation~\eqref{eq:population-existence-complementarity} gives
$c_{\rm res}(W,\zeta,\alpha)=0$ if $\beta>0$ and
$c_{\rm var}(W,\alpha)=0$ if $\tau_\nu>0$.
When $\tau_\nu=0$, stationarity in $W$ is
\begin{equation}
\label{eq:population-existence-W-kkt-zero}
 -\sqrt n\,\E\left[x\xi\right]-\beta G\in\partial\rho(W).
\end{equation}
When $\tau_\nu,\alpha>0$, it is
\begin{equation}
\label{eq:population-existence-W-kkt}
 -\frac{\tau_\nu}{\alpha}\Sigma_x(W-\E_G\left[W\right])
       -\sqrt n\,\E\left[x\xi\right]-\beta G\in\partial\rho(W).
\end{equation}
The convex sum rule applies because the constraints are continuous;
testing its integral subgradient inequality on measurable subsets gives
these pointwise inclusions.

Positivity follows before any division. If $\beta=0$, then $\xi=\tau_\nu=0$,
$W=0$, and feasibility gives $\zeta=\alpha Z$, contradicting
$\|L_y'(\alpha Z)\|_{L^2}>0$ by strict increase when $\alpha>0$ and
\eqref{eq:population-existence-conditions} when $\alpha=0$.
If $\alpha=0$, then $W=\sqrt n\mu$, $\zeta=x^\top\mu$ and $\tau_\nu=0$.
Equation~\eqref{eq:population-existence-W-kkt-zero} is impossible:
$\partial\rho(\sqrt n\mu)$ is bounded because $\rho$ is finite, whereas
$\beta G$ is nondegenerate Gaussian.
Finally, if $\tau_\kappa=0$, Gaussian integration by parts yields
$\tau_\nu=\alpha\E\left[L_y''(U)\right]>0$ and
$\E_G\left[G^\top W\right]=0$. Jensen and
\eqref{eq:population-existence-W-kkt} contradict this:
\[
 0\le\frac{\E_G\left[\rho(W)\right]-\rho(\E_G\left[W\right])}{n}
 \le-\frac{\tau_\nu}{\alpha}\frac{\tr(\Sigma_x C_W)}n
       -\frac\beta n\E_G\left[G^\top W\right]
 =-\tau_\nu\alpha<0.
\]

\paragraph{Identification.}
Now set $\kappa:=\tau_\kappa/\beta$. Equation
\eqref{eq:population-existence-multipliers} and scalar Gaussian integration
by parts give
\[
 \zeta=\prox_{\kappa L_y}(U),\qquad
 \E\left[\xi^2\right]=\beta^2,\qquad
 \tau_\nu=\alpha\E\left[\frac{L_y''(\zeta)}{1+\kappa L_y''(\zeta)}\right]
 =\alpha\nu>0.
\]
The equalities $c_{\rm var}(W,\alpha)=c_{\rm res}(W,\zeta,\alpha)=0$ give,
respectively, the variance and Gaussian-correlation equations in
\eqref{eq:fixedpoint}. With
$\eta:=\nu\Sigma_x\mu-\E\left[x\xi\right]$,
\eqref{eq:population-existence-W-kkt} identifies the unique optimizer
in \eqref{eq:effective}. All coefficients are finite. Comparing with zero gives
\begin{equation}
\label{eq:population-energy-bound}
 \frac1n\E_G\left[\rho(W)-\rho(0)\right]+\E\left[L_y(\zeta)\right]
 \le\E\left[L_y(0)\right].
\end{equation}
\end{proof}

The following example shows why the score-nondegeneracy condition in
\eqref{eq:population-existence-conditions} is needed.

\begin{proposition}[The score-nondegeneracy condition cannot be omitted]
\label{prop:logcosh-nonexistence}
Let $\rho(w):=\lambda\|w\|_1$, $\lambda>0$, and $L_y(u):=\log\cosh u$.
If $\E\left[|x|^2\right]<\infty$ and $\Sigma_x\succ0$, then $\widehat W=0$ for every
sample, but the population system has no finite solution with
$\beta,\nu>0$. Here $L_y'(0)=0$ almost surely.
\end{proposition}
\begin{proof}
Loss and penalty are nonnegative and vanish together only at zero.
At a supposed root, multiply the surrogate inclusion by $W$ and average:
\[
 0\le\frac\lambda n\E_G\left[\|W\|_1\right]
 =-\mu^\top\E\left[x\xi\right]-\nu\alpha^2+\kappa\beta^2
 =-\E\left[\zeta\tanh\zeta\right]\le0.
\]
Hence $\zeta=0$ and $\beta^2=\E\left[\xi^2\right]=0$, a contradiction.
\end{proof}

\subsection{Uniqueness}
\label{app:population-uniqueness-proof}
\begin{proof}[Proof of uniqueness in Theorem~\ref{thm:population-uniqueness}]
Under the uniqueness assumptions of Theorem~\ref{thm:population-uniqueness},
take solutions $(\mu,\alpha,\kappa,\nu,\beta,\eta)$ and
$(\mu',\alpha',\kappa',\nu',\beta',\eta')$; primes refer to the second
solution, coupled through the same $G$ and $(x,y,Z)$.
We show that both minimize \eqref{eq:population-convex-program}.
Their feasible midpoint then identifies $\zeta$;
\eqref{eq:losscoeff}--\eqref{eq:fixedpoint} then identify the remaining
coefficients.

For either solution, the surrogate inclusion gives
\[
 \nu\lambda_{\min}(\Sigma_x)|W|\le|\sqrt n\,\eta-\beta G|,\qquad
 0\le\rho(W)-\rho(0)\le(\sqrt n\,\eta-\nu\Sigma_xW-\beta G)^\top W.
\]
Thus $W\in L^2$ and $\rho(W)\in L^1$. If $\alpha=0$, the variance identity
in \eqref{eq:fixedpoint} makes $W$ deterministic. The surrogate optimality
inclusion would then
require the bounded set $\partial\rho(W)$ of the finite convex function
$\rho$ to contain the nondegenerate Gaussian term $\beta G$, which is
impossible. Jensen gives
$0\le\{\E_G\left[\rho(W)\right]-\rho(\E_G\left[W\right])\}/n
\le-\nu\alpha^2+\beta^2\kappa$, so $\kappa>0$.
Both constraints in \eqref{eq:population-convex-program} are tight,
$\beta\kappa=\|U-\zeta\|_{L^2}$,
$\xi=(U-\zeta)/\kappa$, $\alpha\nu=\E\left[Z\xi\right]$, and
\begin{equation}
\label{eq:population-optimality}
 0\in\partial\rho(W)+\nu\Sigma_x(W-\E_G\left[W\right])
       +\sqrt n\,\E\left[x\xi\right]+\beta G.
\end{equation}
Also $\zeta=U-\kappa\xi\in L^2$ and
$|L_y(\zeta)-L_y(0)|\le L_1|\zeta|$, so the centered loss is integrable.
The preceding identities are the KKT conditions of
\eqref{eq:population-convex-program}, with positive multipliers
$\alpha\nu$ and $\beta$; the surrogate inclusion supplies an $L^2$
penalty subgradient. Convex KKT sufficiency therefore makes both roots
minimizers of the same program. At their feasible midpoint the two
nonnegative Jensen gaps sum to zero; in particular,
\begin{equation}
\label{eq:population-uniqueness-jensen-gap}
 \E\left[\frac{L_y(\zeta)+L_y(\zeta')}{2}
           -L_y\!\left(\frac{\zeta+\zeta'}2\right)\right]=0
 \quad\Longrightarrow\quad \zeta=\zeta',\quad\xi=\xi',\quad\beta=\beta'.
\end{equation}
The Jensen gap in \eqref{eq:population-uniqueness-jensen-gap} is integrable
without an assumption on $L_y(0)$.
Subtracting $U-\zeta=\kappa\xi$, which follows from
\eqref{eq:losscoeff}, for the two solutions gives
\[
 x^\top(\mu-\mu')+(\alpha-\alpha')Z=(\kappa-\kappa')\xi.
\]
Conditionally on almost every $(x,y)$, continuity extends this identity to
all $Z$. Boundedness of $\xi$ forces $\alpha=\alpha'$; its strict increase
in $Z$, since $\alpha>0$, then forces $\kappa=\kappa'$.
Finally $\Sigma_x\succ0$ gives $\mu=\mu'$,
\eqref{eq:losscoeff} gives
$\nu=\nu'$ and $\eta=\eta'$, and strict convexity in \eqref{eq:effective}
gives $W=W'$ almost surely.
\end{proof}

\subsection{Stationary-point formulation}
\label{app:population-stationarity}
\begin{lemma}[Stationary-point characterization of the population system]
\label{lem:population-stationarity}
Under the uniqueness assumptions of Theorem~\ref{thm:population-uniqueness},
define, for $\mu,\eta\in\R^p$ and $\alpha,\kappa,\beta,\nu>0$,
\begin{equation}
\label{eq:population-stationary-objective}
\begin{aligned}
 \mathcal J(\mu,\alpha,\kappa;\beta,\nu,\eta)
 :={}&\E\left[e_{L_y}(x^\top\mu+\alpha Z;\kappa)-L_y(0)\right]
       +\frac{\beta^2\kappa}{2}\\
 &-\frac\nu2(\alpha^2+\mu^\top\Sigma_x\mu)+\eta^\top\mu\\
 &+\frac1n\E_G\left[\min_w\left\{
 \rho(w)+\frac\nu2w^\top\Sigma_xw
                 +(\beta G-\sqrt n\,\eta)^\top w\right\}\right].
\end{aligned}
\end{equation}
The stationary points of $\mathcal J$ are exactly the solutions of
\eqref{eq:losscoeff}--\eqref{eq:fixedpoint}.
\end{lemma}
\begin{proof}
Treat the coefficients as free variables and define $W$ by
\eqref{eq:effective} and $U,\zeta,\xi$ by the first line of \eqref{eq:losscoeff}.
The bounds $|e_{L_y}(U;\kappa)-L_y(0)|\le L_1|U|+\kappa L_1^2/2$ and
$\nu\lambda_{\min}(\Sigma_x)|W|\le|\sqrt n\,\eta-\beta G|$,
together with bounded loss derivatives, make $\mathcal J$ finite and justify
envelope differentiation under the expectations.
Stationarity in $\eta,\nu,\beta$ gives the mean, variance and
Gaussian-correlation equations in \eqref{eq:fixedpoint}, using the mean
identity to center the quadratic moment.
Stationarity in $\kappa,\mu$ gives the equations for $\beta^2,\eta$ in
\eqref{eq:losscoeff}; stationarity in $\alpha$ and Gaussian integration by
parts give
\[
 0=\partial_\alpha\mathcal J
 =\alpha\left(\E\left[\frac{L_y''(\zeta)}{1+\kappa L_y''(\zeta)}\right]-\nu\right),
\]
hence the equation for $\nu$. Conversely, the population equations make
all these derivatives vanish.
\end{proof}

\subsection{Population coefficient bounds}
\label{app:population-coefficients}

\paragraph{Population moment and Gaussian-correlation identities.}
At any finite population root satisfying
$\underline c_xI_p\preceq\Sigma_x\preceq KI_p$,
\begin{equation}
\label{eq:population-root-second-moment}
 \E_G\left[|W_\ast|^2\right]/n=|\mu_\ast|^2+\tr C_{W_\ast}/n,
 \qquad
 \underline c_x\tr C_{W_\ast}/n\le\alpha_\ast^2\le K\tr C_{W_\ast}/n,
\end{equation}
\begin{equation}
\label{eq:population-root-correlation}
 \beta_\ast\kappa_\ast
 =-\frac1n\E_G\left[G^\top(W_\ast-\sqrt n\mu_\ast)\right]
 \le\sqrt{p/n}\,\alpha_\ast.
\end{equation}
The first identities follow from the mean--covariance decomposition and
spectral bounds; the correlation bound is Cauchy--Schwarz in the
$\Sigma_x^{-1}$ and $\Sigma_x$ metrics, using
$\E_G[G^\top\Sigma_x^{-1}G]=p$.
Thus bounded normalized surrogate second moments are equivalent to bounded
$|\mu_\ast|^2+\alpha_\ast^2$.
Under the existence assumptions of Theorem~\ref{thm:population-uniqueness},
if $\E[L_y(0)]\le L_0$, \eqref{eq:population-energy-bound} also gives
\[
\begin{aligned}
 \rho(w)-\rho(0)\ge\gamma_0|w|^2/2,\quad\gamma_0>0
 &\ \Longrightarrow\ \E_G[|W_\ast|^2]/n\le2L_0/\gamma_0,\\
 \rho(w)=\lambda\|w\|_1
 &\ \Longrightarrow\ \lambda\E_G[\|W_\ast\|_1]/n\le L_0.
\end{aligned}
\]
The Lasso bound alone does not control normalized second moments.

Proposition~\ref{prop:population-nondegeneracy} gives uniform coefficient
and surrogate-moment bounds on the region
\begin{equation}
\label{eq:population-box}
 \mathcal Q:=\{q=(\mu,\alpha^2,\kappa):
 |\mu|\le M,\ 0\le\alpha^2\le\alpha_+^2,\ 0\le\kappa\le K_\kappa\},
 \qquad \alpha\ge0.
\end{equation}
Use $\|q\|_1:=|\mu|+|\alpha^2|+|\kappa|$, also for differences.
For $q\in\mathcal Q$, compute $\zeta_q,\xi_q,\nu_q,\beta_q,\eta_q$ by
\eqref{eq:losscoeff} and $W_q$ by \eqref{eq:effective}, with independent
$G\sim N(0,\Sigma_x)$. No fixed-point equation is imposed on $q$.

\begin{proposition}[Uniform coefficient and surrogate-moment bounds]
\label{prop:population-nondegeneracy}
Suppose $\underline c_xI_p\preceq\Sigma_x\preceq KI_p$, $p/n\le K$, and
Assumption~\ref{ass:analytic} holds. Uniformly for $q\in\mathcal Q$,
\[
 \nu_0\le\nu_q\le L_2,\qquad 0\le\beta_q\le L_1,\qquad
 \beta_q\ge\nu_q\alpha,
 \qquad
 |\E\left[x\xi_q\right]|+|\eta_q|\le C.
\]
The surrogate has a unique minimizer, including when $\beta_q=0$, and
$\E_G\left[|W_q|^2\right]/n\le C$.
The constants $\nu_0>0,C$ depend only on
$\underline c_x,K,L_1,L_2$, the parameters $M,\alpha_+,K_\kappa$ defining
\eqref{eq:population-box}, and the lower-curvature function $c_L(\cdot)$ in
Assumption~\ref{ass:analytic}. In particular,
$\alpha_\ast\ge\underline\alpha>0$ implies
$\beta_\ast\ge\nu_0\underline\alpha>0$ whenever $q_\ast\in\mathcal Q$.
\end{proposition}
\begin{proof}
Since $\E\left[U^2\right]\le KM^2+\alpha_+^2$, choose a fixed $B$ with
$\mathbb P(|U|\le B)\ge1/2$. The equation
$\zeta_q+\kappa L_y'(\zeta_q)=U$, obtained from \eqref{eq:losscoeff}, gives
$|\zeta_q|\le|U|+K_\kappa L_1$, hence
\[
 \nu_q\ge\frac{c_L(B+K_\kappa L_1)}{2(1+K_\kappa L_2)}=: \nu_0>0,
 \qquad \E\left[Z\xi_q\right]=\alpha\nu_q\le\beta_q.
\]
The last identity and inequality follow from Gaussian integration by
parts and Cauchy--Schwarz, also when $\alpha=0$.
For every unit vector $u$,
$|u^\top\E\left[x\xi_q\right]|\le\sqrt K L_1$, which bounds $\eta_q$.
Finally, positive curvature and monotonicity against
$0\in\partial\rho(0)$ give a unique optimizer with
$\nu_0\underline c_x|W_q|\le|\sqrt n\,\eta_q-\beta_qG|$.
Taking second moments and using $\tr(\Sigma_x)/n\le Kp/n$ proves
$\E_G\left[|W_q|^2\right]/n\le C$.
\end{proof}

\subsection{Population moment bounds for logistic Lasso}
\label{app:logistic-lasso-moments}

\begin{proof}[Proof of Proposition~\ref{prop:logistic-lasso-moments}]
Suppress stars and let $C$ denote constants independent of $n$.
By \eqref{eq:population-energy-bound},
\[
 \frac{\lambda_1}{n}\E_G\|W\|_1+\E[L_y(\zeta)]\le\log2.
\]
Projection anti-concentration, Jensen's inequality in $Z$, and Gaussian
symmetry give
$\E|U|\ge C^{-1}(\sqrt{\mu^\top\Sigma_x\mu}+\alpha)$.
Hence \eqref{eq:logistic-label-overlap} implies
\[
 \E[L_y(U)]\ge\E[(-yU)_+]\ge\epsilon\E|U|
 \ge C^{-1}\|U\|_{L^2}.
\]
Since logistic loss is $1$-Lipschitz and
$\|U-\zeta\|_{L^2}=\kappa\beta$, truncating each Gaussian coordinate
at $R>0$ and applying Cauchy--Schwarz yields
\begin{align*}
 \|U\|_{L^2}
 &\le C(1+\kappa\beta)\\
 &\le C(1+R)+C\left(\frac1n\sum_{j=1}^p
      \E[G_j^2\mathbf1_{\{|G_j|>R\}}]\right)^{1/2}\|U\|_{L^2}.
\end{align*}
Here the constants do not depend on $R$; we used
$\kappa\beta=-\E_G[G^\top W]/n$, $\lambda_1\ge\lambda_-$, and
$\E_G|W|^2/n\le\|U\|_{L^2}^2/\underline c_x$.
The tail coefficient tends to zero uniformly in $n$ because $p/n$
and the Gaussian coordinate variances are bounded.
Choosing a sufficiently large fixed $R$ and absorbing the last term gives
\[
 \|U\|_{L^2}\le C,\qquad \frac1n\E_G|W|^2\le C.
\]

Next, \eqref{eq:population-root-correlation} gives
$\kappa\beta\le\sqrt{p/n}\,\alpha\le C$, hence
$\|\zeta\|_{L^2}\le\|U\|_{L^2}+\kappa\beta\le C$.
Chebyshev's inequality places at least half the mass of $\zeta$ in a
fixed bounded interval. Since $|L_y'(\zeta)|=(1+e^{y\zeta})^{-1}$,
this gives $\beta\ge C^{-1}$ and hence $\kappa\le C$.

Finally, independent surrogate copies and their $\ell_1$ subgradients satisfy
\[
 \nu\Sigma_x(W-W')+\lambda_1(s-s')=-\beta(G-G'),\qquad
 |s_j-s_j'|\le2.
\]
The bounds $\beta\ge C^{-1}$, $\lambda_1\le\lambda_+$ and
$G_j-G_j'\sim N(0,2(\Sigma_x)_{jj})$, with
$(\Sigma_x)_{jj}\ge\underline c_x$, imply
\[
 C^{-1}\le\E\bigl[(\beta|G_j-G_j'|-2\lambda_1)_+^2\bigr]
 \le\E\bigl[|[\nu\Sigma_x(W-W')]_j|^2\bigr].
\]
Summing and using $\nu\le1/4$ gives
\[
 \frac pC\le\E|\nu\Sigma_x(W-W')|^2\le C\tr(C_W),\qquad
 \alpha^2\ge\frac{\underline c_x}{n}\tr(C_W)\ge\frac p{Cn}.
\]
Together with $p/n\to c>0$, this proves
\eqref{eq:population-moment-regime}.
\end{proof}

\section{Moment bounds, derivative identities and coupling}
\label{app:calculus}
This appendix collects three sets of tools.
Appendix~\ref{app:inputs} establishes feature moments, positive-probe
response and optimizer moment bounds
(Lemmas~\ref{lem:qf-moments}--\ref{lem:interior-moments}).
Appendix~\ref{app:transport-calculus} derives value and Hamiltonian
derivatives (Lemmas~\ref{lem:interior-envelope}--\ref{lem:derivative-calculus}).
Appendix~\ref{app:coefficients} gives coupling bounds for moments and fresh
scores (Proposition~\ref{prop:gc-transfer}), used in the endpoint arguments
of Appendices~\ref{app:proof-main-theorem}--\ref{app:prediction-laws}.

\subsection{Positive-probe response and moments}
\label{app:inputs}
For the moment bounds and the comparison in Appendix~\ref{app:pred-proof},
fix dimension-independent constants $0<v_-\le v_+<\infty$,
$0<b_-\le b_+<\infty$, and $0\le H<\infty$, and set
\begin{equation}
\label{eq:prediction-positive-domain}
 \mathcal Z_n:=[v_-,v_+]\times[b_-,b_+]
             \times\{h\in\R^p:|h|\le H\}.
\end{equation}
\[
 \mathcal D_+:=[0,1]\times\mathcal Z_n.
\]
The constants depend on the probe and field bounds, with no lower bound on
$t$ or $1-t$. The interior domain $\mathcal D$ remains in use for finite
differences. The calculus requires only $0\in\partial\rho(0)$,
not loss-score nondegeneracy.

\begin{lemma}[Feature moment bounds]
\label{lem:qf-moments}
Suppose Assumption~\ref{ass:feature-concentration} holds,
$\Sigma_x\preceq KI_p$, and $p=O(n)$. Then
\eqref{eq:feature-quadratic-concentration} holds with
$K_Q^2=4C_{\rm P}K$.
For deterministic symmetric $A$ and $u\in\R^p$,
\[
 \E_x|x^\top Ax-\tr(\Sigma_xA)|\le K_Q\|A\|_{\hs},\quad
 \E_x|u^\top x|^4\le (K_Q^2+K^2)|u|^4,\quad
 \E_x|x|^4\le (K_Q^2+K^2)p^2.
\]
Consequently,
$\sup_n\E_x\left[(1+|x|^2/n)\mathbf1_{\{|x|>R\sqrt n\}}\right]\to0$
as $R\to\infty$.
\end{lemma}
\begin{proof}
Apply \eqref{eq:data-poincare} to smooth bounded truncations of
$x^\top Ax$ (whose label energy is zero), then pass to the limit by
dominated convergence of the means and Fatou's lemma for centered squares:
\[
 \Var_x(x^\top Ax)\le4C_{\rm P}\E_x|Ax|^2
 =4C_{\rm P}\tr(A^2\Sigma_x)\le4C_{\rm P}K\|A\|_{\hs}^2.
\]
This proves \eqref{eq:feature-quadratic-concentration} without assuming
fourth moments.
Cauchy--Schwarz applied to \eqref{eq:feature-quadratic-concentration} gives
$\E_x\left[|x^\top Ax-\tr(\Sigma_xA)|\right]\le K_Q\|A\|_{\hs}$.
The choice $A=uu^\top$ gives
$\E_x\left[|u^\top x|^4\right]\le K_Q^2|u|^4+(u^\top\Sigma_xu)^2$.
Summing coordinate fourth moments by Cauchy--Schwarz gives
$\E_x\left[|x|^4\right]\le (K_Q^2+K^2)p^2$;
Markov's inequality bounds the tail by $C(R^{-4}+R^{-2})$.
\end{proof}

\begin{lemma}[Well-posedness and response under positive curvature]
\label{lem:interior-wellposedness}
Suppose Assumption~\ref{ass:analytic} holds and
$\Sigma_x\succeq\underline c_xI_p$. Fix observations and selectors and put
\[
 f_v(w):=\rho(w)+\sum_iB_iL_{y_i}(x_i^\top w/\sqrt n)
                       +\frac v2w^\top\Sigma_xw,
 \qquad \ell:=\sqrt n\,h-\sqrt{2b}\,G.
\]
For $v>0$, $W(\ell):=\argmin_w\{f_v(w)-\ell^\top w\}$ exists uniquely and
\begin{equation}
\label{eq:monotone-lip}
 v\underline c_x|W(\ell)-W(\ell')|^2
 \le(\ell-\ell')^\top(W(\ell)-W(\ell')),
 \qquad |W(\ell)-W(\ell')|\le\frac{|\ell-\ell'|}{v\underline c_x}.
\end{equation}
Its almost-everywhere normalized field response
$J:=D_hW/\sqrt n=D_\ell W$ is symmetric with
$0\preceq J\preceq(v\underline c_x)^{-1}I_p$; for $b>0$ it exists at the
Gaussian-perturbed field almost surely.
\end{lemma}
\begin{proof}
The positive quadratic term makes every linear tilt coercive (tending to
$+\infty$ as $|w|\to\infty$) and $v\underline c_x$-strongly convex. Thus the minimizer exists uniquely,
and strong monotonicity of $\ell\in\partial f_v(W(\ell))$ gives
\eqref{eq:monotone-lip}. Since $W=\nabla f_v^*$, its almost-everywhere
Jacobian satisfies $0\preceq D_\ell W\preceq(v\underline c_x)^{-1}I_p$.
The nondegenerate Gaussian field avoids its null exceptional set;
Fubini also justifies differentiation along almost every scalar force.
\end{proof}

\begin{lemma}[Moment bounds on positive-probe domains]
\label{lem:interior-moments}
Let $\mathcal D_+$ be the positive-probe domain above, $W:=W_{t,z}$ the optimizer in
\eqref{eq:interior-estimator}, and $W_-$ its unforced cavity obtained by
omitting any fixed number of observation terms from \eqref{eq:energy},
with the remaining selectors, Gaussian probe, and normalization $n$ unchanged.
Assume that the training pairs are independent copies of $(x,y)$.
Suppose Assumption~\ref{ass:analytic} holds,
$\underline c_xI_p\preceq\Sigma_x\preceq KI_p$, and $p/n\le K$.
\begin{enumerate}[label=(\roman*),ref=\roman*]
\item \label{item:interior-second-moments}We have
\begin{equation}
\label{eq:moment}
 \sup_{(t,z)\in\mathcal D_+}
 \left(\E\left[|W|^2\right]+\E\left[|W_-|^2\right]\right)
 \le C_{\mathcal D_+}n.
\end{equation}
\item \label{item:interior-fourth-moments}If, in addition,
Assumption~\ref{ass:feature-concentration} holds and
$p/n\to c\in(0,\infty)$, then
\begin{equation}
\label{eq:interior-fourth-moments}
 \sup_{(t,z)\in\mathcal D_+}
 \left(\E\left[|W|^4\right]+\E\left[|W_-|^4\right]\right)
 \le C_{\mathcal D_+}n^2.
\end{equation}
\end{enumerate}
Expectations are over the training sample, selectors, and probe.
We write $\|U\|_{L^2}:=(\E|U|^2)^{1/2}$.
Both bounds are uniform in the omitted indices and in $t\in[0,1]$. In \textup{(i)}, $C_{\mathcal D_+}$
depends only on $\underline c_x,K,L_1$ and the bounds defining $\mathcal Z_n$; in \textup{(ii)},
it may also depend on $C_{\rm P}$.
\end{lemma}
\begin{proof}
Put $z_i:=B_iL_{y_i}'(0)x_i$. Monotonicity against zero gives
\begin{equation}
\label{eq:interior-zero-comparison}
 v\underline c_x|W|
 \le\left|n^{-1/2}\sum_i z_i+\sqrt{2b}\,G-\sqrt n\,h\right|,
\end{equation}
with the corresponding summands omitted for $W_-$.
Since $\E\left[|z_i|^2\right]\le L_1^2\tr(\Sigma_x)\le Cn$ and
$|\E\left[z_i\right]|\le L_1\sqrt K$, independence gives
$\E\left[|n^{-1/2}\sum_i z_i|^2\right]\le Cn$.
Together with $\E\left[|G|^2\right]=\tr(\Sigma_x)=O(n)$ and bounded $|h|$,
\eqref{eq:interior-zero-comparison} proves \textup{(i)}.

Under the additional assumptions of \textup{(ii)},
Lemma~\ref{lem:qf-moments} gives $\E\left[|z_i|^4\right]\le Cn^2$.
The fourth-moment inequality for independent centered vectors yields
\[
 \E\left[\left|n^{-1/2}\sum_i(z_i-\E\left[z_i\right])\right|^4\right]
 \le\frac C{n^2}\left\{
       \sum_i\E\left[|z_i|^4\right]
       +\left(\sum_i\E\left[|z_i|^2\right]\right)^2\right\}
 \le Cn^2.
\]
Restoring the means preserves this order. Since
$\E\left[|G|^4\right]=O(n^2)$ and $n^2|h|^4=O(n^2)$,
\eqref{eq:interior-zero-comparison} also proves \textup{(ii)}.
Independence and $0\le B_i(t)\le1$ give the same constants for every
$t\in[0,1]$.
\end{proof}

\subsection{Value and Hamiltonian derivatives}
\label{app:transport-calculus}
The envelope identities turn value gradients into moments; scalar Moreau
derivatives give the corresponding Hamiltonian coefficients.
Recall $z=(v,b,h)$ from Section~\ref{sec:perturbed-optimizers},
and write
\begin{equation}
\label{eq:normalized-coordinates}
 \pi:=\nabla_z\phi=(\partial_v\phi,\partial_b\phi,\nabla_h\phi).
\end{equation}

\begin{lemma}[Value derivatives and moment identities]
\label{lem:interior-envelope}
\label{lem:normalized-identities}
Under the assumptions of Lemma~\ref{lem:interior-moments}\textup{(\ref{item:interior-second-moments})},
writing $\pi=\nabla_z\phi$, the moment dictionary is
\begin{equation}
\label{eq:envelope}
\begin{aligned}
 \mu&=-\pi_h=-\nabla_h\phi,\qquad
 \alpha^2=2\pi_v-\pi_h^\top\Sigma_x\pi_h
          =2\partial_v\phi-\mu^\top\Sigma_x\mu,\\
 \kappa&=-\pi_b=-\partial_b\phi
          =-\frac{\E[G^\top W]}{n\sqrt{2b}}.
\end{aligned}
\end{equation}
Moreover, $\partial_t\phi=\E[\mathcal I]$, where $\mathcal I$ is the
minimum-value increment from reinserting one independent observation
into its unforced cavity.
\end{lemma}
\begin{proof}
The unique-minimizer envelope theorem gives
\[
 \nabla_h\phi=-\frac{\E\left[W\right]}{\sqrt n},\qquad
 \partial_v\phi=\frac{\E\left[W^\top\Sigma_xW\right]}{2n},\qquad
 \partial_b\phi=\frac{\E\left[G^\top W\right]}{n\sqrt{2b}}.
\]
Estimate~\eqref{eq:moment} provides local domination. Since
$D_GW=-\sqrt{2b}J$, Gaussian integration by parts gives
$\E\left[G^\top W\right]=-\sqrt{2b}\E\left[\tr(\Sigma_xJ)\right]$.
Centering the quadratic moment proves
\eqref{eq:envelope}.

For $F(B):=\min_w\cE(w)$, differentiating the finite product Bernoulli law gives
\[
 \partial_t\phi=\frac1n\sum_i
 \E\left[F(B^{i,1})-F(B^{i,0})\right]=\E\left[\mathcal I\right],
\]
where $B^{i,j}$ fixes selector $i$ to $j$. Estimate~\eqref{eq:moment},
together with $\cE(W)\le\cE(0)$, gives integrability, and exchangeability
gives the last equality.
\end{proof}

\begin{lemma}[Moreau and Hamiltonian derivatives]
\label{lem:derivative-calculus}
Under the loss conditions of Assumption~\ref{ass:analytic}, the following hold.
\begin{enumerate}[label=(\roman*),ref=\roman*]
\item \label{item:scalar-derivative-calculus}Put
$\zeta:=\prox_{\kappa L_y}(u)$ and $\xi:=L_y'(\zeta)$. For $\kappa\ge0$,
\begin{equation}
\label{eq:moreau-derivs}
 \partial_u e_{L_y}=\xi,\qquad
 \partial_{uu}e_{L_y}=\frac{L_y''(\zeta)}{1+\kappa L_y''(\zeta)},\qquad
 \partial_\kappa e_{L_y}=-\frac{\xi^2}{2},
\end{equation}
with respective absolute bounds $L_1,L_2,L_1^2/2$.
For $f\in C^2(\R)$ with bounded first and second derivatives and
$Z\sim N(0,1)$,
\begin{equation}
\label{eq:gaussian-variance-derivative}
 \partial_{\alpha^2}\E_Z\left[f(z_0+\alpha Z)\right]
 =\frac12\E_Z\left[f''(z_0+\alpha Z)\right].
\end{equation}
In particular,
\begin{equation}
\label{eq:heat-bound}
 \left|\E_Z\left[e_{L_y}(z_0+\alpha_1Z;\kappa)
                    -e_{L_y}(z_0+\alpha_2Z;\kappa)\right]\right|
 \le\frac{L_2}{2}|\alpha_1^2-\alpha_2^2|.
\end{equation}
Derivatives in $\kappa$ and $\alpha^2$ at zero are right derivatives.
\item \label{item:hamiltonian-derivative-calculus}If also
$\E\left[|x|^2\right]<\infty$, the Hamiltonians
\eqref{eq:H}--\eqref{eq:Hhat} satisfy \eqref{eq:Hgradient}, with coefficients
from \eqref{eq:losscoeff} and continuous extension to $\alpha=0$ and
$\kappa=0$.
\end{enumerate}
\end{lemma}
\begin{proof}
Implicit differentiation of $\zeta+\kappa L_y'(\zeta)=u$ and the envelope
theorem give \eqref{eq:moreau-derivs}. For $\alpha>0$, differentiation and
Gaussian integration by parts give
$\partial_{\alpha^2}\E_Z[f(z_0+\alpha Z)]
=(2\alpha)^{-1}\E_Z[Zf'(z_0+\alpha Z)]
=\frac12\E_Z[f''(z_0+\alpha Z)]$,
proving \eqref{eq:gaussian-variance-derivative}.
For $\kappa\ge0$,
$0\le e_{L_y}(u;\kappa)\le L_y(u)\le L_0+L_1|u|$,
so $\cH$ is finite. The derivative bounds and $\E|x|^2<\infty$
justify differentiation under the expectation, using Gaussian averaging
for the variance derivative. Conditioning on $(x,y)$, these identities and \eqref{eq:losscoeff} yield
\[
 \dd\cH=\E\left[x\xi\right]^\top\dd\mu
              +\frac\nu2\dd(\alpha^2)-\frac{\beta^2}{2}\dd\kappa.
\]
Under the change of variables in \eqref{eq:Hhat},
$\dd\mu=-\dd\pi_h$, $\dd\kappa=-\dd\pi_b$, and
$\dd(\alpha^2)=2\dd\pi_v+2\mu^\top\Sigma_x\dd\pi_h$.
The chain rule therefore gives
\[
 \dd\widetilde\cH=\nu\dd\pi_v+\frac{\beta^2}{2}\dd\pi_b
       +(\nu\Sigma_x\mu-\E\left[x\xi\right])^\top\dd\pi_h,
\]
so $\nabla_\pi\widetilde\cH=(\nu,\beta^2/2,\eta)$.
Gaussian integration by parts also gives $\E[Z\xi]=\alpha\nu$.
The derivative bounds also give continuous extension to the boundary,
with right derivatives at zero.
Finally, integrating \eqref{eq:gaussian-variance-derivative} for
$f(u)=e_{L_y}(u;\kappa)$ and using $|\partial_{uu}e_{L_y}|\le L_2$
proves \eqref{eq:heat-bound}.
\end{proof}

\subsection{Coupling bounds}
\label{app:coefficients}
The bounds below turn an $L^2$ coupling of optimizers into comparisons
of their moments and fresh-score laws.
For $U$ independent of a fresh pair $(x,y)$, define its prediction error by
\[
 \mathcal G(U):=\E_x\left[\Wass\!\left(
 \Law\!\left(x^\top U/\sqrt n\mid x\right),
 N\!\left(x^\top\E\left[U\right]/\sqrt n,\tr(\Sigma_xC_U)/n\right)
 \right)\right].
\]
\begin{proposition}[Coupling bounds for moments and fresh scores]
\label{prop:gc-transfer}
Assume $\Sigma_x\preceq KI_p$. For coupled square-integrable $U,W$ jointly
independent of $(x,y)$,
\begin{equation}
\label{eq:gc-coupling-transfer}
 |\mathcal G(W)-\mathcal G(U)|
 \le\left(2+\sqrt{2/\pi}\right)\sqrt{\|\Sigma_x\|_{\op}}
       \frac{\|W-U\|_{L^2}}{\sqrt n}.
\end{equation}
Writing $d:=\|W-U\|_{L^2}/\sqrt n$, one also has
\begin{equation}
\label{eq:coupling-moments}
 \begin{gathered}
 |\E\left[W-U\right]|/\sqrt n\le d,\qquad
 \E_x\left[\Wass\!\left(\Law(x^\top W/\sqrt n\mid x),
                  \Law(x^\top U/\sqrt n\mid x)\right)\right]\le\sqrt Kd,\\
 \frac{\|C_W-C_U\|_{\hs}+|\tr(\Sigma_x[C_W-C_U])|}{n}
 \le(1+K)d\frac{\|W\|_{L^2}+\|U\|_{L^2}}{\sqrt n}.
 \end{gathered}
\end{equation}
\end{proposition}
\begin{proof}
Couple both scores using $(U,W)$ and both comparison Gaussians using one
standard normal. The mean and standard-deviation differences are bounded by
\begin{equation}
\label{eq:coupling-mean-standard-deviation}
 \frac{|\E\left[W-U\right]|}{\sqrt n}\le d,\qquad
 \left|\sqrt{\tr(\Sigma_xC_W)/n}-\sqrt{\tr(\Sigma_xC_U)/n}\right|
 \le\sqrt{\|\Sigma_x\|_{\op}}d.
\end{equation}
The standard-deviation bound in \eqref{eq:coupling-mean-standard-deviation}
follows from the reverse triangle inequality in centered weighted $L^2$.
Averaging the score difference and applying the Wasserstein triangle inequality
gives \eqref{eq:gc-coupling-transfer}.
Expanding the difference of the centered outer products and using
$\|ab^\top\|_{\hs}=|a||b|$ and Cauchy--Schwarz gives
\eqref{eq:coupling-moments}.
\end{proof}

\section{Positive-probe prediction and the strongly convex full path}
\label{app:pred-proof}

We first prove the quantitative positive-probe comparison in
Proposition~\ref{prop:prediction-positive-probes}, using the response and
moment bounds of Appendix~\ref{app:calculus}. We then prove
Theorem~\ref{thm:prediction-directions} for uniformly strongly convex
penalties; its pure-Lasso case is completed in
Appendix~\ref{app:proof-main-theorem}.
The positive-probe argument uses only the marginal Poincar\'e inequality
for $x$, which follows from Assumption~\ref{ass:feature-concentration}.
No penalty/design alternative
or population moment regime is required for the proposition.
All its estimates are uniform on $\mathcal D_+=[0,1]\times\mathcal Z_n$
and over the full and unforced cavity families. Constants may depend on
the fixed positive lower bounds for $v$ and $b$.
The bounds from Appendix~\ref{app:calculus} are uniform in $t\in[0,1]$.
The selectors enter the argument below only through independence and
$0\le B_i(t)\le1$; no differentiation in $t$ is used.
Set $\Delta:=W-\E\left[W\right]$, holding the deterministic vector $\E\left[W\right]$ fixed when
deleting or reinserting an observation. Matrix suprema below are over
deterministic matrices and are outside the expectation.

\subsection{Concentration of the optimizer}
\label{app:optimizer-concentration}

\begin{lemma}[Covariance and quadratic forms]
\label{lem:pred-proof-concentration}
Under Assumptions~\ref{ass:analytic} and~\ref{ass:feature-concentration},
the spectral bounds \eqref{eq:spectral-bounds}, and $p/n\to c\in(0,\infty)$,
uniformly on
$\mathcal D_+$ and for the full and unforced cavity optimizers,
\begin{equation}
\label{eq:pred-proof-concentration}
 \|C_W\|_{\op}\le C_{\mathcal D_+},
 \qquad
 \sup_{\substack{A=A^\top\\\|A\|_{\op}\le1}}
 \left\|
   \frac{\Delta^\top A\Delta-\tr(AC_W)}n
 \right\|_{L^2}
 \le C_{\mathcal D_+}n^{-1/2}.
\end{equation}
In particular, the normalized covariance satisfies
\begin{equation}
\label{eq:interior-normalized-covariance}
 \sup_{\mathcal D_+}\frac{\|C_W\|_{\op}}n\longrightarrow0,\qquad
 \sup_{\mathcal D_+}\frac{\|C_W\|_{\hs}}n\longrightarrow0.
\end{equation}
\end{lemma}

\begin{proof}
Choose a fixed $R$ so that $R^2n\ge2\tr\Sigma_x$ for all sufficiently
large $n$, and put $A_R:=\{|x|\le R\sqrt n\}$.
On $A_R^c$, $|x|^2\le
4(|x|^2-\E\left[|x|^2\right])^2/(R^2n)$.
Thus \eqref{eq:feature-quadratic-concentration}, applied with $A=I_p$,
gives the sharper fixed-radius tail estimate
\begin{equation}
\label{eq:pred-proof-tail}
 \E\!\left[|x|^2\ind_{A_R^c}\right]
 \le \frac{4\Var(|x|^2)}{R^2n}\le C.
\end{equation}

For $\sigma:=\sqrt{2b}$, adding a scalar force
$-rx^\top w/\sqrt n$ translates $G$ to $G-rx/(\sigma\sqrt n)$.
The density of this translated Gaussian relative to the original probe is
\[
 \exp\!\left(-\frac{rx^\top\Sigma_x^{-1}G}{\sigma\sqrt n}
             -\frac{r^2x^\top\Sigma_x^{-1}x}{2\sigma^2n}\right).
\]
On $A_R$ its conditional $L^2$ norm is bounded uniformly for
$|r|\le L_1$. Conditional Cauchy--Schwarz, then averaging, gives for
every nonnegative measurable $H$
\begin{equation}
\label{eq:unforced-shift-bound}
 \E\left[\ind_{A_R}H\!\left(G-\frac{rx}{\sigma\sqrt n}\right)\right]
 \le C_{\mathcal D_+,R}
       \left(\E\left[\ind_{A_R}H(G)^2\right]\right)^{1/2}.
\end{equation}
Here $H$ may also depend on the independent fresh feature and cavity data,
which are held fixed during the Gaussian change of measure.

Fix a deterministic index not omitted from the objective defining $W$,
and denote its feature by $x$.
Let $W_-(r)$ be the optimizer with this term removed and with the additional
linear term $-rx^\top w/\sqrt n$.
Write $J_-(r)$ for its response to the physical linear field and put
$\Delta_r:=W_-(r)-\E\left[W\right]$; the deterministic vector $\E\left[W\right]$ is held fixed
throughout this proof.
Lemmas~\ref{lem:interior-wellposedness} and
\ref{lem:interior-moments}\textup{(\ref{item:interior-fourth-moments})}
give, for $|r|\le L_1$,
\[
 \|J_-(r)\|_{\op}\le C_{\mathcal D_+},\qquad
 |\Delta_r|\le |\Delta_0|+C_{\mathcal D_+}|x|/\sqrt n,\qquad
 \E\left[|\Delta_0|^4\right]\le C_{\mathcal D_+}n^2.
\]
At zero force, $(\Delta_0,J_-(0))$ is independent of the deleted observation.
For deterministic $|e|=1$ and $A=A^\top$ with $\|A\|_{\op}\le1$,
the Gaussian shift \eqref{eq:unforced-shift-bound} and
Lemma~\ref{lem:qf-moments} therefore yield
\begin{align*}
 \E\!\left[\ind_{A_R}|e^\top J_-(r)x|^2\right]
 &\le C_{\mathcal D_+}
       \left(\E\left[|e^\top J_-(0)x|^4\right]\right)^{1/2}
 \le C_{\mathcal D_+},\\
 \E\!\left[\ind_{A_R}|\Delta_r^\top A J_-(r)x|^2\right]
 &\le C_{\mathcal D_+}
       \left(\E\left[|\Delta_0^\top A J_-(0)x|^4\right]\right)^{1/2}
 \le C_{\mathcal D_+}n.
\end{align*}
Here the shift replaces both $\Delta_0(G)$ and $J_-(0;G)$ by their
forced versions, since $G$ is replaced by $G-rx/\sqrt{2bn}$.
On $A_R^c$, the response bound and \eqref{eq:pred-proof-tail} give the
first estimate with the same order. For the second, independence at
zero force gives
\[
 \E\!\left[\ind_{A_R^c}|\Delta_r^\top A J_-(r)x|^2\right]
 \le C_{\mathcal D_+}\left\{
     \E\left[|\Delta_0|^2\right]\,\E\left[|x|^2\ind_{A_R^c}\right]
     +\frac{\E\left[|x|^4\right]}{n}\right\}
 \le C_{\mathcal D_+}n.
\]
Consequently, uniformly for $|r|\le L_1$,
\begin{equation}
\label{eq:pred-proof-force}
 \E\left[|e^\top J_-(r)x|^2\right]\le C_{\mathcal D_+},
 \qquad
 \E\left[|\Delta_r^\top A J_-(r)x|^2\right]\le C_{\mathcal D_+}n.
\end{equation}

Set $F_A(w):=(w-\E\left[W\right])^\top A(w-\E\left[W\right])/n$.
Insertion of observation $i$ gives
$W=W_-(r_i)$ with
$r_i=-B_iL_{y_i}'(x_i^\top W/\sqrt n)\in[-L_1,L_1]$.
The absolute continuity and Gaussian-field argument in
Lemma~\ref{lem:interior-wellposedness} give, almost everywhere in $r$,
\[
 \frac{\dd}{\dd r}e^\top W_-(r)
    =\frac{e^\top J_-(r)x}{\sqrt n},
 \qquad
 \frac{\dd}{\dd r}F_A(W_-(r))
    =\frac{2\Delta_r^\top A J_-(r)x}{n\sqrt n}.
\]
Integrate over the force interval and bound it by $[-L_1,L_1]$
before taking expectations. By Cauchy--Schwarz and
\eqref{eq:pred-proof-force},
\[
 \E\left[|e^\top(W-W_-(0))|^2\right]\le C_{\mathcal D_+}/n,\qquad
 \E\left[|F_A(W)-F_A(W_-(0))|^2\right]\le C_{\mathcal D_+}/n^2.
\]
Resampling the entire block $(B_i,x_i,y_i)$ uses two such insertions
into a common cavity with the same Gaussian probe.

Apply Efron--Stein
to the independent observation blocks and Gaussian Poincar\'e to the
probe. The covariance-weighted probe gradients of $e^\top W$ and
$F_A(W)$ are respectively
\[
 -\sigma\Sigma_x^{1/2}Je,\qquad
 -\frac{2\sigma}{n}\Sigma_x^{1/2}JA\Delta.
\]
Their expected squared norms are at most $C_{\mathcal D_+}$ and
$C_{\mathcal D_+}/n$.
The same orders follow from summing the resampling bounds, so
$\Var(e^\top W)\le C_{\mathcal D_+}$ and
$\Var(F_A(W))\le C_{\mathcal D_+}/n$.
Since $\E\left[F_A(W)\right]=\tr(AC_W)/n$, this proves
\eqref{eq:pred-proof-concentration}.
For a cavity optimizer the common background omits two observations;
Lemma~\ref{lem:interior-moments}\textup{(\ref{item:interior-fourth-moments})} covers this case with the
original normalization $n$. Finally,
$\|C_W\|_{\hs}\le\sqrt p\,\|C_W\|_{\op}$ proves
\eqref{eq:interior-normalized-covariance}.
\end{proof}

The Gaussian translation above uses the positive probe $b$, and the
response bound uses $v$; the insertion force stays in $[-L_1,L_1]$.
Efron--Stein applies to the independent blocks $(B_i,x_i,y_i)$ even when
the selectors are degenerate. At $t=0$ the actual insertion increments
vanish and full and cavity optimizers coincide; at $t=1$ block resampling
simply resamples the observations. Thus these concentration estimates
hold on all of $\mathcal D_+$.

\subsection{Two tools for Gaussian approximation}
\label{app:gaussian-approximation-tools}
The next two lemmas isolate the spectral-gap construction and the scalar
Stein equation used below.

\begin{lemma}[Spectral-gap construction of a Stein kernel]
\label{lem:pred-proof-kernel}
Let a centered random vector $\bar x\in\R^p$ have finite second moment,
covariance $C_x$, and satisfy
\[
 \Var(f(\bar x))\le C_{\rm P}\E\left[|\nabla f(\bar x)|^2\right]
\]
for bounded $C^1$ functions with bounded gradient. There exists a measurable
matrix $\mathcal T_x(\bar x)$ with square-integrable entries such that
\begin{equation}
\label{eq:pred-proof-kernel}
\begin{gathered}
 \E\left[\bar x f(\bar x)\right]=\E\left[\mathcal T_x(\bar x)\nabla f(\bar x)\right],\qquad
 \E\left[\mathcal T_x\right]=C_x,\\
 \E\left[\|\mathcal T_x\|_{\hs}^2\right]\le C_{\rm P}\tr C_x.
\end{gathered}
\end{equation}
The identity extends to $C^1$ tests with bounded gradient. This is the
spectral-gap construction of
\citet[Theorem~2.4, equation~(6)]{CourtadeFathiPananjady2019}, with the
density-free formulation proved here.
\end{lemma}
\begin{proof}
Let $\mathcal G$ be the closure in
$L^2(\Law(\bar x);\R^p)$ of gradients of bounded $C^1$ functions with
bounded gradient. Centering and the Poincar\'e inequality imply, for each
coordinate $i$,
\[
 |\E\left[\bar x_i f(\bar x)\right]|
 \le (C_{\rm P}(C_x)_{ii})^{1/2}\|\nabla f\|_{L^2}.
\]
Thus $\nabla f\mapsto\E\left[\bar x_i f(\bar x)\right]$ is well defined and extends
continuously to $\mathcal G$. Its Riesz representative is the $i$th row of
$\mathcal T_x$. Summing the squared row-norm bounds proves the last bound
in \eqref{eq:pred-proof-kernel}. Smooth truncations of tests with bounded
gradient converge in $L^2$, and their gradients converge in $L^2$ by
dominated convergence, giving the stated extension. Applying it to the
coordinate functions yields $\E\left[\mathcal T_x\right]=C_x$.
\end{proof}

\begin{lemma}[Stein equation and Wasserstein comparison]
\label{lem:pred-proof-stein-equation}
Let $Z\sim N(0,1)$ and $a>0$. For every smooth $1$-Lipschitz function
$F:\R\to\R$, the Stein equation
\begin{equation}
\label{eq:pred-proof-stein-equation}
 a^2 f'(u)-uf(u)=F(u)-\E\left[F(aZ)\right]
\end{equation}
has a solution with
\[
 \|f'\|_\infty\le\sqrt{2/\pi}\,a^{-1},\qquad
 \|f''\|_\infty\le2a^{-2}.
\]
Consequently, if a centered scalar $T$ has $\E\left[T^2\right]=a^2$ and
\[
 |\E\left[Tf(T)\right]-a^2\E\left[f'(T)\right]|
 \le A_{\rm St}\|f'\|_\infty+B_{\rm St}\|f''\|_\infty
\]
for all $C^2$ functions with bounded first and second derivatives, then
\begin{equation}
\label{eq:pred-proof-stein-criterion}
 \Wass(\Law(T),N(0,a^2))
 \le\sqrt{2/\pi}\,\frac{A_{\rm St}}a+
                          2\frac{B_{\rm St}}{a^2}.
\end{equation}
For every $a\ge0$, the direct bound
$\Wass(\Law(T),N(0,a^2))\le2a$ also holds.
\end{lemma}
\begin{proof}
The standard normal Stein-equation bounds are given in
\citet[Lemma~4.2]{Chatterjee2008Normal}. Apply them to $u\mapsto F(au)$,
whose Lipschitz constant is at most $a$, and rescale its solution by
$f(u):=a^{-1}f_0(u/a)$. This proves the displayed bounds and
\eqref{eq:pred-proof-stein-equation}. Insert this solution into the
discrepancy bound; Gaussian mollification of $1$-Lipschitz tests and
Wasserstein duality give \eqref{eq:pred-proof-stein-criterion}.
Finally, $\E\left[|T|\right]+a\E\left[|Z|\right]\le2a$ bounds the distance by an independent
coupling, including $a=0$.
\end{proof}

\subsection{Quantitative conditional Gaussian approximation}
\label{app:gaussian-approximation-proof}

\begin{proposition}[Conditional Gaussian prediction with positive probes]
\label{prop:prediction-positive-probes}
Along a sequence in which the training pairs are independent copies of
$(x,y)$, $p/n\to c\in(0,\infty)$,
and $\underline c_xI_p\preceq\Sigma_x\preceq KI_p$, suppose
Assumptions~\ref{ass:analytic} and~\ref{ass:feature-concentration} hold.
Let $x$ be an independent fresh feature, and let $\mathcal Z_n$ be any
fixed positive-probe set from \eqref{eq:prediction-positive-domain}.
Let $W$ be either $W_{t,z}$ or its unforced cavity optimizer, with its
corresponding moments $\mu,\alpha^2$ from \eqref{eq:target-moments}.
Uniformly over the omitted index,
\begin{equation}
\label{eq:prediction-directions}
 \sup_{\substack{t\in[0,1]\\z\in\mathcal Z_n}}
 \E_x\left[\Wass\!\left(
   \Law\!\left(\frac{x^\top W}{\sqrt n}\,\middle|\,x\right),
   N(x^\top\mu,\alpha^2)\right)\right]\le Cn^{-1/8}.
\end{equation}
The constant depends only on the fixed parameter bounds and the uniform
constants in the assumptions.
\end{proposition}
\begin{proof}
Fix $(t,z)$ and either optimizer family.
Recall $m_x$, $\bar x$, and $C_x$ from Section~\ref{sec:intro}, and put
\[
 C_x=\E\left[\bar x\bar x^\top\right],\qquad
 T:=\frac{\bar x^\top\Delta}{\sqrt n},\qquad
 a^2:=\E\left[T^2\right]=\frac{\tr(C_x C_W)}n.
\]
The vectors $\bar x$ and $\Delta$ are independent and centered.
Moreover, $\|C_x\|_{\op}\le K$, $|m_x|^2\le K$, and
$a^2\le C_{\mathcal D_+}$ by
Lemma~\ref{lem:pred-proof-concentration}.

By Lemma~\ref{lem:pred-proof-kernel}, the centered feature has a Stein
kernel $\mathcal T_x(\bar x)$ satisfying \eqref{eq:pred-proof-kernel}.
Only quadratic forms of $\mathcal T_x$ are needed. Set
$\mathcal S_x:=(\mathcal T_x+\mathcal T_x^\top)/2$, so that
$\Delta^\top\mathcal T_x\Delta=\Delta^\top \mathcal S_x\Delta$,
$\E\left[\mathcal S_x\right]=C_x$, and $\E\left[\|\mathcal S_x\|_{\hs}^2\right]\le C_{\rm P}\tr C_x$.
Define
\[
 Q:=\frac{\Delta^\top \mathcal S_x(\bar x)\Delta}{n},
 \qquad H_x(\bar x):=\frac{\tr(C_W\mathcal S_x(\bar x))}{n}.
\]
For $M_{\rm sp}\ge1$, split $\mathcal S_x$ into the eigenvalues of magnitude at most $M_{\rm sp}$
and the remaining eigenvalues. Conditional on $\bar x$, the first part
has operator norm at most $M_{\rm sp}$ and is independent of $\Delta$.
Lemma~\ref{lem:pred-proof-concentration} applies to it, while the
second part is bounded using
$|\lambda|\ind_{\{|\lambda|>M_{\rm sp}\}}\le\lambda^2/M_{\rm sp}$.
Hence
\begin{equation}
\label{eq:pred-proof-kernel-error}
 \E\left[|Q-H_x|\right]
 \le C_{\mathcal D_+}\frac{M_{\rm sp}}{\sqrt n}
    +\frac{2\|C_W\|_{\op}}{nM_{\rm sp}}\E\left[\|\mathcal S_x\|_{\hs}^2\right]
 \le C_{\mathcal D_+}\left(\frac{M_{\rm sp}}{\sqrt n}+\frac1{M_{\rm sp}}\right).
\end{equation}
Taking $M_{\rm sp}:=n^{1/4}$ gives $\E\left[|Q-H_x|\right]\le C_{\mathcal D_+}n^{-1/4}$.
Also, $\E\left[H_x\right]=a^2$ and
\[
 \Var(H_x)
 \le \frac{\|C_W\|_{\hs}^2}{n^2}\E\left[\|\mathcal S_x-C_x\|_{\hs}^2\right]
 \le C_{\mathcal D_+}.
\]

\paragraph{Cancellation of the trace term.}
Let $f\in C^2(\R)$ have bounded first and second derivatives, and set
$\psi_f(\bar x):=\E_{\Delta}\left[f'(\bar x^\top\Delta/\sqrt n)\right]$.
Differentiation under the expectation and Cauchy--Schwarz in each
unit direction imply
\[
 \nabla\psi_f(\bar x)
 =\frac1{\sqrt n}\E_{\Delta}\left[
       \Delta f''(\bar x^\top\Delta/\sqrt n)\right],
 \qquad
 |\nabla\psi_f(\bar x)|^2
 \le\frac{\|C_W\|_{\op}}n\|f''\|_\infty^2.
\]
The marginal Poincar\'e inequality therefore gives
$\Var(\psi_f)\le C_{\mathcal D_+}\|f''\|_\infty^2/n$.
Since $\E\left[H_x\right]=a^2$,
\[
 |\E\left[(H_x-a^2)f'(T)\right]|
 =|\Cov(H_x,\psi_f)|
 \le C_{\mathcal D_+}n^{-1/2}\|f''\|_\infty.
\]
Applying \eqref{eq:pred-proof-kernel} conditionally on $\Delta$ yields
$\E\left[Tf(T)\right]=\E\left[Qf'(T)\right]$.
The identity is integrable: $T$ has a finite second moment and
$\E\left[|Q|\right]\le n^{-1}\|C_W\|_{\op}\E\left[\tr|\mathcal S_x|\right]<\infty$.
Together with \eqref{eq:pred-proof-kernel-error}, it follows that
\begin{equation}
\label{eq:pred-proof-stein}
 |\E\left[Tf(T)\right]-a^2\E\left[f'(T)\right]|
 \le C_{\mathcal D_+}\left(
     n^{-1/4}\|f'\|_\infty+n^{-1/2}\|f''\|_\infty\right).
\end{equation}
Thus no concentration of $H_x$ itself is required.

Apply Lemma~\ref{lem:pred-proof-stein-equation} to
\eqref{eq:pred-proof-stein}. For $a>0$ this gives
\[
 \Wass(\Law(T),N(0,a^2))
 \le C_{\mathcal D_+}\left(\frac{n^{-1/4}}a+
                              \frac{n^{-1/2}}{a^2}\right).
\]
If $a\ge n^{-1/8}$ this is $O(n^{-1/8})$.
If $a<n^{-1/8}$, including $a=0$, use instead
$\Wass(\Law(T),N(0,a^2))\le\E\left[|T|\right]+a\E\left[|Z|\right]\le2a$.
Consequently,
\begin{equation}
\label{eq:pred-proof-unconditional}
 \Wass(\Law(T),N(0,a^2))\le C_{\mathcal D_+}n^{-1/8}.
\end{equation}

\paragraph{Conditioning on the fresh feature.}
Write $P_{\bar x}:=\Law_\Delta(\bar x^\top\Delta/\sqrt n)$ and
$\bar P:=\E_{\bar x}\left[P_{\bar x}\right]=\Law(T)$.
For $\varepsilon_{\rm sm}>0$, define
$P_{\bar x}^{\varepsilon_{\rm sm}}:=P_{\bar x}*N(0,\varepsilon_{\rm sm}^2)$ and
$\bar P^{\varepsilon_{\rm sm}}:=\bar P*N(0,\varepsilon_{\rm sm}^2)$.
The CDF of the first law is
\[
 F_{\bar x}^{\varepsilon_{\rm sm}}(u)
 =\E_{\Delta}\left[\Phi_{\rm N}\!\left(
      \frac{u-\bar x^\top\Delta/\sqrt n}{\varepsilon_{\rm sm}}\right)\right],
\]
where $\Phi_{\rm N}$ and $\varphi_{\rm N}$ are the standard normal CDF and density.
Since
\[
 \nabla_{\bar x} F_{\bar x}^{\varepsilon_{\rm sm}}(u)
 =-\frac1{\varepsilon_{\rm sm}\sqrt n}
    \E_{\Delta}\!\left[\Delta\,
       \varphi_{\rm N}\!\left(\frac{u-\bar x^\top\Delta/\sqrt n}{\varepsilon_{\rm sm}}\right)\right],
 \qquad
 |\nabla_{\bar x} F_{\bar x}^{\varepsilon_{\rm sm}}(u)|^2
 \le\frac{\|C_W\|_{\op}}{2\pi n\varepsilon_{\rm sm}^2},
\]
Poincar\'e gives, uniformly in $u$,
\[
 \E_{\bar x}\left[\left|
 F_{\bar x}^{\varepsilon_{\rm sm}}(u)-F_{\bar P^{\varepsilon_{\rm sm}}}(u)\right|\right]
 \le \frac{C_{\mathcal D_+}}{\varepsilon_{\rm sm}\sqrt n}.
\]
Use the one-dimensional identity
$\Wass(P,Q)=\int_\R|F_P-F_Q|\,\dd u$
\citep{Vallender1974}.
On $[-R_{\rm cdf},R_{\rm cdf}]$ the preceding bound contributes
$C_{\mathcal D_+}R_{\rm cdf}/(\varepsilon_{\rm sm}\sqrt n)$.
Outside this interval, the averaged second moment
$\E_{\bar x}\left[\int u^2P_{\bar x}^{\varepsilon_{\rm sm}}(\dd u)\right]=a^2+\varepsilon_{\rm sm}^2$
bounds the two CDF tails by $2(a^2+\varepsilon_{\rm sm}^2)/R_{\rm cdf}$.
Finally, adding independent $\varepsilon_{\rm sm} Z$ costs at most
$\varepsilon_{\rm sm}\E\left[|Z|\right]$ for each law. Thus, for $0<\varepsilon_{\rm sm}\le1$ and $R_{\rm cdf}\ge1$,
\[
 \E_{\bar x}\left[\Wass(P_{\bar x},\bar P)\right]
 \le C_{\mathcal D_+}\left(
      \varepsilon_{\rm sm}+\frac{R_{\rm cdf}}{\varepsilon_{\rm sm}\sqrt n}+\frac1{R_{\rm cdf}}\right).
\]
Taking $\varepsilon_{\rm sm}=n^{-1/6}$ and $R_{\rm cdf}=n^{1/6}$, and using
\eqref{eq:pred-proof-unconditional}, gives
\begin{equation}
\label{eq:pred-proof-conditional-centered}
 \E_{\bar x}\left[\Wass(P_{\bar x},N(0,a^2))\right]
 \le C_{\mathcal D_+}(n^{-1/6}+n^{-1/8}).
\end{equation}

\paragraph{Restoring the feature mean.}
Conditioning on $x$ is equivalent to conditioning on $\bar x=x-m_x$.
The coupling with the same $\Delta$ gives
\[
 \E\left[\left|\frac{m_x^\top\Delta}{\sqrt n}\right|\right]
 \le\sqrt{\frac{m_x^\top C_Wm_x}{n}}
 \le C_{\mathcal D_+}n^{-1/2}.
\]
Moreover, $\Sigma_x=C_x+m_xm_x^\top$, so
$\alpha^2-a^2=m_x^\top C_Wm_x/n\le C_{\mathcal D_+}/n$.
Coupling $aZ$ and $\alpha Z$ therefore costs at most
$\sqrt{2/\pi}\,|\alpha-a|\le C_{\mathcal D_+}n^{-1/2}$.
Apply these two couplings to
\eqref{eq:pred-proof-conditional-centered}, and translate both laws by
$x^\top\mu$. This proves \eqref{eq:prediction-directions}.
All estimates were uniform on $\mathcal D_+$ and in the omitted index;
taking the supremum proves the proposition.
\end{proof}

\subsection{Completion under uniform strong convexity}
\begin{proof}[Proof of Theorem~\ref{thm:prediction-directions}, strongly convex case]
Suppose \eqref{eq:strongly-convex-penalty} holds. Appendix~\ref{app:population-stability}
and \eqref{eq:population-moment-regime} give a uniform bound on $z_\ast$:
\[
 |\mu_\ast|^2\le\E_G\left[|W_\ast|^2\right]/n,\qquad
 0\le\nu_\ast\le L_2,\qquad 0\le\beta_\ast\le L_1,
 \qquad |\eta_\ast|\le L_2K|\mu_\ast|+\sqrt K L_1.
\]
Here the last bound follows directly from \eqref{eq:losscoeff}, without
empirical moment identification or transport. Write
$z(t)=(v(t),b(t),h(t))=(1-t)z_\ast$, and let
$\ell_t:=\sqrt n\,h(t)-\sqrt{2b(t)}G$. Strong convexity and the minimum
of $\rho$ at zero imply
$\rho(w)-\rho(0)\ge\gamma_{\rm sc}|w|^2/2$, also for nonsmooth penalties.
Comparison with zero and nonnegative losses give
\[
 \frac{\gamma_{\rm sc}}2|W_t|^2-|\ell_t||W_t|\le nL_0,
 \qquad \sup_{n,\,t\in[0,1]}\frac{\E\left[|W_t|^2\right]}n\le C,
\]
using Young's inequality, the bound on $z_\ast$, and
$\tr\Sigma_x=O(n)$. At $t=1$ the strongly convex empirical objective
has a unique measurable minimizer by the same response calculus.

Fix $0<\delta\le1$ and take $G'\sim N(0,\Sigma_x)$ independent of all
original randomness. Add
\[
 \frac\delta2w^\top\Sigma_xw+\sqrt{2\delta}\,{G'}^\top w
\]
to the objective of $W_t$, and denote the optimizer by $W_t^\delta$.
Combining the two independent Gaussian fields shows that its marginal
law is that of $W_{t,(v(t)+\delta,b(t)+\delta,h(t))}$.
Subgradient strong monotonicity gives
\[
 |W_t^\delta-W_t|\le\gamma_{\rm sc}^{-1}
       |\delta\Sigma_xW_t+\sqrt{2\delta}\,G'|.
\]
Independence, the moment bound, and the spectral bounds therefore imply
\begin{equation}
\label{eq:strongly-convex-probe-coupling}
 \sup_{t\in[0,1]}\frac{\E\left[|W_t^\delta-W_t|^2\right]}n
 \le C(\delta^2+\delta)\le C\delta.
\end{equation}
At fixed $\delta>0$ these perturbed parameters lie in a fixed
positive-probe set. Propositions~\ref{prop:prediction-positive-probes} and~\ref{prop:gc-transfer}, with each vector's own moments, give
\[
 \sup_{t\in[0,1]}\mathcal G(W_t)
 \le C_\delta n^{-1/8}+C\sqrt\delta.
\]
First let $n\to\infty$, then $\delta\downarrow0$. This proves
Theorem~\ref{thm:prediction-directions} in the strongly convex branch.
No bound on $C_\delta$ as $\delta\downarrow0$ is asserted, so this argument
supplies no uniform rate in $n$ for the full path.
\end{proof}

\paragraph{Lean verification scope.}
The companion Lean~4 package verifies the quantitative comparison in
Proposition~\ref{prop:prediction-positive-probes} and its concentration
tools on arbitrary parameter domains with possibly degenerate selectors,
under uniform positive-probe bounds. It also verifies the strongly convex
minimizer perturbation and mean-square estimates used in
\eqref{eq:strongly-convex-probe-coupling}, under explicit moment and
measurability hypotheses, and an abstract cutoff-removal lemma.
The \href{https://anonymous.4open.science/r/ERM_non_smooth_penalty-7437/}{companion repository}
gives the exact statement correspondence, pins Lean and mathlib, and
supplies checks for \texttt{sorry} and custom axioms. The complete
calibrated-path instantiation and the pure-Lasso proof in
Appendix~\ref{app:proof-main-theorem} remain outside the formal scope;
Theorem~\ref{thm:prediction-directions} is not claimed as formally verified.

\section{Three-step proof of Theorem~\ref{thm:main-interior-HJ}}
\label{app:four}
The first three subsections decompose and bound the value error.
Appendix~\ref{app:active-coordinate-calculus} then records a shared
response formula for Lemma~\ref{lem:lasso-support-upper-bound} and
Corollary~\ref{cor:correlated-sparse-response}.

\subsection{Cavity reduction and the three signed errors}
The proof follows the chain
\[
 \text{insertion cost}\ \longrightarrow\ \text{scalar Moreau envelope}
 \ \longrightarrow\ \text{cavity Hamiltonian}
 \ \longrightarrow\ \text{full Hamiltonian}.
\]
Only the scalar-envelope-to-cavity-Hamiltonian step uses
Proposition~\ref{prop:prediction-positive-probes}.
Write $\varepsilon_{\rm pred}(\mathcal D)$ for the supremum of the error
in \eqref{eq:prediction-directions} restricted to $\mathcal D$, including
the full/cavity family.
For the insertion-cost-to-scalar-envelope step,
Lemma~\ref{lem:derived-scalar-response} gives unforced scalar-response
concentration. Proposition~\ref{prop:r1} transfers it to bounded forces
and compares the insertion cost with the scalar Moreau envelope.

Delete observation 1 and write $\Psi_-$ for its background objective,
$W_-$ for its minimizer, and $\mu_-,\alpha_-^2,\kappa_-$ for its deterministic
moments and response. The deleted pair $(x,y)$ is independent of this background.
Set
\begin{equation}
\label{eq:increment}
 \mathcal I:=\min_w\{\Psi_-(w)+L_y(x^\top w/\sqrt n)\}-\min_w\Psi_-(w).
\end{equation}
For $|r|\le L_1$, retain the cavity data and Gaussian probe used to define
$W_-$, and set
\[
 \begin{gathered}
 W_-(r):=\argmin_w\{\Psi_-(w)-rx^\top w/\sqrt n\},\\
 q_x(r):=x^\top W_-(r)/\sqrt n,\qquad
 J_-(r):=\frac1{\sqrt n}D_hW_-(r),\qquad J_-:=J_-(0).
 \end{gathered}
\]
The finite-force error used in Proposition~\ref{prop:r1} is
\begin{equation}
\label{eq:response-quantities}
 \omega_{\rm lin}:=\E\left[\sup_{|r|\le L_1}
             |q_x(r)-q_x(0)-\kappa_-r|\right].
\end{equation}
For cavity-to-full deletion, write
$\omega_\mu:=\sqrt{(\mu-\mu_-)^\top\Sigma_x(\mu-\mu_-)}$,
$\omega_{\alpha^2}:=|\alpha^2-\alpha_-^2|$, and
$\omega_\kappa:=|\kappa-\kappa_-|$.
With $u:=q_x(0)$, define the three signed errors by
\begin{align}
\label{eq:r1}
 r^{(1)}&:=\E\left[\mathcal I-e_{L_y}(u;\kappa_-)\right],\\
\label{eq:r2}
 r^{(2)}&:=\E\left[e_{L_y}(u;\kappa_-)\right]-\cH(\mu_-,\alpha_-,\kappa_-),\\
\label{eq:r3}
 r^{(3)}&:=\cH(\mu_-,\alpha_-,\kappa_-)-\widetilde\cH(\nabla\phi).
\end{align}
Summing \eqref{eq:r1}--\eqref{eq:r3} cancels the intermediate terms.
Since $\E\left[\mathcal I\right]=\partial_t\phi$ by
Lemma~\ref{lem:interior-envelope}, the value defect satisfies
\begin{equation}
\label{eq:three-error-decomposition}
 r:=\partial_t\phi-\widetilde\cH(\nabla\phi)=\sum_{j=1}^3r^{(j)}.
\end{equation}

\subsection{Scalar-response concentration}
\label{app:interior-concentration}
Throughout this section, $\mathcal D$ is a fixed interior domain introduced
in Section~\ref{sec:HJ-pseudo-equation}. Finite differences are taken in a
fixed enlargement with fixed positive lower bounds for $v$ and $b$. All
estimates are uniform on $\mathcal D$, and constants may depend on its compact
ranges, the bound on $|h|$, and the margins of the enlarged domain. The
response argument is
\[
 \begin{gathered}
 \text{minimum-value concentration}
 \ \Longrightarrow\ \text{amplitude-derivative concentration}\\
 \Longrightarrow\ \text{response-trace concentration}.
 \end{gathered}
\]
The quadratic-form estimate \eqref{eq:feature-quadratic-concentration} and
independence yield \eqref{eq:derived-scalar-response};
the finite-force step in Proposition~\ref{prop:r1} transfers this estimate
to bounded forces.

\begin{lemma}[Scalar response from feature concentration]
\label{lem:derived-scalar-response}
Under the assumptions of Lemma~\ref{lem:interior-moments}\textup{(\ref{item:interior-fourth-moments})},
full and unforced cavity responses satisfy
\begin{equation}
\label{eq:derived-trace-concentration}
 \sup_{\mathcal D}\left\|\tr(\Sigma_xJ)/n-\kappa\right\|_{L^2}\longrightarrow0.
\end{equation}
Consequently, for a fresh feature independent of $J$,
\begin{equation}
\label{eq:derived-scalar-response}
 \varepsilon_{\rm resp}(\mathcal D):=
 \sup_{\mathcal D}\E\left[|x^\top Jx/n-\kappa|\right]\longrightarrow0.
\end{equation}
Each family uses its own deterministic $\kappa$. Neither prediction
comparison nor projection anti-concentration is used.
\end{lemma}
\begin{proof}
Keep the observations, selectors and probe fixed, and put $\sigma:=\sqrt{2b}$ and
\[
 \mathcal V(\sigma):=\frac1n\min_w\left\{\rho(w)
 +\sum_iB_iL_{y_i}(x_i^\top w/\sqrt n)
 +\frac v2w^\top\Sigma_xw+\sigma G^\top w-\sqrt n\,h^\top w\right\}.
\]
Comparison with zero and the supporting loss hyperplanes bound
$|\mathcal V(\sigma)-\rho(0)/n|$ by
$L_0+|n^{-1/2}\sum_iB_iL_{y_i}'(0)x_i+\sigma G-\sqrt n\,h|^2/
      (2v\underline c_xn)$,
which has bounded $L^2$ norm by the comparison-field estimate in the proof
of Lemma~\ref{lem:interior-moments}\textup{(\ref{item:interior-fourth-moments})}.
The insertion cost $D_i$ satisfies
\[
 0\le D_i\le L_0+L_1|x_i^\top W_-|/\sqrt n,
 \qquad \E\left[D_i^2\right]\le C_{\mathcal D},
\]
because the cavity $W_-$ is independent of $(x_i,y_i)$.
Replacing one observation block costs two insertions into a common cavity.
Efron--Stein and Gaussian Poincar\'e therefore give
\[
 \Var(\mathcal V(\sigma))\le\frac{C_{\mathcal D}}n
       +\frac{\sigma^2}{n^2}\E\left[W^\top\Sigma_xW\right]
 \le\frac{C_{\mathcal D}}n.
\]
The envelope theorem gives $\mathcal V'(\sigma)=G^\top W/n$; differentiating
$\mathcal V'(\sigma)=G^\top W/n$ and using $D_GW=-\sigma J$ gives
$\mathcal V''(\sigma)=-G^\top JG/n$ almost everywhere, with
$\|\mathcal V''\|_{L^2}\le C_{\mathcal D}$.
Comparing $\mathcal V'(\sigma)$ with an amplitude difference of length
$n^{-1/4}$, which lies inside the buffered domain for large $n$, gives
\begin{equation}
\label{eq:probe-amplitude-derivative-concentration}
 \left\|\frac{G^\top W-\E\left[G^\top W\right]}n\right\|_{L^2}
 \le C_{\mathcal D}n^{-1/4}.
\end{equation}

Conditional on the data, $D_GW=-\sigma J$ and the second-order Gaussian
Stein identity \citep{BellecZhangSecondOrderStein} gives
\begin{equation}
\label{eq:second-order-stein-response}
 \E\left[(G^\top W+\sigma\tr(\Sigma_xJ))^2\right]
 =\E\left[W^\top\Sigma_xW\right]
  +\sigma^2\E\left[\tr((\Sigma_x^{1/2}J\Sigma_x^{1/2})^2)\right]
 \le C_{\mathcal D}n.
\end{equation}
The response bound in Lemma~\ref{lem:interior-wellposedness} justifies the
identity \eqref{eq:second-order-stein-response} by Gaussian $W^{1,2}$
approximation. First-order integration by parts centers the left-hand side
of \eqref{eq:second-order-stein-response}. Since $\sigma$ is bounded below
on $\mathcal D$,
\eqref{eq:probe-amplitude-derivative-concentration} proves
\eqref{eq:derived-trace-concentration}.
Finally, condition on the unforced $J$, which is independent of $x$:
\[
 \E\left[|x^\top Jx/n-\kappa|\right]
 \le\frac{K_Q}{n}\E\left[\|J\|_{\hs}\right]
       +\left\|\tr(\Sigma_xJ)/n-\kappa\right\|_{L^2}
 \le\frac{K_Q\sqrt p}{v\underline c_xn}+o(1).
\]
This proves the scalar response; the quadratic-form estimate
\eqref{eq:feature-quadratic-concentration} is never conditioned on a response
depending on $x$ through a nonzero force.
\end{proof}

\subsection{Bounding the three errors}
\label{app:r1}
\begin{proposition}[Scalar-response error]
\label{prop:r1}
Under the assumptions of Lemma~\ref{lem:interior-wellposedness}, with
$v,b>0$,
\[
 |r^{(1)}|\le L_1\omega_{\rm lin}.
\]
If, in addition, the training pairs are independent copies of $(x,y)$,
$p/n\to c\in(0,\infty)$, $\Sigma_x\preceq KI_p$, and
Assumption~\ref{ass:feature-concentration} holds, then on every fixed
interior domain $\mathcal D$,
\[
 \sup_{\mathcal D}|r^{(1)}|\longrightarrow0.
\]
\end{proposition}
\begin{proof}
\emph{Finite-force linearization.}
For convergence, work under the additional assumptions.
Lemma~\ref{lem:derived-scalar-response} gives
\eqref{eq:derived-scalar-response}. A bounded force is a translation of
the actual Gaussian probe. Use
the Gaussian translation inequality \eqref{eq:unforced-shift-bound}, proved
in Lemma~\ref{lem:pred-proof-concentration}, with
$A_R:=\{|x|\le R\sqrt n\}$.
With
$H(G):=|x^\top J_-(G)x/n-\kappa_-|$,
\eqref{eq:derived-scalar-response} and
$H\le C_{\mathcal D}(1+R^2)$ on $A_R$ give
\[
 \sup_{|r|\le L_1}\E\left[|x^\top J_-(r)x/n-\kappa_-|\right]
 \le C_{\mathcal D,R}\sqrt{\varepsilon_{\rm resp}(\mathcal D)}
   +C_{\mathcal D}\E_x\left[(1+|x|^2/n)\mathbf1_{A_R^c}\right]=o(1).
\]
Choose one fixed sufficiently large $R$ as in
\eqref{eq:pred-proof-tail}. That estimate and
$\Pp(A_R^c)\le\E[|x|^2\mathbf1_{A_R^c}]/(R^2n)$ make the tail term
$O(n^{-1})$.
Absolute continuity and Fubini give
$q_x(r)-q_x(0)=\int_0^r x^\top J_-(s)x/n\,\dd s$ almost surely. Thus
\[
 \omega_{\rm lin}\le\int_{-L_1}^{L_1}
       \E\left[|x^\top J_-(s)x/n-\kappa_-|\right]\,\dd s.
\]
The supremum in $\omega_{\rm lin}$ also covers a random insertion force;
no forced-optimizer prediction comparison is needed.
Consequently,
\begin{equation}
\label{eq:finite-scalar-response}
 \sup_{\mathcal D}\omega_{\rm lin}\longrightarrow0.
\end{equation}

\emph{Fenchel comparison.}
This step uses only the basic response assumptions.
The envelope derivative of
$\min_w\{\Psi_-(w)-rx^\top w/\sqrt n\}$ is $-q_x(r)$.
Fenchel duality, applicable since the loss is finite and continuous, gives
\[
 \mathcal I=\sup_z\left\{\int_0^zq_x(-s)\,\dd s-L_y^*(z)\right\},\qquad
 e_{L_y}(u;\kappa_-)=\sup_z\{zu-\kappa_-z^2/2-L_y^*(z)\}.
\]
Both suprema use $\operatorname{dom}L_y^*\subset[-L_1,L_1]$.
Their integrands differ by at most
$L_1\sup_{|r|\le L_1}|q_x(r)-u-\kappa_-r|$; take expectations.
This proves the first bound; \eqref{eq:finite-scalar-response} gives the
uniform convergence under the additional assumptions.
\end{proof}

\begin{proposition}[Gaussian-comparison error]
\label{prop:r2}
Under the assumptions of Proposition~\ref{prop:prediction-positive-probes},
$|r^{(2)}|\le L_1\varepsilon_{\rm pred}(\mathcal D)$ uniformly on
$\mathcal D$.
\end{proposition}
\begin{proof}
Condition on the independent pair $(x,y)$. Independence of $W_-$ from
this pair makes its score law depend only on $x$. Since
$e_{L_y}(\cdot;\kappa_-)$ is $L_1$-Lipschitz uniformly in $y$,
Wasserstein duality and averaging give
\[
 |r^{(2)}|\le L_1\E_x\left[\Wass\!\left(
 \Law(x^\top W_-/\sqrt n\mid x),N(x^\top\mu_-,\alpha_-^2)\right)\right].
\]
\end{proof}

\begin{proposition}[Cavity-to-full error]
\label{prop:r3}
Suppose Assumption~\ref{ass:analytic} holds,
$\underline c_xI_p\preceq\Sigma_x\preceq KI_p$, and $p/n\le K$. Let
$\mathcal D$ be a fixed interior domain introduced in
Section~\ref{sec:HJ-pseudo-equation}, and suppose its defining bounds remain
valid on a fixed enlargement used for finite differences. Then, uniformly on
$\mathcal D$,
\begin{equation}
\label{eq:cavity-full-error-bound}
 |r^{(3)}|\le L_1\omega_\mu+\frac{L_2}{2}\omega_{\alpha^2}
                         +\frac{L_1^2}{2}\omega_\kappa.
\end{equation}
All three discrepancies vanish uniformly, without prediction comparison
or feature quadratic-form concentration.
\end{proposition}
\begin{proof}
The scalar derivative bounds and Gaussian variance estimate in
Lemma~\ref{lem:derivative-calculus}\textup{(\ref{item:scalar-derivative-calculus})}
bound changes of the Hamiltonian in mean, variance and response by the
three terms on the right-hand side of \eqref{eq:cavity-full-error-bound}.
Couple $W,W_-$ using the same data and probe. Strong monotonicity of the
cavity objective gives, when the deleted observation is present,
\begin{equation}
\label{eq:nonsmooth-deletion}
 |W_+-W_-|\le L_1|x|/(v\underline c_x\sqrt n).
\end{equation}
Consequently $\|W-W_-\|_{L^2}\le C_{\mathcal D}$, and
\eqref{eq:coupling-moments} gives
$\omega_\mu+\omega_{\alpha^2}\le C_{\mathcal D}/\sqrt n$.
For the response, use the actual-probe correlation formula:
\begin{equation}
\label{eq:kappa-deletion-stein}
 \omega_\kappa=\frac{|\E\left[G^\top(W-W_-)\right]|}{n\sqrt{2b}}
 \le\frac{\sqrt{\tr(\Sigma_x)/n}\,\|W-W_-\|_{L^2}}{\sqrt{2bn}}
 \le\frac{C_{\mathcal D}}{\sqrt n}.
\end{equation}
Thus all three discrepancies vanish uniformly on $\mathcal D$.
\end{proof}

\begin{proof}[Proof of Theorem~\ref{thm:main-interior-HJ}]
Propositions~\ref{prop:r1}--\ref{prop:r3} and
\eqref{eq:three-error-decomposition} give
\begin{equation}
\label{eq:interior-value-rate}
 \sup_{\mathcal D}|\partial_t\phi-\widetilde\cH(\nabla\phi)|
 \le\sum_{j=1}^3\sup_{\mathcal D}|r^{(j)}|\longrightarrow0.
\end{equation}
\end{proof}

\subsection{A shared active-coordinate response formula}
\label{app:active-coordinate-calculus}
The response at the observation-selected insertion force is treated in
the proof of Lemma~\ref{lem:lasso-support-upper-bound}.

\begin{lemma}[Local active-coordinate response]
\label{lem:active-coordinate-response}
Let $W$ minimize $\rho_0(w)+\lambda\|w\|_1-\ell^\top w$, with $\rho_0$ convex
and $C^2$ and $\lambda>0$. Suppose its support $S$ and active signs $\sigma_S$
satisfy strict inactive slack $|\ell_j-\nabla_j\rho_0(W)|<\lambda$ for $j\notin S$,
and $[\nabla^2\rho_0(W)]_{SS}\succ0$. Then the global minimizer is unique
and, locally in $\ell$, has the same support and signs. It is continuously
differentiable with
\begin{equation}
\label{eq:active-background}
 D_\ell W=P_S([\nabla^2\rho_0(W)]_{SS})^{-1}P_S^\top,
\end{equation}
where $P_S$ embeds the active coordinates.
\end{lemma}
\begin{proof}
Apply the implicit-function theorem to
$[\nabla\rho_0(W)]_S+\lambda\sigma_S=\ell_S$, $W_{S^c}=0$.
In a small neighborhood, signs, strict slack and restricted positivity
persist, so this branch satisfies the full KKT conditions and is a global
minimizer. To identify it as the unique optimizer, let $a$ be its
$\ell_1$ subgradient. For any displacement $d$, convexity gives
\[
 \rho_0(W+d)+\lambda\|W+d\|_1-\ell^\top(W+d)
 -\{\rho_0(W)+\lambda\|W\|_1-\ell^\top W\}
 \ge\lambda\sum_{j\notin S}(1-|a_j|)|d_j|.
\]
Another minimizer must therefore have $d_{S^c}=0$. Convexity would make
its joining segment constant, contradicting
$d_S^\top[\nabla^2\rho_0(W)]_{SS}d_S>0$ unless $d=0$.
Differentiating the active KKT equations gives \eqref{eq:active-background}
(with zero derivative for empty support). The formula concerns local
support constancy, as in partial-smoothness sensitivity
\citep{LewisLiangTian2022,LewisZhang2013}; finite forces may cross thresholds.
\end{proof}

\section{Direct comparison and calibrated transport}
\label{app:value-transport}
We transport population moments along $z(t)=(1-t)z_\ast$ to each fixed
cutoff $t=1-\epsilon$ using a penalized value increment and short-time coupling.
The Gaussian input is only Proposition~\ref{prop:prediction-positive-probes},
through Theorem~\ref{thm:main-interior-HJ}.

\begin{lemma}[Forward transport from the calibrated population solution]
\label{lem:lasso-calibrated-transport}
Along a sequence in which the training pairs are independent copies of
$(x,y)$, suppose $p/n\to c\in(0,\infty)$,
$\underline c_xI_p\preceq\Sigma_x\preceq KI_p$,
Assumptions~\ref{ass:analytic} and~\ref{ass:feature-concentration} and
\eqref{eq:population-moment-regime} hold. For fixed $0<\epsilon<1$, set
\begin{equation}
\label{eq:lasso-calibrated-cutoff}
 q^\epsilon:=(\mu^\epsilon,(\alpha^\epsilon)^2,\kappa^\epsilon).
\end{equation}
The components of $q^\epsilon$ are the actual moments and response of
$W_{1-\epsilon,\,\epsilon z_\ast}$ defined in \eqref{eq:target-moments}--\eqref{eq:orders}, and
\begin{equation}
\label{eq:lasso-calibrated-moments}
 \|q^\epsilon-q_\ast\|_1\longrightarrow0,
 \qquad q_\ast=(\mu_\ast,\alpha_\ast^2,\kappa_\ast).
\end{equation}
No penalty/design alternative from Assumption~\ref{ass:penalty-structure}
is required.
\end{lemma}
\begin{proof}
\emph{Regularity on the calibrated tube.}
Equations~\eqref{eq:population-root-second-moment}--\eqref{eq:population-root-correlation},
Proposition~\ref{prop:population-nondegeneracy}, and
\eqref{eq:population-moment-regime} give constants $0<c<C<\infty$ such that,
for all sufficiently large sizes,
\[
 |\mu_\ast|+\alpha_\ast+\kappa_\ast+\nu_\ast+\beta_\ast
       +|\eta_\ast|\le C,
 \qquad \min\{\nu_\ast,\beta_\ast\}\ge c.
\]
Put $\pi_\ast:=\nabla_z\phi(0,z_\ast)$. Since $W_{0,z_\ast}=W_\ast$,
Lemmas~\ref{lem:normalized-identities} and
\ref{lem:derivative-calculus}\textup{(\ref{item:hamiltonian-derivative-calculus})}
give $\nabla\widetilde\cH(\pi_\ast)=z_\ast$.
Fix $\epsilon$ and choose a dimension-independent $r_0>0$ so that the tube
$\{(t,z(t)+d):0\le t\le1-\epsilon,\ |d|\le r_0\}$ and a fixed enlargement
have bounded parameters and $v,b$ bounded away from zero.
For optimizers coupled with the same data, selectors and probe,
\eqref{eq:monotone-lip} gives, with $v,v'\ge v_\circ>0$,
\[
 \underline c_xv_\circ\frac{|W-W'|}{\sqrt n}
 \le|h-h'|+|\sqrt{2b}-\sqrt{2b'}|\frac{|G|}{\sqrt n}
                  +K|v-v'|\frac{|W'|}{\sqrt n}.
\]
Lemma~\ref{lem:interior-moments}\textup{(\ref{item:interior-second-moments})}
and \eqref{eq:envelope} yield $\|W-W'\|_{L^2}/\sqrt n\le C|z-z'|$
and dimension-independent bounds $M,L$ on $|\nabla_z\phi|$ and its spatial
Lipschitz constant, with $|\pi_\ast|\le M$.
The same domination and finite Bernoulli averaging give joint $C^1$
regularity, with one-sided time derivatives at the endpoints.

The gradients and $\pi_\ast$ lie in the convex physical domain
\[
 \mathcal P:=\{\pi:-\pi_b\ge0,\quad
                  2\pi_v\ge\pi_h^\top\Sigma_x\pi_h\}.
\]
For two moment triples with $|\mu|,|\mu'|\le M$,
$0\le\kappa,\kappa'\le M$ and $\alpha,\alpha'\ge0$, put
$\zeta(u,\kappa,y):=\prox_{\kappa L_y}(u)$ and
$\xi(u,\kappa,y):=L_y'(\zeta(u,\kappa,y))$.
Implicit differentiation of $u=\zeta+\kappa L_y'(\zeta)$ gives, by
Assumption~\ref{ass:analytic} and $1+\kappa L_y''(\zeta)\ge1$,
\begin{equation}
\label{eq:population-xi-derivative-bounds}
 \sup_{u,y,\,0\le\kappa\le M}
 \left(|\xi|+|\partial_u\xi|+|\partial_u^2\xi|+|\partial_u^3\xi|
 +|\partial_\kappa\xi|+|\partial_\kappa\partial_u\xi|\right)\le C.
\end{equation}
Interpolate $(\mu_r,\alpha_r^2,\kappa_r)$ linearly from the primed to the
unprimed triple, and write $U_r:=x^\top\mu_r+\alpha_rZ$ and
$(\Delta_\mu,\Delta_{\alpha^2},\Delta_\kappa)
 :=(\mu-\mu',\alpha^2-(\alpha')^2,\kappa-\kappa')$.
For $f=f(u,\kappa,y)$ with bounded derivatives, conditioning on $(x,y)$
and differentiating the Gaussian variance give
\begin{equation}
\label{eq:population-coefficient-path}
\begin{aligned}
 &\E[f(U_1,\kappa_1,y)-f(U_0,\kappa_0,y)]\\
 &\quad=\int_0^1\Bigl\{
 \E[(x^\top\Delta_\mu)\partial_u f(U_r,\kappa_r,y)]
 +\frac{\Delta_{\alpha^2}}2\E[\partial_u^2 f(U_r,\kappa_r,y)]
 +\Delta_\kappa\E[\partial_\kappa f(U_r,\kappa_r,y)]
 \Bigr\}\,\dd r.
\end{aligned}
\end{equation}
At $\alpha^2=0$ and $\kappa=0$, use right derivatives, as in
\eqref{eq:gaussian-variance-derivative}.
Apply \eqref{eq:population-coefficient-path} to
$f=\partial_u\xi,\xi^2,\xi$, retaining $e^\top x$ for the last choice.
Using \eqref{eq:population-xi-derivative-bounds} and, for $|e|=1$,
\[
 \E|x^\top\Delta_\mu|\le\sqrt K|\Delta_\mu|,\quad
 \E[|e^\top x|\,|x^\top\Delta_\mu|]\le K|\Delta_\mu|,\quad
 \E|e^\top x|\le\sqrt K,
\]
the identities $\nu=\E[\partial_u\xi]$, $\beta^2=\E[\xi^2]$ and
$\eta=\nu\Sigma_x\mu-\E[x\xi]$ give, taking the supremum over $|e|=1$,
\[
 |\nu-\nu'|+|\beta^2-(\beta')^2|+|\eta-\eta'|
 \le C\bigl(|\Delta_\mu|+|\Delta_{\alpha^2}|+|\Delta_\kappa|\bigr),
\]
with $C$ depending only on $M,K,L_1,L_2,L_3,L_4$.
Thus \eqref{eq:Hgradient} and \eqref{eq:envelope} give a dimension-independent
Lipschitz constant $L_H$ for $\nabla\widetilde\cH$ on bounded subsets of
$\mathcal P$, and integration along segments gives
\begin{equation}
\label{eq:physical-Hamiltonian-Taylor}
 |\widetilde\cH(\pi)-\widetilde\cH(\pi')-\nabla \widetilde\cH(\pi')^\top(\pi-\pi')|
 \le\frac{L_H}{2}|\pi-\pi'|^2.
\end{equation}
The same bound holds on boundary segments by translating both endpoints
by $\delta(e_v-e_b)$ and letting $\delta\downarrow0$.

\emph{A deterministic gradient-transport estimate.}
Fix $0<s<1-\epsilon$. If the value residual $r(t,z)$ in
\eqref{eq:three-error-decomposition} satisfies $|r(t,z)|\le e$
on the tube over $[s,1-\epsilon]$, then
\begin{equation}
\label{eq:direct-gradient-transport-bound}
 \sup_{s\le t\le1-\epsilon}|\nabla_z\phi(t,z(t))-\pi_\ast|^2
 \le C\bigl(|\nabla_z\phi(s,z(s))-\pi_\ast|^2+e\bigr),
\end{equation}
with $C$ depending only on $M,L,L_H,r_0,\epsilon$, independently of $n,s$.
Choose $A>L$ with $Ar_0>2M$, and define
\[
 \mathcal E(t):=\max_{|d|\le r_0}
 \{\phi(t,z(t)+d)-\phi(t,z(t))-\pi_\ast^\top d-A|d|^2/2\}.
\]
Strong concavity and the boundary upper bound $2Mr_0-Ar_0^2<0$
on the outward derivative make the unique maximizer $d_t$ interior
and continuous, with
\[
 \nabla_z\phi(t,z(t)+d_t)-\pi_\ast=Ad_t,\qquad
 |d_t|\le\frac{|\nabla_z\phi(t,z(t))-\pi_\ast|}{A-L}.
\]
Taylor bounds, testing the admissible point
$d=(\nabla_z\phi(t,z(t))-\pi_\ast)/(A+L)$ for the lower bound, give
\begin{equation}
\label{eq:direct-transport-energy-bounds}
 \frac{|\nabla_z\phi(t,z(t))-\pi_\ast|^2}{2(A+L)}
 \le\mathcal E(t)\le
 \frac{|\nabla_z\phi(t,z(t))-\pi_\ast|^2}{2(A-L)}.
\end{equation}
The envelope theorem, $\dot z=-\nabla\widetilde\cH(\pi_\ast)$, and
\eqref{eq:physical-Hamiltonian-Taylor} give
\[
 \mathcal E'(t)\le\frac{L_H}{2}
   (A^2|d_t|^2+|\nabla_z\phi(t,z(t))-\pi_\ast|^2)+2e
 \le C\mathcal E(t)+2e.
\]
Gr\"onwall's inequality and \eqref{eq:direct-transport-energy-bounds}
prove \eqref{eq:direct-gradient-transport-bound}.

\emph{Initialization and application.}
Couple $W_{s,z_\ast}$ and $W_\ast$ with the same probe for small $s>0$.
Strong convexity, objective comparison and nonnegative retained losses give
\begin{equation}
\label{eq:initial-short-time-coupling}
 \frac1n\E\left[|W_{s,z_\ast}-W_\ast|^2\right]
 \le Cs\,\E_{x,y,G}\left[L_y(x^\top W_\ast/\sqrt n)\right]\le Cs.
\end{equation}
The last bound uses independence, \eqref{eq:population-moment-regime}
and $L_y(u)\le L_0+L_1|u|$.
Changing $z_\ast$ to $z(s)$ costs $O(s)$ in normalized $L^2$ by
the coupling estimate in the first step;
\eqref{eq:envelope} gives
$|\nabla_z\phi(s,z(s))-\pi_\ast|^2\le Cs$.
For fixed $s,\epsilon$, Theorem~\ref{thm:main-interior-HJ} supplies a uniform
value error tending to zero on the remaining tube. Hence
\eqref{eq:direct-gradient-transport-bound} gives
\[
 \limsup_{n\to\infty}|\nabla_z\phi(1-\epsilon,\epsilon z_\ast)-\pi_\ast|^2
 \le C_\epsilon s.
\]
Let $s\downarrow0$ after $n\to\infty$ at fixed $\epsilon$ and apply
\eqref{eq:envelope} to obtain \eqref{eq:lasso-calibrated-moments}.
\end{proof}

\section{Cutoff removal, completion of Gaussian prediction, and endpoint moments}
\label{app:proof-main-theorem}
We complete Theorem~\ref{thm:prediction-directions} for pure Lasso and prove
Theorem~\ref{thm:endpoint}. The transport lemma uses only the positive-probe
comparison through the HJ equation. The support gap, localization and
restricted curvature below give the uniform terminal-interval comparison
\eqref{eq:lasso-calibrated-terminal-MSE}; its endpoint specialization removes
the cutoff. Throughout, the dimension limit precedes the cutoff limit.

Under Assumptions~\ref{ass:analytic} and~\ref{ass:penalty-structure}, the
empirical objective $F$ is continuous and coercive. Its minimizer is unique
in the strongly convex case; for Lasso, choose the minimum-norm minimizer
$\widehat W$ of its compact convex minimizer set.
The minimizers $W_k:=\argmin_w\{F(w)+|w|^2/(2k)\}$ satisfy
$|W_k|\le|\widehat W|$ and $F(W_k)-\min F\le|\widehat W|^2/(2k)$,
hence $W_k\to\widehat W$. Approximation over a countable dense set gives
measurability for strongly convex minimizers, including each $W_k$,
and therefore for this limit.

\subsection{Strongly convex penalties}
\begin{proposition}[Calibrated-cutoff removal under strong convexity]
\label{prop:strongly-convex-calibrated-cutoff}
Under the assumptions of Lemma~\ref{lem:lasso-calibrated-transport}, suppose
that $\rho$ is strongly convex with dimension-independent modulus
$\gamma_{\rm sc}>0$. Couple $\widehat W$ and $W_{1-\epsilon,\,\epsilon z_\ast}$ using the same
observations. For all sufficiently large sizes and $0<\epsilon<1$,
\begin{equation}
\label{eq:strongly-convex-calibrated-cutoff}
 \frac1n\E\left[|\widehat W-W_{1-\epsilon,\,\epsilon z_\ast}|^2\right]\le C\epsilon.
\end{equation}
\end{proposition}
\begin{proof}
We restore the omitted observations and remove the probes separately.
Let $I:=I(1-\epsilon)$ and
$\overline W:=\argmin_w F_I(w)$, where
$F_I(w):=\rho(w)+\sum_{i\in I}L_{y_i}(x_i^\top w/\sqrt n)$.
Comparison with zero and strong convexity give
\begin{equation}
\label{eq:strongly-convex-retained-radius}
 \frac{|\widehat W|^2+|\overline W|^2}{n}\le C,\qquad
 |\widehat W-\overline W|^2
 \le C\sum_{i\notin I}L_{y_i}(x_i^\top\overline W/\sqrt n).
\end{equation}
Conditional on the selectors and retained observations, each omitted
pair is independent of $\overline W$. Since $L_y(z)\le L_0+L_1|z|$,
\[
 \frac1n\E\left[|\widehat W-\overline W|^2\right]
 \le C\epsilon.
\]
Strong monotonicity of $\partial F_I$ compares the probed and unprobed
retained minimizers:
\[
 |W_{1-\epsilon,\,\epsilon z_\ast}-\overline W|
 \le C|\epsilon\nu_\ast\Sigma_x\overline W
       +\sqrt\epsilon\,\beta_\ast G-\epsilon\sqrt n\,\eta_\ast|.
\]
Proposition~\ref{prop:population-nondegeneracy}, \eqref{eq:strongly-convex-retained-radius},
and $\E_G\left[|G|^2\right]\le Kp$ bound the squared expectation divided
by $n$ by $C\epsilon$. The squared triangle inequality proves
\eqref{eq:strongly-convex-calibrated-cutoff}.
\end{proof}

\subsection{Pure Lasso}
We compare every optimizer in the terminal interval to one fixed calibrated
cutoff, as stated in \eqref{eq:lasso-calibrated-terminal-MSE}. The three lemmas below establish
the support gap, localization and restricted curvature used in
Proposition~\ref{prop:lasso-calibrated-cutoff-MSE}.

For Lemmas~\ref{lem:lasso-support-upper-bound}--\ref{lem:lasso-ac-curvature},
assume the training pairs are independent copies of $(x,y)$,
$p/n\to c\in(0,\infty)$, $\underline c_xI_p\preceq\Sigma_x\preceq KI_p$,
and Assumptions~\ref{ass:analytic} and~\ref{ass:feature-concentration} hold,
with constants uniform in $n$.

\begin{lemma}[An upper bound for the interior Lasso support]
\label{lem:lasso-support-upper-bound}
Let $\rho(w)=\lambda\|w\|_1$, $\lambda>0$. For $W:=W_{t,z}$ and its
deterministic response $\kappa$, let $\mathcal D$ be a fixed interior domain
introduced in Section~\ref{sec:HJ-pseudo-equation}, and suppose its defining
bounds remain valid on a fixed enlargement used for finite differences. Then
\begin{equation}
\label{eq:lasso-interior-support-bound}
 \sup_{\mathcal D}\E\left[
 \left(\frac{|\operatorname{supp}(W)|}{n}
       -v\kappa-t\frac{L_2\kappa}{1+L_2\kappa}\right)_+\right]\to0,
 \qquad u_+:=\max\{u,0\}.
\end{equation}
If the assumptions of Lemma~\ref{lem:lasso-calibrated-transport} also hold,
there are $s_0<1$ and $\epsilon_1>0$ such that, for each fixed
$0<\epsilon<\epsilon_1$,
\begin{equation}
\label{eq:lasso-calibrated-support}
 \Pp\{|\operatorname{supp}(W_{1-\epsilon,\,\epsilon z_\ast})|>s_0n\}\to0.
\end{equation}
\end{lemma}
\begin{proof}
\emph{Active-coordinate response and the support trace identity.}
Write $S:=\operatorname{supp}(W)$, $J:=D_hW/\sqrt n$,
$q_i:=x_i^\top W/\sqrt n$, and
$H:=v\Sigma_x+n^{-1}\sum_i B_i(t)L_{y_i}''(q_i)x_ix_i^\top$.
Inactive KKT inequalities are strict almost surely: conditional on
$G_{-j}$, the optimizer constrained by $w_j=0$ is independent of $G_j$,
and its $j$th residual is affine with nonzero slope in a Gaussian of
variance $1/(\Sigma_x^{-1})_{jj}>0$. Lemma~\ref{lem:active-coordinate-response}
therefore gives $J=P_SH_{SS}^{-1}P_S^\top$, hence the samplewise identity
\begin{equation}
\label{eq:lasso-samplewise-active-trace}
 \frac{|S|}{n}=v\frac{\tr(\Sigma_xJ)}n
 +\frac1n\sum_i B_i(t)L_{y_i}''(q_i)\frac{x_i^\top Jx_i}{n}.
\end{equation}

\emph{Response at the optimizing insertion force.}
The assumptions above give \eqref{eq:derived-scalar-response} by
Lemma~\ref{lem:derived-scalar-response}. We transfer this unforced estimate
to the observation-selected force as follows.
Let
$W_+(G)$ include a fresh observation $(x,y)$ and let $W_-(G)$ be its
unforced cavity, with
$G\sim N(0,\Sigma_x)$ and $\sigma:=\sqrt{2b}$. Conditional on the cavity
data and fresh pair, set
\[
 u_-(z):=x^\top W_-(z)/\sqrt n,\qquad
 T(z):=z-\frac{L_y'(u_-(z))}{\sigma\sqrt n}x.
\]
We absorb the inserted loss into the Gaussian field. The resulting map
has a positive rank-one Jacobian, which controls its density and transfers
the unforced response to the optimizing force.
Substitution in the KKT inclusions gives $W_+(T(z))=W_-(z)$, while
\[
 \nabla_z u_-(z)=-\sigma J_-(z)x/\sqrt n,\qquad
 DT(z)=I+L_y''(u_-(z))xx^\top J_-(z)/n.
\]
The map preserves coordinates orthogonal to $x$. On every parallel line,
monotonicity gives a longitudinal slope at least one, and the displacement
is bounded by $L_1|x|/(\sigma\sqrt n)$; hence it is onto and one-to-one.
The Lipschitz dependence of $u_-$ also controls transverse changes, so $T$
and its inverse are Lipschitz for fixed $x$. The rank-one determinant
formula gives
\begin{equation}
\label{eq:optimizing-force-jacobian}
 \det DT(z)=1+L_y''(u_-(z))x^\top J_-(z)x/n\ge1,
 \qquad W_+(T(z))=W_-(z).
\end{equation}

Since the bi-Lipschitz map preserves null sets, differentiating
$W_+(G)=W_-(T^{-1}(G))$ almost everywhere and applying the rank-one
inverse formula to $DT$ gives, for $z=T^{-1}(G)$,
\begin{equation}
\label{eq:optimizing-force-scalar}
 \frac{x^\top J_+(G)x}{n}
 =\frac{x^\top J_-(z)x/n}
        {1+L_y''(u_-(z))x^\top J_-(z)x/n}.
\end{equation}
This change of variables permits active-set changes along the insertion force.

Writing $\varphi_{\Sigma_x}$ for the $N(0,\Sigma_x)$ density, the density
of $T^{-1}(G)$ relative to that law is
\[
 \frac{\varphi_{\Sigma_x}(T(z))\det DT(z)}
                    {\varphi_{\Sigma_x}(z)}.
\]
On $A_R:=\{|x|\le R\sqrt n\}$ the determinant is bounded, and
$|L_y'|\le L_1$ gives a uniform conditional $L^2$ density bound using
Gaussian exponential moments. Thus
\begin{equation}
\label{eq:optimizing-force-change-measure}
 \E\left[\mathbf1_{A_R}F(T^{-1}(G))\right]
 \le C_{\mathcal D,R}
          \left(\E\left[\mathbf1_{A_R}F(G)^2\right]\right)^{1/2}
\end{equation}
for nonnegative $F$, with its other arguments fixed conditionally.
For $F(G):=|x^\top J_-(G)x/n-\kappa_-|$,
\eqref{eq:derived-scalar-response} and
$F\le C_{\mathcal D}(1+R^2)$ on $A_R$ bound the right side by
$C_{\mathcal D,R}\sqrt{\varepsilon_{\rm resp}(\mathcal D)}$.
The complement contributes at most
$C_{\mathcal D}\E_x\left[(1+|x|^2/n)\mathbf1_{A_R^c}\right]$.
For one fixed sufficiently large $R$, \eqref{eq:pred-proof-tail} bounds
this contribution by $O(n^{-1})$. Therefore
\begin{equation}
\label{eq:optimizing-force-response}
 \sup_{\mathcal D}\E\left[
 |x^\top J_-(T^{-1}(G))x/n-\kappa_-|\right]\longrightarrow0.
\end{equation}

\emph{The support bound and the calibrated gap below $n$.}
For observation $i$, the cavity response at the optimizing force is
$J_{-i}(-L_{y_i}'(q_i))$ in the notation of Appendix~\ref{app:four}.
For $f(u):=L_2u/(1+L_2u)$, the scalar response identity
\eqref{eq:optimizing-force-scalar} bounds each summand in
\eqref{eq:lasso-samplewise-active-trace} by
\[
 L_{y_i}''(q_i)\frac{x_i^\top Jx_i}{n}
 \le f\left(\frac{x_i^\top J_{-i}(-L_{y_i}'(q_i))x_i}{n}\right).
\]
The matrix on the right depends on $(x_i,y_i)$. Conditional on
$B_i(t)=1$, \eqref{eq:optimizing-force-response}, the $L_2$-Lipschitz
property of $f$, and \eqref{eq:kappa-deletion-stein} replace this term
by $f(\kappa)$ with uniform $o_{L^1}(1)$ error. Also,
Lemma~\ref{lem:derived-scalar-response} gives
$\tr(\Sigma_xJ)/n=\kappa+o_{L^1}(1)$, while
$n^{-1}\sum_iB_i(t)=t+o_{L^1}(1)$. Substitution proves
\eqref{eq:lasso-interior-support-bound}.

Set $K_\kappa:=1+\sup_n\kappa_\ast$, which is finite by
\eqref{eq:population-moment-regime}.
Lemma~\ref{lem:lasso-calibrated-transport} gives
$\kappa^\epsilon\le K_\kappa$ eventually at each fixed cutoff.
With $\delta:=(1+L_2K_\kappa)^{-1}$ and $\nu_\ast\le L_2$,
\[
 \epsilon\nu_\ast\kappa^\epsilon
 +(1-\epsilon)\frac{L_2\kappa^\epsilon}{1+L_2\kappa^\epsilon}
 \le1-\delta+\epsilon L_2K_\kappa\le1-\delta/2
\]
for sufficiently small fixed $\epsilon$. Markov's inequality in
\eqref{eq:lasso-interior-support-bound}, with $s_0:=1-\delta/4$,
proves \eqref{eq:lasso-calibrated-support}.
\end{proof}

\begin{lemma}[Localization from a bounded Lasso root]
\label{lem:lasso-bounded-root-localization}
Let $\rho(w)=\lambda\|w\|_1$, with
$0<\lambda_-\le\lambda\le\lambda_+<\infty$, and suppose
$\sup_n\{\E_G\left[|W_\ast|^2\right]/n+\kappa_\ast\}<\infty$.
There are fixed $R,r_h>0$ and events $\mathcal L_n$, depending only on
observations and selectors, with
$\Pp(\mathcal L_n^c)\le Ce^{-cn}$ on which, simultaneously for
$t\in[1/2,1]$, $v\ge0$, and $|h|\le r_h$, the objective
\begin{equation}
\label{eq:lasso-localization-objective}
 \lambda\|w\|_1+\sum_{i\in I(t)}L_{y_i}(x_i^\top w/\sqrt n)
       +\frac v2w^\top\Sigma_xw-\sqrt n\,h^\top w
\end{equation}
is coercive and all its minimizers have norm less than $R\sqrt n$.
Here $I(t):=\{i:U_i\le t\}$ for independent uniform selectors $U_i$.
Moreover,
\begin{equation}
\label{eq:lasso-localization-exceptional-tail}
 \E\left[\frac{|\widehat W|^2}{n}\mathbf1_{\mathcal L_n^c}\right]
 \le Cne^{-cn}\to0.
\end{equation}
\end{lemma}
\begin{proof}
Set $f(w):=\E_{x,y,Z}[e_{L_y}(x^\top w/\sqrt n+\alpha_\ast Z;\kappa_\ast)]$.
We use the population radial growth bound
\begin{equation}
\label{eq:lasso-population-linear-growth}
 \E[L_y(x^\top w/\sqrt n)]+\lambda\|w\|_1/n
 \ge c_{\rm lin}|w|/\sqrt n-C_{\rm lin},
\end{equation}
and, for
$\mathcal T_{R,D}:=\{w:|w|\le R\sqrt n,\ \lambda\|w\|_1/n\le D\}$,
the empirical complexity estimate
\begin{equation}
\label{eq:lasso-localization-complexity}
 \E\left[\sup_{w\in\mathcal T_{R,D}}
 \left|\frac1{n\sqrt n}\sum_i\varepsilon_iB_i(1/2)x_i^\top w\right|\right]
 \le\inf_{T>0}\{DT/\lambda_-+CR/T\}
 \le C\sqrt{RD/\lambda_-},
\end{equation}
where the $\varepsilon_i$ are independent Rademacher signs.
Both bounds are verified below. Nonnegativity of the penalty preserves
\eqref{eq:lasso-population-linear-growth}, with $c_{\rm lin}/2$, when the
loss is multiplied by $1/2$. For small fixed $r_h$, put
$D:=L_0+r_hR+2$. On $|w|=R\sqrt n$, penalties exceeding $nD$ already
make \eqref{eq:lasso-localization-objective} exceed $n(L_0+1)$.
For the remaining vectors, fixed-height clipping and
\eqref{eq:lasso-localization-complexity} give, outside an exponentially
small event, the normalized lower bound
\[
 c_{\rm lin}R/2-C_{\rm lin}-r_hR
       -C\sqrt{R(L_0+r_hR+2)/\lambda_-}-C.
\]
For sufficiently small $r_h$ and then large fixed $R$, radial growth
dominates these fluctuations and the bound exceeds $L_0+1$.
Monotonicity of $I(t)$ and nonnegativity of the quadratic term make this
uniform in $t\in[1/2,1]$ and $v\ge0$. The value at zero is at most $nL_0$;
convexity along rays therefore gives coercivity and confines all minimizers
to $|w|<R\sqrt n$.

To verify the bounds, the surrogate KKT condition supplies
$s(G)\in\partial\|W_\ast(G)\|_1$ with
\begin{equation}
\label{eq:lasso-mean-subgradient-support}
 \frac\lambda n\E_G[s(G)]=-\nabla f(\E_GW_\ast),\qquad
 \frac\lambda n\|w\|_1\ge-\nabla f(\E_GW_\ast)^\top w.
\end{equation}
This averaged subgradient need not be a subgradient at the mean.
On $|w-\E_GW_\ast|\le\sqrt n$, local loss curvature gives
\[
 u^\top\nabla^2 f(w)u
 =\frac1n\E\left[(x^\top u)^2
       \frac{L_y''(\zeta)}{1+\kappa_\ast L_y''(\zeta)}\right]
 \ge c_{\rm lin}|u|^2/n,\qquad
 \zeta:=\prox_{\kappa_\ast L_y}(x^\top w/\sqrt n+\alpha_\ast Z).
\]
Indeed, Lemma~\ref{lem:qf-moments} bounds the projection tail by
$\E[(x^\top u)^2\mathbf1_{\{|x^\top w/\sqrt n|>T\}}]\le C|u|^2/T$;
choose large fixed $T$ and restrict also to $|Z|\le1$.
Integration along rays from the mean, \eqref{eq:lasso-mean-subgradient-support},
and $f(w)\le\E[L_y(x^\top w/\sqrt n)]+L_1\alpha_\ast\E|Z|$ prove
\eqref{eq:lasso-population-linear-growth}.
The same projection fourth moments bound the fourth moments of
$n^{-1/2}\sum_i\varepsilon_iB_i(1/2)x_{i,j}$ uniformly; splitting its
coordinates at height $T$ gives \eqref{eq:lasso-localization-complexity}.
For the clipped loss $\min\{L_y,H\}$, symmetrization and contraction give this
complexity bound up to constants, the clipping bias is $C(1+R^2)/H$,
and bounded differences adds at most one except with probability
$2e^{-2n/H^2}$. Choose $r_h$ with
$r_h+C\sqrt{r_h/\lambda_-}<c_{\rm lin}/8$, then $R$ as above and
fixed $H$ making the bias at most one. This defines $\mathcal L_n$ with
the claimed exponential probability. Finally,
$\lambda_-\|\widehat W\|_1\le nL_0$ on every sample, so
\eqref{eq:lasso-localization-exceptional-tail} follows.
\end{proof}

\begin{lemma}[Uniform curvature near coordinate spaces below the rank threshold]
\label{lem:lasso-ac-curvature}
Assume additionally \eqref{eq:projection-anticoncentration}.
For fixed $0\le s<1$ and $R<\infty$, there are
$\tau,\epsilon_0,\gamma_R>0$ and events $\mathcal A_n$ with
$\Pp(\mathcal A_n^c)\to0$ such that
\begin{equation}
\label{eq:lasso-ac-restricted-curvature}
 \frac1n\sum_{i\in I}L_{y_i}''(x_i^\top w/\sqrt n)(x_i^\top d)^2
 \ge\gamma_R|d|^2
\end{equation}
simultaneously whenever $|S|\le sn$,
$\|d_{S^c}\|_1\le\tau\sqrt n|d|$, $|w|\le R\sqrt n$, and
$|I^c|\le\epsilon_0n$. The sets $S$ and $I$ may depend on the data.
\end{lemma}
\begin{proof}
An anti-concentration net argument gives many nonzero projections,
uniformly over sparse directions. We extend it to small off-support
$\ell_1$ perturbations and discard rows with large scores.

First, for a sufficiently large fixed $B<\infty$,
\begin{equation}
\label{eq:lasso-ac-first-absolute-moment}
 \Pp\left\{\sup_{|u|=1}\frac1n\sum_i|x_i^\top u|>B\right\}
 \le Cn^{-1}.
\end{equation}
As a function of $(x_1,\ldots,x_n)$, this supremum is
$n^{-1/2}$-Lipschitz in the product Euclidean norm.
Symmetrization and contraction bound its mean by
$\sqrt K+4\sqrt{\tr(\Sigma_x)/n}$.
Tensorizing \eqref{eq:data-poincare} over the independent observations
and applying it to smooth bounded approximations with the same Lipschitz
bound gives variance at most $C_{\rm P}/n$. Chebyshev's inequality proves
\eqref{eq:lasso-ac-first-absolute-moment}.

On the event
\[
 \left\{\sup_{|u|=1}\frac1n\sum_i|x_i^\top u|\le B\right\},
\]
whose complement is bounded in
\eqref{eq:lasso-ac-first-absolute-moment}, choose $s<s'<1$ and
$1-4\omega_{\rm ac}>s'$. Nets of radius
$\delta:=\omega_{\rm ac}r/(4B)$ over all coordinate unit spheres of
dimension at most $s'n$ have total size at most
$2^p(1+2/\delta)^{s'n}$. For a fixed net point, anti-concentration and
independence bound the probability that at least
$(1-4\omega_{\rm ac})n$ projections have magnitude at most $2r$ by
$2^n(Cr)^{(1-4\omega_{\rm ac})n}$. Thus the union bound is
\[
 2^{p+n}(1+2/\delta)^{s'n}(Cr)^{(1-4\omega_{\rm ac})n}
 \le e^{-cn}
\]
for fixed small $r$. The strict exponent inequality is why every
$s<1$ is allowed. The first-absolute-moment bound extends the net
conclusion to all sparse unit vectors, each with at least
$3\omega_{\rm ac}n$ projections exceeding $r$ in magnitude.
This proves the needed data-selected uniformity directly; see
\citet{KoltchinskiiMendelson2013} for small-ball methods more generally.

For $|d|=1$, augment $S$ by the largest $(s'-s)n$ off-support coordinates.
The remaining Euclidean norm is at most $C\tau/\sqrt{s'-s}$.
Choose $\tau$ small: \eqref{eq:lasso-ac-first-absolute-moment} then leaves
at least $2\omega_{\rm ac}n$ projections of $d$ with magnitude at least
$r/4$. Uniformly in $|w|\le R\sqrt n$, the fraction of rows with
$|x_i^\top w/\sqrt n|>M$ is at most $BR/M$.
Choose $BR/M\le\omega_{\rm ac}/2$ and
$\epsilon_0\le\omega_{\rm ac}/2$. At least $\omega_{\rm ac}n$ rows remain
with squared projection at least $r^2/16$ and loss curvature at least
$c_L(M)$, proving \eqref{eq:lasso-ac-restricted-curvature}.
\end{proof}

\begin{proposition}[Uniform mean-square comparison with a calibrated cutoff]
\label{prop:lasso-calibrated-cutoff-MSE}
Under the pure-Lasso assumptions of
Theorem~\ref{thm:prediction-directions}, define
\[
 V_\delta:=W_{1-\delta,\,\delta z_\ast}\quad(0<\delta\le1),
 \qquad V_0:=\widehat W.
\]
Use the same observations, nested selectors
$I(t)=\{i:U_i\le t\}$ from independent uniforms, and the same Gaussian
probe $G$ for every $\delta$, with the minimum-norm selection at zero.
For every sufficiently small fixed $\epsilon>0$,
\begin{equation}
\label{eq:lasso-calibrated-terminal-MSE}
 \limsup_{n\to\infty}\sup_{0\le\delta\le\epsilon}
       \frac1n\E\left[|V_\delta-V_\epsilon|^2\right]
 \le C\epsilon^{1/4}.
\end{equation}
In particular,
\begin{equation}
\label{eq:lasso-calibrated-cutoff-MSE}
 \limsup_{n\to\infty}\frac1n\E\left[|\widehat W-W_{1-\epsilon,\,\epsilon z_\ast}|^2\right]
 \le C\epsilon^{1/4}.
\end{equation}
The supremum is over deterministic targets and remains outside expectation.
\end{proposition}
\begin{proof}
We compare every target to the fixed reference $V_\epsilon$ through its
nonnegative convexity remainder. The response and support estimates below
are needed only at this fixed positive cutoff.
Equations~\eqref{eq:population-root-second-moment}--\eqref{eq:population-root-correlation},
Proposition~\ref{prop:population-nondegeneracy}, and
\eqref{eq:population-moment-regime} give, for all sufficiently large $n$,
\[
 0<\nu_-\le\nu_\ast\le L_2,\qquad
 0<\beta_-\le\beta_\ast\le L_1,\qquad |\eta_\ast|\le C.
\]
These dimension-uniform bounds use only the population system.
At a fixed reference cutoff, put
$J^\epsilon:=D_hV_\epsilon/\sqrt n$ and
$\kappa^\epsilon:=\E[\tr(\Sigma_xJ^\epsilon)]/n$.
The already proved transport lemma, whose Gaussian input is only
Proposition~\ref{prop:prediction-positive-probes}, gives
\begin{equation}
\label{eq:lasso-cutoff-uniform-moment}
 \sup_{0<\epsilon<\epsilon_1}\limsup_{n\to\infty}
 \left\{\E\left[|V_\epsilon|^2\right]/n+\kappa^\epsilon\right\}<\infty.
\end{equation}
No response bound for the moving target is used.

\paragraph{Uniform moments on a terminal interval.}
For $\delta>0$, square completion gives
\[
\begin{aligned}
 &\frac{\delta\nu_\ast}{2}w^\top\Sigma_xw
       +(\sqrt\delta\,\beta_\ast G-\delta\sqrt n\,\eta_\ast)^\top w\\
 &\quad\ge-\frac1{2\nu_\ast}
       |\beta_\ast\Sigma_x^{-1/2}G
             -\sqrt{\delta n}\,\Sigma_x^{-1/2}\eta_\ast|^2
 \ge-C(|G|^2+n).
\end{aligned}
\]
Comparison with zero and nonnegative losses imply
$\lambda\|V_\delta\|_1\le nL_0+C(|G|^2+n)$.
At zero, $\lambda\|V_0\|_1\le nL_0$. Consequently,
\begin{equation}
\label{eq:lasso-uniform-polynomial-envelope}
 |V_\delta|\le C(n+|G|^2),\qquad 0\le\delta\le1.
\end{equation}
Take $R,r_h,\mathcal L_n$ from
Lemma~\ref{lem:lasso-bounded-root-localization}. Choose a fixed $M$ so
large that $\Pp(|G|>M\sqrt n)\le Ce^{-cn}$, and then a fixed
$\epsilon_1<1/2$ so small that
$\epsilon_1C+\sqrt{\epsilon_1}L_1M\le r_h$.
On the common event
\[
 \mathcal B_n:=\mathcal L_n\cap\{|G|\le M\sqrt n\},
\]
the effective linear field
$\delta\eta_\ast-\sqrt\delta\,\beta_\ast G/\sqrt n$ has norm at most
$r_h$ for all $\delta\le\epsilon_1$. The localization lemma thus gives
$|V_\delta|\le R\sqrt n$ simultaneously on this interval.
The event $\mathcal L_n$ is independent of $G$. Therefore
\eqref{eq:lasso-uniform-polynomial-envelope} gives
\begin{equation}
\label{eq:lasso-terminal-localization-tail}
\begin{aligned}
 &\sup_{0\le\delta\le\epsilon_1}
 \E\left[\frac{|V_\delta|^2}{n}\mathbf1_{\mathcal B_n^c}\right]\\
 &\quad\le \frac Cn\left\{
 \E[(n+|G|^2)^2]\Pp(\mathcal L_n^c)
 +\E[(n+|G|^2)^2\mathbf1_{\{|G|>M\sqrt n\}}]\right\}
 \longrightarrow0.
\end{aligned}
\end{equation}
Indeed, $\E[(n+|G|^2)^2]=O(n^2)$, and Gaussian polynomial tail moments
bound the second term in braces by $Cn^2e^{-cn}$, after enlarging $M$.
It follows that
\[
 \sup_{0\le\delta\le\epsilon_1}\frac{\E\left[|V_\delta|^2\right]}n\le C
\]
for all sufficiently large $n$.

\paragraph{An unconditional bound for the reference remainder.}
Fix $0<\epsilon<\epsilon_1$. For $0\le\delta\le\epsilon$, set
\[
\begin{aligned}
 F_\epsilon(w)&:=\lambda\|w\|_1+
       \sum_{i\in I(1-\epsilon)}L_{y_i}(x_i^\top w/\sqrt n),\\
 Q(w)&:=\frac{\nu_\ast}{2}w^\top\Sigma_xw-\sqrt n\,\eta_\ast^\top w,\\
 A_\delta(w)&:=\sum_{i\in I(1-\delta)\setminus I(1-\epsilon)}
                         L_{y_i}(x_i^\top w/\sqrt n).
\end{aligned}
\]
Then $V_\delta$ minimizes
$F_\epsilon+A_\delta+\delta Q+\sqrt\delta\,\beta_\ast G^\top w$.
Put $d_\delta:=V_\delta-V_\epsilon$ and take the optimality subgradient
$g_\epsilon=-\epsilon\nabla Q(V_\epsilon)-\sqrt\epsilon\,\beta_\ast G$
of $F_\epsilon$ at $V_\epsilon$. Its convexity remainder satisfies
\begin{equation}
\label{eq:lasso-anchor-remainder}
\begin{aligned}
 D_{\epsilon,\delta}
 &:=F_\epsilon(V_\delta)-F_\epsilon(V_\epsilon)
                           -g_\epsilon^\top d_\delta\ge0,\\
 D_{\epsilon,\delta}
 &\le A_\delta(V_\epsilon)
       +(\epsilon-\delta)\nabla Q(V_\epsilon)^\top d_\delta\\
 &\qquad-\frac{\delta\nu_\ast}{2}d_\delta^\top\Sigma_xd_\delta
       +(\sqrt\epsilon-\sqrt\delta)\beta_\ast G^\top d_\delta.
\end{aligned}
\end{equation}
This follows by optimality of $V_\delta$, quadratic expansion of $Q$,
and discarding $-A_\delta(V_\delta)\le0$.

The omitted losses obey
\[
 A_\delta(V_\epsilon)\le
 \sum_{i\notin I(1-\epsilon)}L_{y_i}(x_i^\top V_\epsilon/\sqrt n).
\]
Let $V_\epsilon^{-i}$ be the unforced cavity optimizer with observation
$i$ deleted. It is independent of $(x_i,y_i)$ and $B_i(1-\epsilon)$;
when $B_i(1-\epsilon)=0$, it equals $V_\epsilon$.
The deletion estimate \eqref{eq:nonsmooth-deletion} gives
$\|V_\epsilon-V_\epsilon^{-i}\|_{L^2}\le C/\epsilon$, and hence
\[
 \frac{\|V_\epsilon^{-i}\|_{L^2}}{\sqrt n}
 \le C+\frac{C}{\epsilon\sqrt n}.
\]
Using $L_y(u)\le L_0+L_1|u|$ and $\Sigma_x\preceq KI_p$, before any
restriction to favorable events, we obtain
\[
 \frac1n\E\sum_{i\notin I(1-\epsilon)}L_{y_i}(x_i^\top V_\epsilon/\sqrt n)
 \le\epsilon\left(L_0+L_1\sqrt K
             \frac{\|V_\epsilon^{-i}\|_{L^2}}{\sqrt n}\right)
 \le C\epsilon+\frac C{\sqrt n}.
\]
The same bound holds uniformly for all target $\delta$.

For every deterministic $\delta>0$, $\E[G^\top V_\delta]\le0$.
To see this, condition on the observations and selectors, and take an
independent copy $G'$ of $G$. Subgradient monotonicity gives
\[
 (G-G')^\top\bigl(V_\delta(G)-V_\delta(G')\bigr)\le0.
\]
Taking expectations and using centering and independence proves the sign;
\eqref{eq:lasso-uniform-polynomial-envelope} ensures integrability.
At $\delta=0$, the minimum-norm optimizer is independent of $G$, so the
expectation is zero. This auxiliary copy is used only for the sign; the
actual comparison of $V_\delta$ and $V_\epsilon$ uses the same $G$.
At the fixed reference, Gaussian integration by parts gives
\[
 -\E[G^\top V_\epsilon]
       =\sqrt\epsilon\,\beta_\ast n\kappa^\epsilon.
\]
Keeping the sign, without conditioning on any event depending on $G$,
\[
\begin{aligned}
 &(\sqrt\epsilon-\sqrt\delta)\beta_\ast
               \E[G^\top(V_\delta-V_\epsilon)]\\
 &\quad\le(\sqrt\epsilon-\sqrt\delta)\sqrt\epsilon\,
                 \beta_\ast^2 n\kappa^\epsilon
 \le\epsilon\beta_\ast^2n\kappa^\epsilon.
\end{aligned}
\]
The gradient term in \eqref{eq:lasso-anchor-remainder} is at most
$C\epsilon n$ in expectation by Cauchy--Schwarz and the uniform second
moments. Discard the nonpositive quadratic term. By
\eqref{eq:lasso-cutoff-uniform-moment}, these unconditional estimates give
\begin{equation}
\label{eq:lasso-anchor-remainder-expectation}
 \limsup_{n\to\infty}\sup_{0\le\delta\le\epsilon}
             \frac{\E D_{\epsilon,\delta}}n\le C\epsilon.
\end{equation}

\paragraph{Support, slack, and curvature at the reference.}
Let $s^\epsilon\in\lambda\partial\|V_\epsilon\|_1$ be the actual
penalty subgradient in its optimality condition and put
$S_a:=\{j:|s_j^\epsilon|>\lambda-a\}$, with $0<a<\min\{\lambda_-,1\}$.
Active coordinates belong to $S_a$. For $\sigma:=\sqrt\epsilon\beta_\ast$,
\begin{equation}
\label{eq:lasso-nearly-active-count}
 \E\left[\#\{j:V_{\epsilon,j}=0,\ |s_j^\epsilon|>\lambda-a\}\right]
 \le\frac{2pa}{\sigma\sqrt{2\pi\underline c_x}}.
\end{equation}
Indeed, conditional on the data and $G_{-j}$, the optimizer constrained
by $w_j=0$ is independent of $G_j$. Its $j$th residual is affine with
slope $-\sigma$; the two threshold intervals have total length $2a$,
and $\Var(G_j\mid G_{-j})\ge\underline c_x$.
For $s_0<s_1<1$ from the fixed-cutoff support gap
\eqref{eq:lasso-calibrated-support}, Markov's inequality gives
\begin{equation}
\label{eq:lasso-nearly-active-probability}
 \limsup_{n\to\infty}\Pp\{|S_a|>s_1n\}\le Ca/\sqrt\epsilon.
\end{equation}
Every support and slack estimate concerns $V_\epsilon$ alone.

Apply Lemma~\ref{lem:lasso-ac-curvature} with $s=s_1$ and the common
radius $R$, and denote its retained-row tolerance by
$\epsilon_{\rm ac}>0$. Shrink $\epsilon_1$ below
$\epsilon_{\rm ac}/2$ and the cutoff required for the support gap.
The common favorable event is
\[
 \mathcal H_{n,\epsilon,a}:=\mathcal A_n\cap\mathcal B_n
       \cap\{|S_a|\le s_1n\}
       \cap\{|I(1-\epsilon)^c|\le\epsilon_{\rm ac}n\}.
\]
Here $\mathcal A_n$ is the curvature event. At fixed
$\epsilon<\epsilon_1$, the omitted-row count has the
$\operatorname{Binomial}(n,\epsilon)$ law, so its displayed failure
probability tends to zero, as does $\Pp(\mathcal A_n^c)$.

On $\mathcal H_{n,\epsilon,a}$ the curvature-or-slack argument applies
simultaneously to every target. If
$\|(d_\delta)_{S_a^c}\|_1\le\tau\sqrt n|d_\delta|$, integrate the loss
curvature on the segment $V_\epsilon+ud_\delta$, $0\le u\le1$, using
the fixed retained set $I(1-\epsilon)$ and the common localization ball.
Otherwise the inactive coordinates outside $S_a$ contribute penalty
slack at least $a\|(d_\delta)_{S_a^c}\|_1$. Since
$|d_\delta|\le2R\sqrt n$, both cases yield
\begin{equation}
\label{eq:lasso-calibrated-remainder-lower}
 D_{\epsilon,\delta}\ge
       \min\{\gamma_R/2,a\tau/(2R)\}|d_\delta|^2
       \ge ca|d_\delta|^2.
\end{equation}
The event is uniform over supports, directions and score segments; no
support estimate for $V_\delta$ is used.

Use $D_{\epsilon,\delta}\ge0$ to restrict
\eqref{eq:lasso-anchor-remainder-expectation} to this event. On the rest
of $\mathcal B_n$, $|d_\delta|^2/n\le4R^2$, so curvature and retained-row
failures contribute at most $4R^2$ times their vanishing probabilities;
\eqref{eq:lasso-nearly-active-probability} controls the remaining failure.
Outside $\mathcal B_n$, use
\eqref{eq:lasso-terminal-localization-tail} for both $V_\delta$ and
$V_\epsilon$. Thus for fixed $\epsilon,a$,
\[
 \limsup_{n\to\infty}\sup_{0\le\delta\le\epsilon}
       \frac1n\E\left[|d_\delta|^2\right]
 \le C\left(\frac\epsilon a+\frac a{\sqrt\epsilon}\right).
\]
After the dimension limit, choose $a:=\epsilon^{3/4}$ for sufficiently
small $\epsilon$. This proves \eqref{eq:lasso-calibrated-terminal-MSE};
setting $\delta=0$ gives \eqref{eq:lasso-calibrated-cutoff-MSE}.
\end{proof}

\begin{proof}[Completion of the proof of Theorem~\ref{thm:prediction-directions}, pure-Lasso case]
For a fixed $\epsilon>0$, Equations~\eqref{eq:population-root-second-moment}--\eqref{eq:population-root-correlation},
Proposition~\ref{prop:population-nondegeneracy}, and
\eqref{eq:population-moment-regime} give dimension-uniform positive lower
bounds on $\nu_\ast,\beta_\ast$ and upper bounds on every component of
$z_\ast$. Thus the calibrated path over $[0,1-\epsilon]$ lies in a
fixed positive-probe set, and
Proposition~\ref{prop:prediction-positive-probes} gives
\[
 \sup_{0\le t\le1-\epsilon}\mathcal G(W_{t,(1-t)z_\ast})\longrightarrow0.
\]
For $t\in[1-\epsilon,1]$, put $\delta=1-t$ and use the coupling above.
Proposition~\ref{prop:gc-transfer}, always using the actual optimizer
moments, yields
\[
 \sup_{0\le\delta\le\epsilon}\mathcal G(V_\delta)
 \le\mathcal G(V_\epsilon)+C
      \sup_{0\le\delta\le\epsilon}
          \left(\frac{\E|V_\delta-V_\epsilon|^2}{n}\right)^{1/2}.
\]
By \eqref{eq:lasso-calibrated-terminal-MSE},
\[
 \limsup_{n\to\infty}\sup_{1-\epsilon\le t\le1}\mathcal G(W_t)
 \le C\epsilon^{1/8},\qquad W_1:=\widehat W.
\]
Combine the two intervals and let $\epsilon\downarrow0$ after the
dimension limit. This proves the pure-Lasso case and completes
Theorem~\ref{thm:prediction-directions}. The conclusion is
qualitative; no uniform rate in $n$ at vanishing probes is asserted.
\end{proof}

\begin{proof}[Proof of Theorem~\ref{thm:endpoint}]
For either penalty branch, put
$D_\epsilon:=\|\widehat W-W_{1-\epsilon,\,\epsilon z_\ast}\|_{L^2}/\sqrt n$.
Propositions~\ref{prop:strongly-convex-calibrated-cutoff}
and~\ref{prop:lasso-calibrated-cutoff-MSE} give
$\lim_{\epsilon\downarrow0}\limsup_n D_\epsilon=0$ in their respective branches.
For the cutoff moments in Lemma~\ref{lem:lasso-calibrated-transport},
the coupling bounds \eqref{eq:coupling-mean-standard-deviation} give
\[
 \left|\frac{\E\widehat W}{\sqrt n}-\mu^\epsilon\right|\le D_\epsilon,
 \qquad
 \left|\sqrt{\frac{\tr(\Sigma_xC_{\widehat W})}{n}}
              -\alpha^\epsilon\right|\le\sqrt K D_\epsilon.
\]
By \eqref{eq:lasso-calibrated-moments},
$|\mu^\epsilon-\mu_\ast|+|\alpha^\epsilon-\alpha_\ast|\to0$
at each fixed $\epsilon$. Taking $n\to\infty$ and then
$\epsilon\downarrow0$ proves convergence of the normalized mean and
prediction standard deviation. Since \eqref{eq:population-moment-regime}
and \eqref{eq:spectral-bounds} bound $\alpha_\ast$, this also gives
variance convergence and proves \eqref{eq:main-moment-comparison}.
\end{proof}

\section{Prediction laws and classification}
\label{app:prediction-laws}
Theorem~\ref{thm:prediction-directions} gives the empirical
score comparison, and Theorem~\ref{thm:endpoint} identifies its moments.
Cutoff coupling also transfers the normalized covariance bound.
Classification then follows by smoothing
the threshold, using $\liminf_n\alpha_\ast^2>0$ from
\eqref{eq:population-moment-regime}.

\begin{proof}[Proof of Corollary~\ref{cor:conditional-prediction}]
Write $\widehat\mu_n:=\E[\widehat W]/\sqrt n$,
$\widehat\alpha_n^2:=\tr(\Sigma_xC_{\widehat W})/n$, and
$\widehat\alpha_n:=\sqrt{\widehat\alpha_n^2}$. Coupling the comparison
Gaussians with a common standard normal gives
\begin{equation}
\label{eq:conditional-prediction-bound}
\begin{aligned}
 \mathfrak W&:=\E_x\left[\Wass\!\left(
 \Law(x^\top\widehat W/\sqrt n\mid x),
 N(x^\top\mu_\ast,\alpha_\ast^2)\right)\right]\\
 &\le\mathcal G(\widehat W)
       +\sqrt K|\widehat\mu_n-\mu_\ast|
       +\sqrt{2/\pi}\,|\widehat\alpha_n-\alpha_\ast|
 \longrightarrow0.
\end{aligned}
\end{equation}
The first term vanishes by
Theorem~\ref{thm:prediction-directions} at $t=1$; the last
two vanish by Theorem~\ref{thm:endpoint} and
$|\sqrt u-\sqrt v|\le\sqrt{|u-v|}$. Thus
\begin{equation}
\label{eq:conditional-prediction}
 \E_x\left[\Wass\!\left(
 \Law(x^\top\widehat W/\sqrt n\mid x),
 N(x^\top\mu_\ast,\alpha_\ast^2)\right)\right]\longrightarrow0.
\end{equation}
Independence from the entire fresh pair makes the score law unchanged if
the fresh label is also fixed, even when $y$ depends on $x$.
For each fixed $(x,y)$, Wasserstein duality bounds the discrepancy for
$f_y$ by its common Lipschitz constant $\Lambda$ times the score-law
distance conditional on $x$. Averaging over $(x,y)$ gives
$\Lambda\mathfrak W$, proving \eqref{eq:main-risk}.
\end{proof}

\begin{proof}[Proof of Corollary~\ref{cor:classification}]
The comparison Gaussian has density at most
$(\sqrt{2\pi}\,\underline\alpha)^{-1}$.
Upper and lower ramps of width $\ell_{\rm ramp}>0$ around a half-line
boundary give, for any scalar laws $P,Q$ with this density bound on $Q$,
\begin{equation}
\label{eq:wasserstein-halfline}
 |P(I)-Q(I)|\le\Wass(P,Q)/\ell_{\rm ramp}
                         +\ell_{\rm ramp}/(\sqrt{2\pi}\,\underline\alpha),
\end{equation}
for every open or closed half-line $I$; continuity of $Q$ handles either
endpoint convention. With $\operatorname{sign}_+(0)=1$, the error half-lines
are $(-\infty,0)$ for $y=1$ and $[0,\infty)$ for $y=-1$.
Condition on $(x,y)$ to select the error half-line. By independence from
the fresh pair, the score laws still depend only on $x$, so averaging
bounds the error discrepancy by
$\mathfrak W/\ell_{\rm ramp}+\ell_{\rm ramp}/(\sqrt{2\pi}\,\underline\alpha)$.
Corollary~\ref{cor:conditional-prediction}, followed by
$\ell_{\rm ramp}\downarrow0$, proves convergence. The Gaussian conditional
error is $\Phi_{\rm N}(-yx^\top\mu_\ast/\alpha_\ast)$, giving
\eqref{eq:classification-risk}, including possible fitted-score ties.
\end{proof}

\begin{remark}[Covariance bound at the empirical endpoint]
\label{rem:boundary-covariance}
Under the assumptions of Theorem~\ref{thm:endpoint},
$\|C_{\widehat W}\|_{\hs}/n\to0$.
Indeed, couple $W_{1-\epsilon,\epsilon z_\ast}$ and $\widehat W$ through
the same observations, and put
$D_\epsilon:=\|\widehat W-W_{1-\epsilon,\epsilon z_\ast}\|_{L^2}/\sqrt n$.
Lemma~\ref{lem:lasso-calibrated-transport} and the cutoff bounds in
Propositions~\ref{prop:strongly-convex-calibrated-cutoff} and
\ref{prop:lasso-calibrated-cutoff-MSE}, in their respective penalty branches,
bound the normalized second moments uniformly in the ordered limit.
Thus \eqref{eq:coupling-moments} and Lemma~\ref{lem:pred-proof-concentration} give
\[
 \limsup_{n\to\infty}\frac{\|C_{\widehat W}\|_{\hs}}n
 \le C\limsup_{n\to\infty}D_\epsilon\longrightarrow0
 \qquad(\epsilon\downarrow0).
\]
\end{remark}

\section{Penalty formulas and a directional counterexample}
\label{app:penalty-specializations}
We first specialize the surrogate to quadratic penalties, then to Lasso
and elastic net, and check consistency with Gaussian linear-model
formulas. The last subsection distinguishes Gaussian fresh scores from
Gaussian coefficient projections.
The penalty formulas use a finite population solution
$(\mu_\ast,\alpha_\ast,\kappa_\ast,\beta_\ast,\nu_\ast,\eta_\ast)$,
with $\beta_\ast,\nu_\ast>0$, and the surrogate $W_\ast(G)$.
Write $U_\ast,\zeta_\ast,\xi_\ast$ for the quantities in
\eqref{eq:losscoeff} evaluated at this solution.
Throughout, $D_\eta W_\ast$ denotes differentiation with respect to the
linear field $\eta$, evaluated at $\eta_\ast$, with all other coefficients
fixed; it does not differentiate the population system.
Transfer of prediction moments requires all assumptions of
Theorem~\ref{thm:endpoint}; support transfer remains open.
The squared-loss calculations are algebraic consistency checks, since
squared loss does not satisfy Assumption~\ref{ass:analytic}.

\subsection{Quadratic penalties}
\label{app:quadratic-specialization}
\label{sec:applications}
For $\rho(w):=w^\top Hw/2$, with $H\succeq0$, put
$Q_\ast:=(\nu_\ast\Sigma_x+H)^{-1}$.
Ridge regularization corresponds to $H=\lambda_2I_p$.
Completing the square gives
\begin{equation}
\label{eq:quadratic-effective}
\begin{gathered}
 W_\ast(G)=\sqrt n\,Q_\ast\eta_\ast-\beta_\ast Q_\ast G,\\
 \E_G\left[W_\ast\right]=\sqrt n\,Q_\ast\eta_\ast,\qquad
 C_{W_\ast}=\beta_\ast^2Q_\ast\Sigma_xQ_\ast,
\end{gathered}
\end{equation}
and the population system reduces to
\begin{equation}
\label{eq:quadratic-closure}
\begin{gathered}
 H\mu_\ast+\E\left[x\xi_\ast\right]=0,\qquad
 \kappa_\ast=\frac1n\tr(\Sigma_xQ_\ast),\\
 \alpha_\ast^2=\frac{\beta_\ast^2}{n}
                  \tr(\Sigma_xQ_\ast\Sigma_xQ_\ast).
\end{gathered}
\end{equation}
These identities hold for $H\succeq0$ whenever a population solution
exists. For $H\succ0$, existence follows from
Theorem~\ref{thm:population-uniqueness} under its loss assumptions.
The strongly convex branch of Theorem~\ref{thm:endpoint} additionally
requires a uniform positive lower bound on the eigenvalues of $H$.

\begin{proposition}[Quadratic surrogate and loss consistency]
\label{prop:quadratic-surrogate-consistency}
For $\rho(w):=w^\top Hw/2$, $H\succeq0$, every finite population solution
satisfies \eqref{eq:quadratic-effective}--\eqref{eq:quadratic-closure}.
For squared loss $L_y(u):=(u-y)^2/2$ and
$y:=x^\top\theta_0+\varepsilon$, with deterministic $\theta_0$ and
independent centered noise of variance $\sigma^2$, put
\[
 T_\ast:=\frac1n\tr(\Sigma_xQ_\ast\Sigma_xQ_\ast),\qquad
 \Delta_\ast:=(\mu_\ast-\theta_0)^\top
                    \Sigma_x(\mu_\ast-\theta_0).
\]
If $1-\nu_\ast^2T_\ast>0$, then
\begin{equation}
\label{eq:quadratic-risk}
\begin{gathered}
 \nu_\ast=\frac1{1+\kappa_\ast},\qquad
 \mu_\ast=\nu_\ast Q_\ast\Sigma_x\theta_0,\\
 \frac12\E\left[(x^\top\mu_\ast+\alpha_\ast Z-y)^2\right]
 =\frac{\sigma^2+\Delta_\ast}{2(1-\nu_\ast^2T_\ast)},
\end{gathered}
\end{equation}
where $Z\sim N(0,1)$ is independent of $(x,y)$.
\end{proposition}
\begin{proof}
Completing the square gives \eqref{eq:quadratic-effective} and
$D_\eta W_\ast=\sqrt n\,Q_\ast$.
Substitution into the mean, variance, response, and field equations gives
\eqref{eq:quadratic-closure}.
The relation $U_\ast-\zeta_\ast=\kappa_\ast\xi_\ast$ and Gaussian
integration by parts give
\[
 \beta_\ast^2=\frac{\E\left[(U_\ast-\zeta_\ast)^2\right]}{\kappa_\ast^2},\qquad
 \nu_\ast\alpha_\ast=
 \frac{\E\left[Z(U_\ast-\zeta_\ast)\right]}{\kappa_\ast}.
\]
Together with \eqref{eq:quadratic-closure}, these recover the quadratic
equations of \citet[Theorem~4.4]{Yaakoubi2026} for centered features.
For squared loss,
\[
 \xi_\ast=\frac{U_\ast-y}{1+\kappa_\ast},\qquad
 \eta_\ast=\nu_\ast\Sigma_x\theta_0,\qquad
 \beta_\ast^2=\nu_\ast^2
       (\sigma^2+\Delta_\ast+\alpha_\ast^2).
\]
Combining these identities with
$\alpha_\ast^2=\beta_\ast^2T_\ast$ proves
\eqref{eq:quadratic-risk}.
For centered features, the risk formula agrees with
\citet[Corollary~5.1]{Yaakoubi2026}.
\end{proof}

\subsection{Lasso and elastic-net formulas}
\label{app:sparse-specialization}
We prove the formulas in Example~\ref{sec:sparse-applications} for
$\rho(w)=\lambda_1\|w\|_1+\lambda_2|w|^2/2$, with $\lambda_1>0$
and $\lambda_2\ge0$.

\begin{proposition}[Diagonal sparse surrogate and moment equations]
\label{prop:diagonal-sparse-formulas}
For the above penalty and diagonal
$\Sigma_x=\operatorname{diag}(c_1,\ldots,c_p)\succ0$, the surrogate
$W_\ast$ and its moments satisfy
\eqref{eq:soft}--\eqref{eq:diagonal-sparse-moments}.
\end{proposition}
\begin{proof}
Write $G_j=\sqrt{c_j}Z_j$, with independent $Z_j\sim N(0,1)$, and put
the surrogate optimality condition in each coordinate in scalar form:
$q_{\ast,j}-\Gamma_{\ast,j}W_{\ast,j}\in\lambda_1\partial|W_{\ast,j}|$, hence
\[
 W_{\ast,j}=\frac{\soft(q_{\ast,j},\lambda_1)}{\Gamma_{\ast,j}}.
\]
Taking means and variances gives the equations for $\mu_\ast$ and
$\alpha_\ast^2$.
Almost everywhere,
$\partial_{\eta_j}W_{\ast,j}
 =\sqrt n\,\Gamma_{\ast,j}^{-1}\ind_{\{|q_{\ast,j}|>\lambda_1\}}$.
Gaussian integration by parts yields
\begin{equation}
\label{eq:penalty-response-identity}
 \kappa_\ast
 =\frac1{n\sqrt n}\E_G\left[\tr(\Sigma_xD_\eta W_\ast)\right]
 =-\frac1{\beta_\ast n}\E_G\left[G^\top W_\ast\right],
\end{equation}
which gives the response equation in
\eqref{eq:diagonal-sparse-moments}.
\end{proof}

\begin{corollary}[Correlated sparse surrogate response]
\label{cor:correlated-sparse-response}
For the same penalty and a general $\Sigma_x\succ0$, put
$S_\ast:=\operatorname{supp}(W_\ast(G))$, and let $P_{S_\ast}$ embed its
coordinates into $\R^p$.
Almost surely,
\begin{equation}
\label{eq:active}
 D_\eta W_\ast=\sqrt n\,
 P_{S_\ast}\bigl(\nu_\ast(\Sigma_x)_{S_\ast S_\ast}
           +\lambda_2 I_{|S_\ast|}\bigr)^{-1}P_{S_\ast}^\top.
\end{equation}
For pure Lasso,
$\nu_\ast\kappa_\ast=\E_G\left[|S_\ast|\right]/n$.
\end{corollary}
\begin{proof}
The optimizer is piecewise affine in
$\sqrt n\,\eta_\ast-\beta_\ast G$, whose nondegenerate Gaussian law avoids
its finitely many activation boundaries.
Lemma~\ref{lem:active-coordinate-response} gives \eqref{eq:active}, with
zero derivative when $S_\ast$ is empty.
For $\lambda_2=0$,
$\tr(\Sigma_xD_\eta W_\ast)=\sqrt n\,|S_\ast|/\nu_\ast$;
\eqref{eq:penalty-response-identity} proves the support-size identity.
At $\lambda_1=0$, \eqref{eq:active} follows directly from
\eqref{eq:quadratic-effective} with $H=\lambda_2 I_p$, and the surrogate
has full support almost surely.
The identity is related to active-set degrees of freedom
\citep{TibshiraniTaylor2012}.
\end{proof}

\subsection{Gaussian linear-model consistency}
\begin{proposition}[Gaussian linear-model consistency]
\label{prop:gaussian-linear-consistency}
Let $L_y(u):=(u-y)^2/2$, $x\sim N(0,I_p)$, and
$y:=x^\top\theta_0+\varepsilon$, with deterministic $\theta_0$ and
independent centered noise of variance $\sigma^2$.
For $\rho(w):=\lambda_1\|w\|_1+\lambda_2\|w\|_2^2/2$, with
$\lambda_1>0$ and $\lambda_2\ge0$, put
\[
 \tau_\ast:=\frac{\beta_\ast}{\nu_\ast},\qquad
 \mathsf S_\ast(z):=
 \frac{\soft(z,\lambda_1/\nu_\ast)}{1+\lambda_2/\nu_\ast}.
\]
Every finite population solution satisfies
\begin{equation}
\label{eq:sparse-linear-bias-variance}
\begin{gathered}
 \tau_\ast^2=\sigma^2+\frac1n\sum_{j=1}^p
       \E_G\left[(W_{\ast,j}-\sqrt n\,\theta_{0,j})^2\right],\\
 1-\nu_\ast=\frac1n\sum_{j=1}^p
       \E_Z\left[\mathsf S_\ast'(\sqrt n\,\theta_{0,j}+\tau_\ast Z)\right].
\end{gathered}
\end{equation}
Suppose $\lambda_1,\lambda_2$ are fixed,
$\nu_\ast\to\nu_\infty>0$, $\tau_\ast\to\tau_\infty>0$, $p/n\to c$,
and the empirical law of $\sqrt n\,\theta_{0,j}$ converges with second
moments to the law of a random variable $B$.
Writing
$\mathsf S_\infty(z):=
\soft(z,\lambda_1/\nu_\infty)/(1+\lambda_2/\nu_\infty)$, one obtains
\begin{equation}
\label{eq:sparse-amp-comparison}
\begin{gathered}
 \tau_\infty^2=\sigma^2+
       c\E\left[(\mathsf S_\infty(B+\tau_\infty Z)-B)^2\right],\\
 \nu_\infty=1-c\E\left[\mathsf S_\infty'(B+\tau_\infty Z)\right],
\end{gathered}
\end{equation}
where $Z\sim N(0,1)$ is independent of $B$.
\end{proposition}
\begin{proof}
For squared loss,
$\xi_\ast=(U_\ast-y)/(1+\kappa_\ast)$ and
$\nu_\ast=(1+\kappa_\ast)^{-1}$.
Since $\E\left[x\xi_\ast\right]=\nu_\ast(\mu_\ast-\theta_0)$,
$\eta_\ast=\nu_\ast\theta_0$.
The soft-threshold formula and Gaussian symmetry give
\[
 W_{\ast,j}\stackrel{\mathrm{law}}{=}
 \mathsf S_\ast(\sqrt n\,\theta_{0,j}+\tau_\ast Z).
\]
Expanding $\E\left[(U_\ast-y)^2\right]$ and using
$\mu_\ast=\E_G\left[W_\ast\right]/\sqrt n$ and
$\alpha_\ast^2=\tr(C_{W_\ast})/n$ gives the first equation in
\eqref{eq:sparse-linear-bias-variance}.
Equation~\eqref{eq:penalty-response-identity} gives
\[
 \nu_\ast\kappa_\ast=\frac1n\sum_{j=1}^p
 \E_Z\left[\mathsf S_\ast'(\sqrt n\,\theta_{0,j}+\tau_\ast Z)\right],
\]
and $\nu_\ast\kappa_\ast=1-\nu_\ast$ gives the second equation.
Gaussian averaging makes the response continuous and bounded, while
the squared error has quadratic growth. The assumed convergence with
second moments therefore gives \eqref{eq:sparse-amp-comparison}.
For pure Lasso, these are the state-evolution and threshold-calibration
formulas of \citet{BayatiMontanari2012}, with sample-to-dimension ratio
$1/c$; the elastic-net map is the corresponding specialization of
\citet{Thrampoulidis2018}.
\end{proof}

\subsection{Why Gaussianity in every deterministic direction fails}
\label{app:directional-counterexample}
The following example distinguishes the conditional fresh-score
comparison \eqref{eq:prediction-directions} from a Gaussian approximation
of every fixed coefficient projection, at the same positive interpolation
time. Its data are Gaussian, so the obstruction is not a lack of
Poincar\'e concentration.

\begin{proposition}[A non-Gaussian coordinate under Gaussian data]
\label{prop:coordinate-nongaussian}
Let $p=n$, let $x_i\sim N(0,I_n)$ be independent, and take $y_i\equiv1$,
binary logistic loss $L_y(u):=\log(1+e^{-yu})$, and pure Lasso
$\rho(w):=\lambda_1\|w\|_1$ with fixed $\lambda_1>0$ and $\lambda_2=0$.
Set $W_t:=W_{t,1,1/2,0}$ as in \eqref{eq:interior-estimator}.
The data law satisfies a dimension-independent Poincar\'e inequality,
and there are constants $t_0\in(0,1)$ and $c_0>0$, independent of $n$,
such that, for every coordinate $j$,
\begin{equation}
\label{eq:poincare-coordinate-obstruction}
 \Wass\!\left(\Law((W_{t_0})_j),
 N\!\left(\E\left[(W_{t_0})_j\right],\Var((W_{t_0})_j)\right)\right)
 \ge c_0/2.
\end{equation}
For the same $t_0$ and an independent fresh $x\sim N(0,I_n)$,
\begin{equation}
\label{eq:positive-time-fresh-score}
 \E_x\!\left[\Wass\!\left(
 \Law(x^\top W_{t_0}/\sqrt n\mid x),
 N\!\left(0,\Var((W_{t_0})_j)\right)
 \right)\right]\longrightarrow0.
\end{equation}
The variance is the same for every coordinate but may depend on $n$.
\end{proposition}
\begin{proof}
Here $\Sigma_x=I_n$.
The Gaussian Poincar\'e inequality gives
$\Var(F(x,1))\le\E_x\left[|\nabla_xF(x,1)|^2\right]$.
Moreover, $\Var(x^\top Ax)=2\|A\|_{\hs}^2$ for symmetric $A$, so
\eqref{eq:feature-quadratic-concentration} holds.
At $t=0$, the positive quadratic probe gives
\[
 \cE_{0,1,1/2,0}(w)=\lambda_1\|w\|_1+
                         \tfrac12\|w\|_2^2+G^\top w,\qquad
 (W_0)_j=\soft(-G_j,\lambda_1),\quad G\sim N(0,I_n).
\]
The coordinates of $W_0$ are independent and symmetric, with
$\Var((W_0)_j)=\E[\soft(-Z,\lambda_1)^2]>0$ for $Z\sim N(0,1)$.
Their common law has an atom at zero of mass $\Pp(|Z|\le\lambda_1)>0$,
whereas its matching Gaussian has no atoms. Since both laws have finite
first moments and are independent of $n$,
\begin{equation}
\label{eq:poincare-zero-time-gap}
 c_0:=\Wass(\Law((W_0)_j),N(0,\Var((W_0)_j)))>0
\end{equation}
is a fixed positive constant.

Couple $W_t$ and $W_0$ using the same Gaussian probe.
Strong monotonicity gives
\[
 |W_t-W_0|\le
 \left|\frac1{\sqrt n}\sum_{i=1}^n
 B_i(t)x_iL_{y_i}'(x_i^\top W_0/\sqrt n)\right|.
\]
Conditional on $W_0$, these summands are independent.
Before the factor $1/\sqrt n$, each has second moment at most $tn$ and
mean of norm at most $t$: for every deterministic unit vector $e$,
\[
 \left|\E_{x,y}\left[(e^\top x)L_y'(x^\top W_0/\sqrt n)\mid W_0\right]\right|
 \le (\E_x\left[(e^\top x)^2\right])^{1/2}=1.
\]
Consequently,
\begin{equation}
\label{eq:poincare-small-time-coupling}
 \frac1n\E\left[|W_t-W_0|^2\right]
 \le t+\frac{n-1}{n}t^2\le t+t^2.
\end{equation}
Permutation symmetry gives
$\|(W_t)_j-(W_0)_j\|_{L^2}\le\sqrt{t+t^2}$, and coordinate-sign symmetry
gives $\E\left[(W_t)_j\right]=\E\left[(W_0)_j\right]=0$.
Coupling their moment-matched Gaussians with one standard normal and
using the reverse triangle inequality in $L^2$ yields
\[
 \Wass(\Law((W_t)_j),N(0,\Var((W_t)_j)))
 \ge c_0-(1+\sqrt{2/\pi})\sqrt{t+t^2}.
\]
Choose
\[
 t_0:=\min\{1/4,c_0^2/[16(1+\sqrt{2/\pi})^2]\}\in(0,1/4].
\]
This choice is independent of $n$, makes the subtracted term at most
$c_0/2$, and lies in a fixed positive-probe interior domain.
It proves \eqref{eq:poincare-coordinate-obstruction}.

At this same $t_0$, coordinate-sign and permutation symmetry give
\[
 \E[W_{t_0}]=0,\qquad
 C_{W_{t_0}}=\Var((W_{t_0})_j)I_n,\qquad
 \frac1n\tr(\Sigma_xC_{W_{t_0}})=\Var((W_{t_0})_j).
\]
The Gaussian data with constant labels and logistic Lasso satisfy
Assumptions~\ref{ass:analytic} and~\ref{ass:feature-concentration}, with
$p/n=1$ and $\Sigma_x=I_n$. Applying
Proposition~\ref{prop:prediction-positive-probes} on a fixed positive-probe
domain containing $(t_0,1,1/2,0)$ proves
\eqref{eq:positive-time-fresh-score}. Thus the conditional fresh-score
Gaussian approximation and the coordinate obstruction hold for the same
optimizer $W_{t_0}$.
\end{proof}

\section{Prediction experiment details}
\label{app:prediction-experiments}
\subsection{Synthetic data and conditional Gaussian comparison}
\label{app:synthetic-prediction-experiments}
\paragraph{Model.}
Set $n=200k$, $p=400k$, with $10k$ independent size-ten blocks and
$300k$ independent standard Gaussian singleton coordinates. Within a block,
\[
 x_{bj}=\sqrt{0.95}\frac{H_b+8S_b}{\sqrt{65}}+\sqrt{0.05}\,E_{bj},
\]
where $H_b,E_{bj}$ are independent standard Gaussians and $S_b$ are
independent symmetric signs. Replacing $(H_b+8S_b)/\sqrt{65}$ by $H_b$
gives the covariance-matched Gaussian design. The first $2k$ blocks and
$5k$ singletons have teacher coefficients $\sqrt{0.9/(191k)}$ and
$\sqrt{0.1/(5k)}$, respectively, so
$\theta_0^\top\Sigma_x\theta_0=1$. We generate
$y=\tanh(x^\top\theta_0+\varepsilon)$ with independent
$\varepsilon\sim N(0,0.01^2)$ and use pseudo-Huber loss with scale $0.2$.
Figure~\ref{fig:synthetic-prediction}(a) uses $k=1,\ldots,10$ and
$(\lambda_1,\lambda_2)=(0.1,0.1),(0.2,0.1),(0.2,0)$;
panel~(b) uses $k=1$. Penalties refer to the sum-loss objective in $W$;
the mean-loss $\ell_1$ coefficient in $\theta=W/\sqrt n$ is
$\lambda_1/\sqrt n$.

\paragraph{Conditional Wasserstein distances.}
At each dimension and penalty, recalibrate the population solution and
fit $W_{t,z(t)}$ at $t\in\mathcal T=\{0,0.1,\ldots,1\}$, with
$z(t)=(1-t)z_\ast$. Hold 128 independent fresh features fixed across
training replications. Independent cohorts of 2,000 fits estimate the
mean $\widehat\mu$ and covariance $\widehat C$ of $\theta$, and the
conditional score laws $\widehat P_j$, respectively. Panel (a) plots
\[
 \max_{t\in\mathcal T}\frac1{128\widehat\alpha}
 \sum_{j=1}^{128}\Wass\!\left(\widehat P_j,
 N((x^{(j)})^\top\widehat\mu,\widehat\alpha^2)\right),
 \qquad \widehat\alpha^2=\tr(\Sigma_x\widehat C),
\]
with moments and laws evaluated separately at each $t$.
Each dotted reference is the maximum over $t$ of the mean standardized
error from 1,024 Gaussian simulations with 2,000 scores, including
independent mean-estimation error and bootstrap trace-variance uncertainty.
These measure finite-sampling Wasserstein error. The full calibrated grid,
including the surrogate and empirical endpoints, illustrates
Theorem~\ref{thm:prediction-directions} when its additional
hypotheses hold. The simulations do not establish these hypotheses or the
asymptotic convergence.

\paragraph{Conditional-score and coordinate histograms.}
Figure~\ref{fig:synthetic-prediction}(b) uses 3,000 independently resampled
training datasets per column, holding the fresh test vector fixed.
The top row shows $x^\top W/\sqrt n\mid x$, and the bottom row shows
$e_1^\top W$ from the same fitted vectors. Coordinate~1 is a predetermined
signal coordinate; zero coefficients are retained. The first two columns
use the original design with pure $\ell_2$ regularization
$(\lambda_1,\lambda_2)=(0,0.1)$ and pure $\ell_1$ regularization
$(0.2,0)$, respectively. The third uses $(0.2,0)$ and replaces the
Gaussian feature primitives by interquartile-range-matched Student
$t_{1.5}$ variables. This heavy-tailed variant has infinite feature
variance and does not satisfy the Poincar\'e assumption. The gray Gaussian
densities use each histogram's sample mean and variance. Solid and dashed
curves mark close agreement and visible departures, respectively.

\begin{figure}[!htb]
 \centering
 \includegraphics[width=0.94\linewidth]{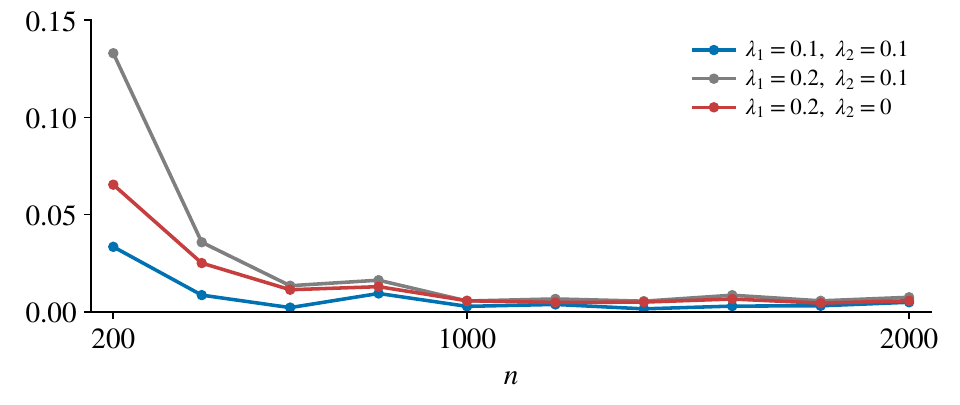}
 \revisioncolor
 \caption[Relative prediction-variance error]{Relative prediction-variance error.
 The vertical axis shows
 $\max_{t\in\mathcal T}|\widehat\alpha_t^2-\alpha_\ast^2|/\alpha_\ast^2$,
 with $\mathcal T=\{0,0.1,\ldots,1\}$ and $z(t)=(1-t)z_\ast$.
 Here $\widehat\alpha_t^2=\tr(\Sigma_x\widehat C_t)$ uses the same
 2,000-fit moment cohorts as Figure~\ref{fig:synthetic-prediction}(a),
 where $\widehat C_t$ is the sample covariance of $W_{t,z(t)}/\sqrt n$;
 $\alpha_\ast^2$ is calibrated separately for each dimension and penalty.}
 \label{fig:synthetic-alpha-error}
\end{figure}

\paragraph{Prediction variance.}
Figure~\ref{fig:synthetic-alpha-error} compares the estimated prediction
variance along the same calibrated path with its population target.
The maximum relative discrepancy decreases overall with dimension for
all three penalties.

\FloatBarrier

\subsection{Diamonds preprocessing and prediction comparisons}
\label{app:diamonds-experiments}
Figure~\ref{fig:diamonds-prediction} compares classification errors in (a)
and pooled and centered fresh-score laws at $\lambda_1=0.5$ in (b,c).
\paragraph{Frozen data and feature map.}
We use the \href{https://ggplot2.tidyverse.org/reference/diamonds.html}{\texttt{diamonds}
dataset distributed with \texttt{ggplot2}}, containing 53,940 observations.
Removing exact duplicate rows and observations whose length, width, or
depth lies outside $(0,20)$ leaves 53,772 rows. A fixed permutation assigns the first $M=30{,}000$ rows to the reference population.
The label is $y_i=2\mathbf1_{\{\mathrm{price}_i>3000\}}-1$; price is
excluded from the inputs. All preprocessing is fitted on these reference
rows and frozen. The six numerical inputs (carat, depth percentage, table
percentage, length, width, and depth) are log-transformed, clipped at the
reference 1\% and 99\% quantiles, and standardized. Cubic splines with
ten nominal quantile knots and categorical indicators for cut, color,
and clarity form an 83-dimensional dictionary. Principal-component
projection and whitening retain 64 coordinates, denoted $s_i$.

Put $m=M^{-1}\sum_i s_i$, $C=M^{-1}\sum_i(s_i-m)(s_i-m)^\top$, and
$B=(C+\tau^2I_{64})^{-1/2}$, with $\tau=0.5$. Every observation is an
independent draw from
\begin{equation}
\label{eq:diamonds-smoothed-law}
 \begin{gathered}
 x=\bigl(B(s_I-m+\tau g),\ h,\ 1\bigr),\qquad y=y_I,\\
 I\sim\operatorname{Unif}\{1,\ldots,M\},\qquad
 g\sim N(0,I_{64}),\qquad h\sim N(0,I_{p-65}),
 \end{gathered}
\end{equation}
where $I,g,h$ are independent. The label is unchanged when noise is
added. The last coordinate is a penalized constant. Thus
$\Sigma_x=\E[xx^\top]=I_p$, whereas
$\operatorname{Cov}(x)=\operatorname{diag}(I_{p-1},0)$.
Figure~\ref{fig:diamonds-prediction} uses $(n,p)=(512,256)$, with 191
independent noise coordinates. This defines a semisynthetic population.

\paragraph{ERM and risk evaluation.}
Writing $\theta=W/\sqrt n$, each pure-Lasso training fit solves
\begin{equation}
\label{eq:diamonds-lasso-erm}
 \widehat\theta=\argmin_\theta\left\{
 \frac1n\sum_{i=1}^n\log(1+e^{-y_i x_i^\top\theta})
 +\frac{\lambda_1}{\sqrt n}\|\theta\|_1\right\}.
\end{equation}
The 1,000 independent training samples are paired across
$\lambda_1=0.1,0.15,\ldots,0.8$. For each fitted coefficient vector,
classification error is computed by averaging over all reference
centers and analytically integrating the fresh Gaussian noise in
\eqref{eq:diamonds-smoothed-law}. Thus the empirical markers average
population risks over training randomness; they do not use a finite
test-set estimate. Error bars in Figure~\ref{fig:diamonds-prediction}(a)
show $\pm1.96$ training Monte Carlo standard errors. Our curve evaluates
$\E_{x,y}\Phi(-yx^\top\mu_\ast/\alpha_\ast)$ using the calibrated
population equations, without fitting optimizer moments to the markers
or histograms. Its maximum absolute discrepancy from the displayed ERM
means is $0.044$ percentage points.

\paragraph{Figure conventions.}
In Figure~\ref{fig:diamonds-prediction}, blue denotes our pure-Lasso
prediction, red the smooth-penalty equations, and yellow the
Gaussian-mixture equations. Solid/dashed curves use
$\varepsilon=0.01/0.05$ (red) and $K=1/10$ components per class (yellow).
In panel~(a), markers denote pure Lasso on the original population (blue),
smoothed Lasso with $\varepsilon=0.01$ (red), and pure Lasso on the fitted
$K=1$ mixture (yellow). Panels~(b,c) use 2,000 independent
training/fresh-feature pairs per method: the top row compares pure and
smoothed Lasso, and the bottom row compares the original population and
its fitted mixture. Curves give the corresponding predictions.

\paragraph{Comparison equations.}
The smooth-penalty curves use the equations of
\citet[Theorem~4.4]{Yaakoubi2026} with the smooth penalty in the
mean-loss parameterization $\theta=W/\sqrt n$:
\[
 \rho_\varepsilon(\theta)=\frac{\lambda_1}{\sqrt n}
 \sum_j\bigl(\sqrt{\theta_j^2+\varepsilon^2}-\varepsilon\bigr),
 \qquad \varepsilon\in\{0.05,0.01\}.
\]
Their implementation uses the Hessian at zero,
$\nabla^2\rho_\varepsilon(0)=\lambda_1/(\sqrt n\varepsilon)I_p$,
and the centered feature covariance
$\operatorname{diag}(I_{p-1},0)$, with the penalized constant handled in
the mean equation. These are extrapolations: the smooth penalties differ
from exact Lasso and, with zero ridge, are not uniformly strongly convex
as required by that theorem.

The Gaussian-mixture curves use
\citet[equation~(8)]{LoureiroMixtures2021}. Fit $K=1$ or $10$ full-covariance Gaussian components
per label to one frozen noisy realization of the 64-dimensional core,
using the reference class frequencies and three EM starts.

\section{Support recovery and fresh-score comparisons}
\label{app:support-simulation}
\label{sec:support-recovery}

\begin{figure}[!htb]
 \centering
 \includegraphics[width=\linewidth]{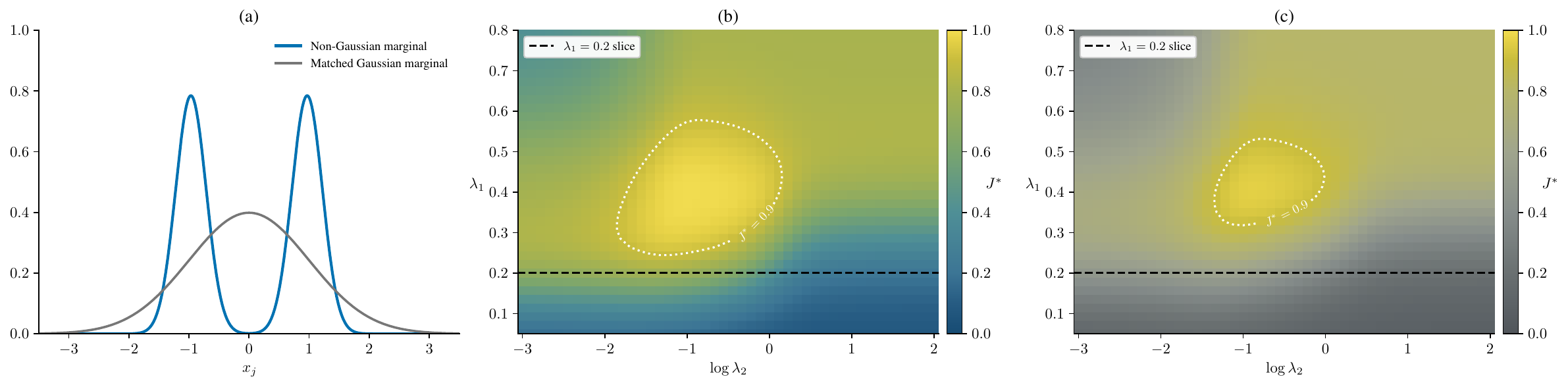}
 \caption{Feature marginals and predicted support recovery.
 (a) Exact block-coordinate marginal densities for the non-Gaussian
 design (blue solid) and matched Gaussian design (gray solid);
 both have mean zero and variance one, and the full designs share
 the same covariance. (b,c) Surrogate-predicted Jaccard $J^*$ on the $36\times46$
penalty grid for the non-Gaussian (blue--yellow) and matched Gaussian
(gray--yellow) designs, respectively. White dotted contours mark
$J^*=0.9$, and horizontal
dashed lines mark the $\lambda_1=0.2$ slice.}
 \label{fig:support-island}
\end{figure}

\paragraph{Setting and metrics.}
We compare support recovery and fresh-score laws using the synthetic
model of Appendix~\ref{app:synthetic-prediction-experiments}
with $k=1$, hence $(n,p)=(200,400)$, under both the non-Gaussian and
covariance-matched Gaussian designs.
Figure~\ref{fig:support-island}(a) compares a representative block-coordinate
marginal with its matched Gaussian counterpart; both have mean zero and
variance one. Let $S_0=\operatorname{supp}(\theta_0)$,
with $S_{0,b}$ and $S_{0,i}$ its 20 active block coordinates and five
active singleton coordinates. For a selected support $S$, define
\[
 J(S,S_0)=\frac{|S\cap S_0|}{|S\cup S_0|},\qquad
 \mathrm{TPR}_{\rm block}(S)=\frac{|S\cap S_{0,b}|}{20},\qquad
 \mathrm{TPR}_{\rm singleton}(S)=\frac{|S\cap S_{0,i}|}{5}.
\]
The surrogate prediction is
\begin{equation}
\label{eq:J*star}
 J^*=\E_G\left[J(\operatorname{supp}(W_\ast(G)),S_0)\right],
\end{equation}
with analogous expectations for the two true-positive rates. Numerical
supports use the threshold $|W_j|/\sqrt n>10^{-8}$ for both ERM and the
surrogate. Support recovery is examined numerically here; it is not a
consequence of Theorem~\ref{thm:endpoint}.

\paragraph{Penalty grid and replication.}
We calibrate \eqref{eq:losscoeff}--\eqref{eq:fixedpoint} separately under
each feature--response law on the $36\times46$ grid
\[
 \lambda_1\in\{0.05+0.75k/35:k=0,\ldots,35\},\qquad
 \lambda_2\in\{10^{-3+k/9}:k=0,\ldots,45\}.
\]
Empirical support metrics are averaged over 500 independent training
datasets per design, reusing each dataset across the grid. The
$\lambda_1=0.2$ slice uses the same 46 ridge values and training fits.

\begin{figure}[!htb]
 \centering
 \includegraphics[width=\linewidth]{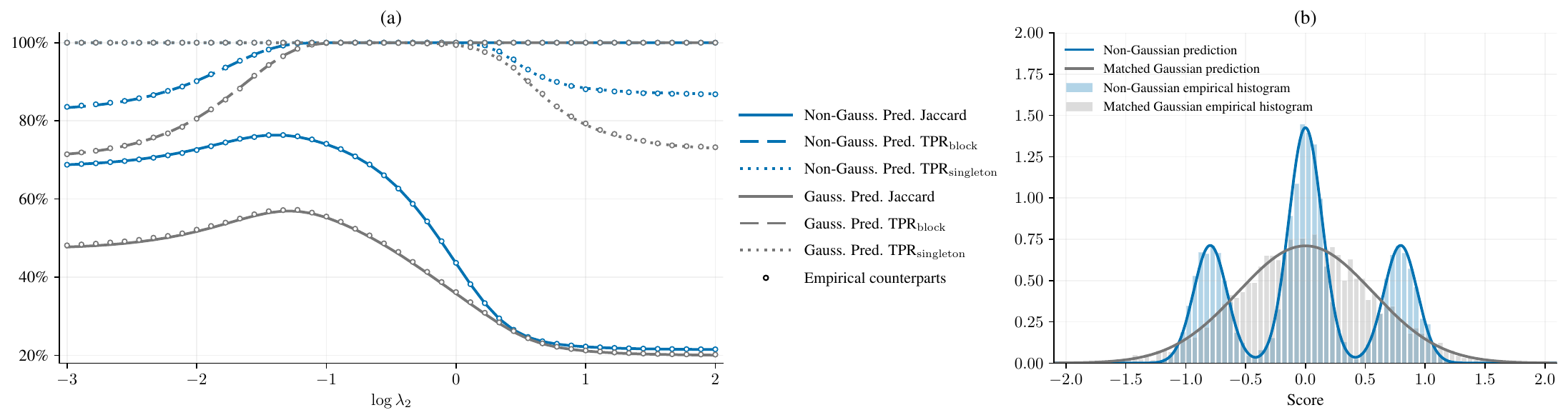}
 \caption{Support recovery and fresh-score distributions for the non-Gaussian
 (blue) and covariance-matched Gaussian (gray) designs.
 (a) Support recovery along the $\lambda_1=0.2$ slice.
 Lines show surrogate expectations: solid for Jaccard, long dashes for
 $\text{TPR}_{\mathrm{block}}$, and dotted for $\text{TPR}_{\mathrm{singleton}}$. ``Pred.'' denotes surrogate
 predictions; ``Gauss.'' denotes the covariance-matched Gaussian design.
 Open circles show empirical means of the corresponding metrics over independent training datasets per design.
 (b) Fresh-score densities at $(\lambda_1,\lambda_2)=(0.2,0.1)$.
 Solid curves show the predicted laws of $x^\top\mu_\ast+\alpha_\ast Z$
 under each design's own calibration, where $Z\sim N(0,1)$ is independent
 of $x$. The blue and gray histograms show empirical scores
 $x^\top\widehat W/\sqrt n$ for the non-Gaussian and matched Gaussian
 designs. Writing
 $\mu_\ast=m_{b*}\mathbf{1}_{S_{0,b}}+m_{i*} \mathbf{1}_{S_{0,i}}$,
 where $\mathbf{1}_S$ is the coordinate-indicator vector of $S$,
 the calibrated values $(m_{b*},m_{i*},\alpha_\ast)$ are approximately
 $(0.04111,0.04924,0.02859)$ for the non-Gaussian design and
 $(0.03961,0.05192,0.04698)$ for the matched Gaussian design.}
 \label{fig:support-slice}
\end{figure}

\paragraph{Results.}
\label{app:support-slice}
Figure~\ref{fig:support-island}(b,c) shows the dependence of predicted recovery
on both penalty coefficients. The maximum absolute difference between
$J^*$ and the empirical mean Jaccard over the grid is $0.74$ percentage
points for the non-Gaussian design and $0.72$ for the matched Gaussian
design. Along the $\lambda_1=0.2$ slice in
Figure~\ref{fig:support-slice}(a), increasing $\lambda_2$ initially improves
block-signal recovery and Jaccard similarity. At larger $\lambda_2$,
increasing false selections and declining singleton recovery reduce
Jaccard despite high block-signal recovery. The two designs yield
different recovery profiles despite their matched feature covariances.

In Figure~\ref{fig:support-slice}(b), the empirical histograms for both
designs closely follow the prediction laws of
Corollary~\ref{cor:conditional-prediction}. Although the designs have the same feature covariance,
their separately calibrated score laws differ, with an unstandardized
$\Wass$ distance of approximately $0.10994$.
\FloatBarrier

\section{Simple criteria for a Poincar\'e inequality}
\label{app:poincare-criteria}

A random vector $U\in\R^d$ satisfies a Poincar\'e inequality with
constant $C$ if
\begin{equation}
\label{eq:appendix-poincare-definition}
 \Var(f(U))\le C\E\left[|\nabla f(U)|^2\right]
\end{equation}
for every bounded $C^1$ function $f$ with bounded gradient.
The inequality extends by smooth truncation to bounded $C^1$ tests
with finite gradient energy. All norms below are Euclidean.

\subsection{Tensorization and Lipschitz transformations}

\begin{lemma}[Stability properties]
\label{lem:poincare-stability}
Let $U_1,\ldots,U_m$ be independent random vectors satisfying
\eqref{eq:appendix-poincare-definition} with constants $C_1,\ldots,C_m$.
Then
\begin{equation}
\label{eq:poincare-weighted-tensorization}
 \Var(F(U_1,\ldots,U_m))
 \le\sum_{i=1}^m C_i\E\left[|\nabla_iF(U_1,\ldots,U_m)|^2\right].
\end{equation}
In particular, their concatenation has constant $\max_i C_i$.
If $U$ has constant $C$ and $\Phi:\R^d\to\R^p$ is
$L$-Lipschitz, then $\Phi(U)$ has constant $CL^2$.
\end{lemma}

\begin{proof}
The conditional-variance tensorization inequality, followed by the
Poincar\'e inequality in each coordinate block, gives
\eqref{eq:poincare-weighted-tensorization}.
For smooth $\Phi$, the chain rule gives
\[
 \Var(f(\Phi(U)))
 \le C\E\left[|D\Phi(U)^\top\nabla f(\Phi(U))|^2\right]
 \le CL^2\E\left[|\nabla f(\Phi(U))|^2\right].
\]
For general Lipschitz $\Phi$, mollify it without increasing its
Lipschitz constant and pass to the limit by dominated convergence.
This also covers singular laws of $U$.
For tensorization, see \citet[Section~4.3]{BakryGentilLedoux2014}.
\end{proof}

\subsection{Generated data}

A standard Gaussian vector has Poincar\'e constant $1$.
Consequently, for a fixed generator $\Phi$,
\[
 Z\sim N(0,I_d),\qquad x=\Phi(Z),\qquad
 \operatorname{Lip}(\Phi)\le L
 \quad\Longrightarrow\quad
 \Var(f(x))\le L^2\E\left[|\nabla f(x)|^2\right].
\]
This allows non-Gaussian features with dependent coordinates.
The bound is dimension-independent when $L$ is uniformly bounded.
For a feedforward network with $1$-Lipschitz activations and weight
matrices $A_\ell$, it suffices that
$\prod_\ell\|A_\ell\|_{\op}$ is uniformly bounded.
The same criterion applies to the joint pair $(x,y)$ when that pair
is a Lipschitz function of a Poincar\'e latent vector.

\citet[Section~2 and Section~3.1]{Seddik2020GAN} use Lipschitz Gaussian generators
to model GAN-generated data and their representations, with
Gaussian-tailed concentration. Indeed, the above construction gives,
for every $1$-Lipschitz $f$ and $L>0$,
\[
 \Pp\{|f(x)-\E_x\left[f(x)\right]|\ge t\}
 \le 2\exp\!\left(-\frac{t^2}{2L^2}\right).
\]
These tails are stronger than the exponential tails generally
provided by a Poincar\'e inequality. However, a tail bound alone
does not imply the gradient inequality
\eqref{eq:appendix-poincare-definition}; here the latter follows
directly from Lemma~\ref{lem:poincare-stability}.

\end{document}